\documentclass{amsart}
\usepackage{amssymb}
\usepackage[margin=1.4in]{geometry}
\usepackage{graphicx}
\usepackage{hyperref}
\usepackage{mathrsfs}
\usepackage{enumitem}
\usepackage{stmaryrd}
\usepackage{hyperref}
\usepackage{mathtools}
\usepackage{microtype}

\newtheorem{theorem}{Theorem}[section]

\theoremstyle{definition}

\theoremstyle{remark}

\numberwithin{equation}{section}

\newcommand{\jump}[1]{\left\llbracket{#1}\right\rrbracket}

\newcommand{\placeholder}{\,\cdot\,}

\begin{document}

 \author[Haziot]{Susanna~V.~Haziot}
 \address{Fakult\"at f\"ur Mathematik, Universit\"at Wien, Oskar-Morgenstern-Platz 1, 1090, Wien, Austria}
 \email{susanna.haziot@univie.ac.at}

 \author[Hur]{Vera~Mikyoung~Hur}
 \address{Department of Mathematics, University of Illinois at Urbana-Champaign, Urbana, IL 61801, USA}
 \email{verahur@math.uiuc.edu}

\author[Strauss]{Walter Strauss}
\address{Department of Mathematics and the Lefschetz Center for Dynamical Systems, Brown University, Box 1917, Providence, RI 02912, USA}
\email{walter\_strauss@brown.edu}

\author[Toland]{J.~F.~Toland}
\address{Department of Mathematical Sciences, University of Bath, Bath,
 BA2 7AY, UK}
 \email{masjft@bath.ac.uk}

 \author[Wahl\'en]{Erik Wahl\'en}
 \address{Centre for Mathematical Sciences, Lund University, PO Box 118, 22100 Lund, Sweden}
 \email{erik.wahlen@math.lu.se}

 \author[Walsh]{Samuel Walsh}
 \address{Department of Mathematics, University of Missouri, 202 Math Sciences Building,
Columbia, MO 65211, USA}
 \email{walshsa@missouri.edu}

 \author[Wheeler]{Miles H.~Wheeler}
 \address{Department of Mathematical Sciences, University of Bath, Bath, BA2 7AY, UK}
 \email{mw2319@bath.ac.uk}

 \title{Traveling water waves --- the ebb and flow of two centuries}

\begin{abstract}
This survey covers 
the mathematical theory of steady water waves with an emphasis on topics that are at the forefront of current research.  These areas include: variational characterizations of traveling water waves; analytical and numerical studies of periodic waves with critical layers that may overhang; existence, nonexistence, and qualitative theory of solitary waves and fronts; traveling waves with localized vorticity or density stratification; and waves in three dimensions.
\end{abstract}

\maketitle

\tableofcontents

\section{Introduction (by W. Strauss)}\label{sec:walter}

\subsection{Fluids}
 The motion of a fluid can be very complicated, as we know whenever we see waves break on a beach,
fly in an airplane, or look at a lake on a windy day.
The first comprehensive mathematical model of a fluid was proposed by Euler in the 1750's.
It is amazing that Euler's equations and its variants are still the basic models
used by mathematicians and engineers today!
The most notable variant was introduced by Navier and Stokes to allow for the fluid to be viscous.
Of course an immense amount of progress has been made over the past two and a half centuries.
But one thing we know is that the difficulties of understanding fluids are so profound
that there are huge problems that are still beyond our reach.
This is well-known both to mathematicians and to more practical scientists like hydraulic
and aeronautical engineers.

Euler  wrote
\begin{quote}
{\it``If it is not permitted to us to penetrate to a complete knowledge concerning the motion
of fluids, it is not to mechanics nor to the insufficiency of the known principles of motion that
we must attribute the cause.  It is analysis itself which abandons us here since all the theory
of the motion of fluids has been reduced to the solution of analytic formulas."}
\end{quote}
Euler's comment was amazingly prescient!  In the many years since his time,
the theory of fluids has had a tremendous effect on the historical development of mathematics itself.
For instance, it was the main motivation for Cauchy to develop his complex function theory.

{The Euler Equations are as follows. }
Denote
${\bf u}({\bf x}, t)$ as the velocity of the water at the spatial point ${\bf x}\in\mathbb{R}^3$ and time $t$.
Denote $\rho({\bf x},t)$ as the water's density.
Now water flows are essentially incompressible.
Conservation of mass is expressed by
\begin{equation} \label{masscons}
\frac{\partial \rho}{\partial t} + \nabla\cdot (\rho\ {\bf u}) =0.
\end{equation}
Under the assumption that the density $\rho$ is constant, it simply reduces to
$\nabla\cdot {\bf u}=0$.
Conservation of momentum (Newton's Law of Motion) is expressed by
\begin{equation}\label{Euler}
\frac{\partial{\bf u}}{\partial t} + ({\bf  u}\cdot\nabla){\bf  u} + \frac1\rho \nabla P =    {\bf F}.
\end{equation}
The first two terms express the acceleration because if we write a fluid particle
path as ${\bf x}={\bf  X}(t)$, then its velocity is ${\bf u}({\bf X}(t),t)$.
Its acceleration is
\begin{equation}  \label{accel}
\frac{d}{dt} {\bf u}({\bf  X}(t),t) =    \nabla_x{\bf u} \cdot  \frac{\partial{\bf  X}}{\partial t} + \frac{\partial {\bf u}}{\partial t}
= \frac{\partial{\bf u}}{\partial t} + ({\bf u}\cdot\nabla){\bf u}.
\end{equation}
$P({\bf x},t)$ is the pressure, and $\rho{\bf  F}$ is any outside force that may be acting on the fluid
(for instance, gravity).   Thus there are 4 scalar equations in 4 unknowns.

\subsection{Water waves}
Now consider a water wave.
The water lies below a body of air.
The interface between the air and the water is a {\it free surface} $S\subset \mathbb{R}^3$.
Gravity points down.

In addition to water waves, there are many other interesting free surface problems, such as
(1) the melting of ice, the free surface being the boundary of the ice,  and
(2) the stretching of a flexible membrane over an obstacle,
the free surface being the boundary of the contact region.

Now consider a water wave below $S$ and above a flat bottom $B$.
We'll assume the air is quiescent and
$S$ is smooth, with no bubbles or jets coming in or going out of the water.
Because $S$ is an unknown, there has to be an extra boundary condition on it.
The  standard boundary conditions are:

BC1: Each water particle on $S$ remains on $S$ for all time.
That is, particles of water do not invade the air, nor vice-versa.
This implies that there are no jets or bubbles.

BC2: On the surface the pressure $P$ equals the
atmospheric pressure $P_{atm}$ of the air.   (We are assuming zero surface tension.)

BC3: Each water particle on the bottom $B$ remains on $B$ for all time.
That is, the bottom is impermeable.

A variant of BC3 is the case of infinite depth, in which case the ``boundary" condition is that
the water becomes quiescent as one goes to deeper depths.  Think of the middle of the ocean.

\begin{figure}
  \centering
  \includegraphics
  [width=.7\linewidth  ]   
  {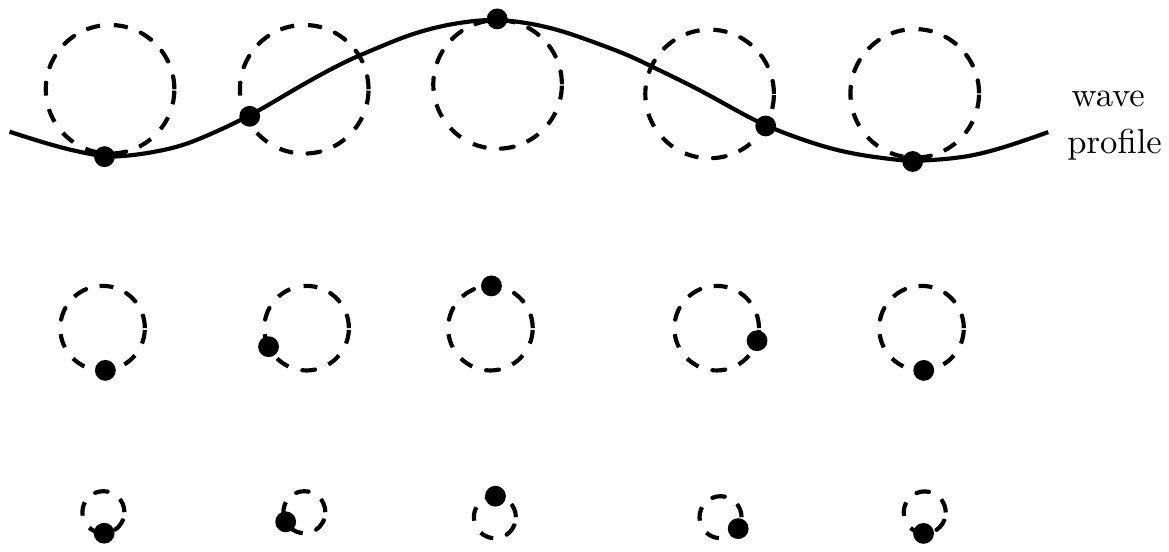}
  \caption{Trajectories in the Gerstner Wave}
  \label{walter:gerstner figure}
\end{figure}

The earliest history of water waves began with
Laplace (1776) and Lagrange (1781), who linearized the equations around quiescence.
Lagrange used coordinates that follow the individual particles, now known as the Lagrangian formulation.
Remarkably, Gerstner (1802) found an {\it exact} family of traveling solutions
of the full nonlinear equations (see Figure~\ref{walter:gerstner figure}).   The Gerstner
wave profile $S$ moves horizontally but every particle moves in a perfect circle.

Early in the 19th century Fourier developed his theory of heat conduction.
Around 1815 Cauchy and Poisson surely knew about Fourier's expansions
and they used them to study water waves near a traveling wave, linearized around quiescence.
Then Cauchy's 2D water wave analysis was a prime motivation to his development of complex analysis (1825).

\begin{figure}
  \centering
  \includegraphics
  [scale=1]   
  {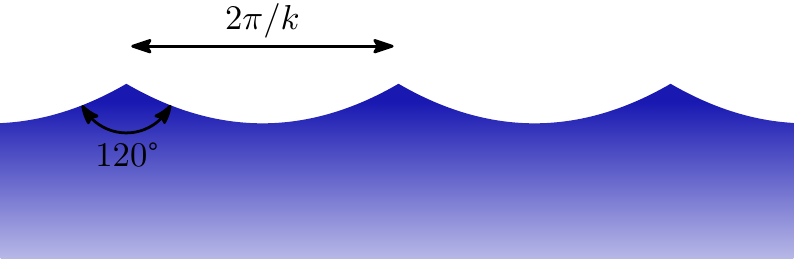}
  \caption{Stokes' Extreme Wave}
    \label{walter:extreme wave figure}
\end{figure}

An important figure was Airy, who in particular found the dispersion relation $c^2k=g\tanh(kd)$ around 1845.
Here $c$ is the speed of the traveling wave, $d$ is the average depth, and $k$ is the wave number
of a wave that is periodic in $x$.
Stokes had a long career starting in the 1840's.  Among many accomplishments,
he computed the expansion of the waves in a power series of the amplitude, discussed  waves
of large amplitude, and discovered the extreme wave of greatest height that must have a $120^\circ$ angle at its crest
(see Figure~\ref{walter:extreme wave figure}).

\subsection{Traveling water waves}

Traveling means that ${\bf u} = {\bf u}(x-ct, y, z),\ P = P(x-ct, y, z)$.
Two-dimensional (2D) means they are independent of one horizontal coordinate.
In the 2D case we'll keep $x$ and make $y$ the vertical variable.
In the 3D case we would take $z$ as the vertical variable.

A large part of this article, as well as most of the literature, is concerned with the 2D traveling WWs.
The fundamental reason is that the equation for the stream function
(see \eqref{stream eq} below) is elliptic
in 2D but not in 3D.   So in 2D we denote ${\bf x}=(x,y),\ {\bf u} = (u,v)$.
The force of gravity is ${\bf F}=(0,-g)$.  We also take $\rho=1$.
If we write $x$ instead of $x-ct$ for a traveling wave,
then $\frac{\partial}{\partial t}$ becomes $-c\frac{\partial}{\partial x}$.  The Euler equations

\begin{equation}  \label{Euler1}
\nabla\cdot {\bf u}=0,\qquad \frac{\partial{\bf u}}{\partial t} + ({\bf u}\cdot\nabla){\bf  u} +  \nabla P =    {\bf F}.
\end{equation}
 become
\begin{equation}   \label{Euler2D}
\begin{array}  {l}
u_x +v_y= 0,   \\
         (u-c)u_x + vu_y  =-P_x,  \quad    (u-c)v_x + vv_y=-P_y-g.
\end{array}   \end{equation}
The
boundary conditions become
\begin{equation}   \label{WWBC}
 \begin{array}  {l}
(u-c)\,\eta_x  + v\,\eta_y=0  \qquad\text{on}\quad S=\{\eta(x,y)=0\}, \\
P = P_{atm}\qquad \quad  \qquad\qquad\text{on}\quad S, \\
v=0\qquad\qquad\qquad\qquad  \  \text{ on}\quad B=\{y=-d\}.
\end{array}  \end{equation}
Here $P_{atm}=$ atmospheric pressure, $g=$
gravitational acceleration (vertical), $d=$  average depth of the water.
It is assumed here that the only force on the fluid particles is gravity (no viscosity).

The {\it vorticity} is $\omega= v_x-u_y$.\ \   (In 3D it is $\nabla\times {\bf u}$.)
In 2D the relative {\it stream function $\psi$} is determined up to a constant by
\[\psi_y=u-c,\quad \psi_x=-v. \]
In terms of $\psi$ the whole problem
\eqref{Euler2D} - \eqref{WWBC}
 nicely reduces to
\begin{equation}\label{stream eq}
\left\{\begin{array}{l}
\Delta \psi=-\omega \quad\hbox{in}\quad \Omega,\\[0.15cm]
\psi=0\quad\hbox{on}\quad S,\\[0.15cm]
\psi=-m\quad\hbox{on}\quad B,\\[0.15cm]
|\nabla\psi|^2 =Q-2gy \quad\hbox{on}\quad S,
\end{array}\right.
\end{equation}
where $m$ and $Q$ are constants.
In fact, the first equation in \eqref{stream eq}  is immediate from the definition of $\omega$.
The first equation in \eqref{WWBC} means that $\psi$ is constant on $S$, which we normalize to be 0.
From the bottom boundary condition we then see that $m = \int_0^{\eta(x)} (u(x,y)-c) dy$
is the (constant) relative mass flux.
The last equation in \eqref{stream eq} is the Bernoulli relation, which is verified by direct differentiation.
Clearly  $2gy\le Q$, which provides the greatest possible height.
$Q$ is the ``total head" or ``Bernoulli constant".

A {\it stagnation point} is a point where $u=c,\ v=0$, which means that $\nabla\psi=0$.
More generally, a {\it critical point} is a point where $u=c$, which means merely
that $\psi_y=0$.
By eliminating the pressure in \eqref{Euler2D},
we see that $(u-c,v)$ is orthogonal to $(\omega_x,\omega_y)$.
So,  away from stagnation points, $\omega$ and $\psi$ are functionally dependent.
Writing $\omega(x,y) = \gamma(\psi(x,y))$, where we assume that $\gamma$ is single-valued,
we get the nonlinear elliptic PDE
\begin{equation}   \label{vorticity eq}
  -\Delta \psi  =   \gamma(\psi)
\end{equation}
  in the fluid domain.
  The unknowns in the system \ref{stream eq} - \ref{vorticity eq}
  are $\psi$ and $S$ and the two constants $m$ and $Q$.

 In the {\it irrotational} case, an alternative formulation can be made
 in terms of the {\it velocity potential} $\phi$.
In 3D the vorticity is $\nabla\times\mathbf{u}$.
If the vorticity vanishes, then of course $\mathbf{ u}=\nabla\phi$ for some scalar $\phi$.
Due to $\nabla\cdot\mathbf{u}=0$, $\phi$ is harmonic.
Specializing to 2D in terms of $\phi$, the traveling irrotational water wave system is
\begin{equation}    \label{irrotWW}
\left\{\begin{array}{l}
\Delta \phi = 0     \quad\hbox{in}\quad \Omega,    \\[0.15cm]

\eta_x (\phi_x -c)  -\phi_y = 0 \quad \text{ on } S,   \\[0.15cm]

-c\phi + \tfrac12 | \nabla\phi |^2 = - P_{atm}  -gy \quad \text{ on } S,   \\[0.15cm]

\phi_y = 0 \quad \text{ on } B.
 \end{array}\right.
\end{equation}

\subsection{The last century}

Most work continued to be both 2D and irrotational.
In that case $\psi$ is harmonic in the plane fluid region.
So, using the Green's function of the Laplacian,
the problem can be formulated as an equation on the curve $S$,
the unknowns being functions of the single variable $x$.
Nekrasov \cite{Nekrasov} did this and was the first to construct many periodic exact solutions
of small amplitude $a$ for infinite depth.  He expanded a presumed solution in powers of $a$
and proved that there is a positive radius of convergence.
Later, this was done independently and more directly without the Green's function by Levi-Civita \cite{lc},
who expanded everything more directly in powers of $a$.
Struik \cite{struik} extended this analysis to the case of finite depth.
Dubreil-Jacotin \cite{dubreil} extended Levi-Civita's approach to the case with vorticity.
All of these rigorously constructed waves had small amplitudes.

The first large-amplitude waves (periodic, irrotational) were constructed by
Krasovskii \cite{krasovskii} by making use of Nekrasov's equation
together with Krasnoselskii's positive operator theory.
Later, Keady and Norbury \cite{keady} employed Dancer's elegant global bifurcation theory to obtain a
smooth curve of solutions of large amplitude.
Amick, Fraenkel and Toland \cite{AFT1982} took the limit of this curve to construct the Stokes  extreme wave.

Due to the great difficulty of the water wave problem,
right from the earliest years many approximate models of water waves have been studied.
They are formally obtained by various scaling limits.
The most famous one is the Korteweg-deVries equation (KdV)
\begin{equation}   \label{KdV}
 \alpha_\tau + \alpha_{xxx} + \alpha\alpha_x = 0,
\end{equation}
first derived by Boussinesq \cite{bouss},
It comes formally from the 2D water wave system by a rescaling in time $\tau=\epsilon t$,
an expansion in the small amplitude $a$ and an assumption that the wave is unidirectional.
The limit is rigorously valid in a region $\{|x|\le O(\epsilon^{-3/2}),\ t\le O(\epsilon^{-3/2})\}$;
see Schneider and Wayne \cite{schn-wayne}.
KdV also occurs as an approximate model in many other physical problems
in situations where there is an appropriate balance between dispersion and nonlinearity,
such as in plasma theory.
KdV  has many remarkable properties, notably the existence of
solitons that are extremely stable traveling waves.
There are many other approximate models of water waves that have been extensively studied, such as
the nonlinear Schr\"odinger equation, several Boussinesq systems, the Camassa-Holm equation and
the Kadomtsev-Petviashvili equation,.   Several of these are ``completely integrable".

\subsection {Recent and future research areas}
 \begin{itemize}
\item Vorticity.  Waves of large amplitude can be  obtained analytically either by
Leray-Schauder degree or by the Dancer-Buffoni-Toland theory of analytic families of operators.
Large periodic waves with vorticity were first constructed in \cite{CS2004}.
\item Besides traveling waves that are periodic, there are those that are solitary or are fronts,
which means that the solution becomes flat at spatial infinity.  A solitary wave has the same height
at both $x\to\infty$ and $x\to-\infty$, while a front has two different heights.  See Section \ref{sec:miles} below.

\item Various aspects of the shape of the surface $S$ have been studied, including
symmetry properties, locations of the maximum and minimum heights,
the maximum steepness of a wave, and the existence of overhanging waves.

\item Regularity of the waves.  It is generally true that every streamline (including $S$) is
an analytic curve, as in \cite{const-escher}.

\item The particle paths (streamlines) are usually not like the Gerstner wave.
Typical ones have drift and can look like loops, as in \cite{CS2010}.

\item Stagnation points and critical layers occur in some traveling waves
if there is vorticity or if the density is variable. See Sections \ref{sec:susanna},
 \ref{sec:sam:localized} and \ref{sec:sam:strat} below.

\item Physical properties of the fluid can have major effects on the nature of the traveling waves.
Examples are surface tension, density, temperature, salinity and viscosity,

\item Internal waves occur in addition to the surface wave in case the density (or temperature, etc.)
has a jump inside the fluid, as in a two-fluid model.

\item Of course the bottom $B$ might not be flat, which would prevent a normal traveling wave.
See e.g.  \cite{CGS}.

\item Traveling waves in 3D are much more subtle.  See
\cite{IoossPlotnikov11} and Section \ref{sec:erik} below.
\end{itemize}

Furthermore,  there are problems that are fundamentally time-dependent.

\begin{itemize}
\item Stablity and instability.  If one considers only spectral stability, then time disappears from the problem.
See \cite{Bridges} and \cite{NgSt}.

\item The initial value problem has been studied intensively, especially since
\cite{Wu97}.
The goal is to solve the water wave system either locally in time, or globally in time with
small initial data.

\item Floating bodies such as a ship.  See \cite{LannesFloating}.

\item Mathematically understanding the nature of turbulence is a major challenge for the 21st century.

\end{itemize}

Thus we see that there are large areas to investigate and a huge number of open problems!

Section \ref{sec:john} is devoted to some variational aspects of irrotational water waves, Section \ref{sec:susanna} is devoted to rotational water waves in stratified media, Section \ref{sec:vera} is devoted to constant vorticity flows, Section \ref{sec:miles} is devoted to solitary waves and fronts, Section \ref{sec:sam:localized} is devoted to localized vorticity, Section \ref{sec:sam:strat} is devoted to stratification and Section \ref{sec:erik} is devoted to 3D waves.

\section{Variational aspects of irrotational water wave theory (by J.~F.~Toland)}\label{sec:john}

\subsection{The  mathematical problem}
The water-wave problem is one of determining the  profile of the unconstrained surface of a heavy liquid  in motion under the effect of  gravity $g$ acting vertically down.   The simplest model assumes that viscosity and surface tension can be neglected,  that the fluid depth is infinite, and that the flow is 2-dimensional and  irrotational. This  means that  at time $t$ the fluid occupies a region
$$\widetilde\Omega_t =\{ (x,y,z):  y < \tilde\eta(x,t)\} \text{  under the surface } \widetilde{S}_t = \{(x,\tilde\eta(x,t),z): x,\,z \in \mathbb{R}\},
$$
where $(x,y,z)$ denotes
  orthogonal Eulerian  coordinates, with $y$ in  the vertical direction and the fluid motion independent of the $z$-coordinate.
Therefore  the velocity at  a point $(x,y,z)$ at time $t$ is the gradient of a scalar velocity potential $\tilde \phi$
$$
\vec v(x,y;t) = \nabla \tilde \phi (x, y;t) \text{ where } \nabla = (\partial_x,\partial_y) ~~ \text{and $\tilde \phi$ is independent of $z$.}
$$
When for convenience the flow is normalised to be $2\pi$-periodic in $x$, $\tilde \phi$ satisfies
\begin{subequations}\label{fens}
\begin{align}
 &\Delta \tilde \phi (x,y;t) = 0 \text{ on } \widetilde\Omega_t;
\\& \nabla \tilde \phi(x,y;t) \text{ is $2\pi$-periodic in $x$};
\\&\nabla \tilde \phi
(x,y;t) \rightarrow (0,0) \text{ as } y \rightarrow -\infty;
\\&\tilde \eta_t(x,t) + \tilde \phi_x(x,\eta(x,t),t)\tilde \eta_x(x,t) - \tilde \phi_y(x,\eta(x,t),t) = 0;
\\&\tilde \phi_t(x,\eta(x,t),t) +
{\scriptstyle{\frac12}} |\nabla \tilde \phi(x,\eta(x,t),t|^2 +g \tilde \eta(x,t) = 0.\label{11e}
\end{align}\end{subequations}
\textbf{Remark.} If $\tilde \phi$ satisfies \eqref{fens}, except that instead of zero on the right in \eqref{11e} there is a time dependent function $c(t)$. Then $\tilde \phi -\int_0^t c(s)ds$ satisfies \eqref{fens}.\qed

\textbf{Stokes Waves } (steady irrotational water waves) are solutions of \eqref{fens} that represent waves that travel, without change of shape and with constant horizontal velocity $c$ in the $x$-direction, on the surface of fluid that is at rest at infinite. Since the profile of the wave  in a vertical plane parallel to the $x$-direction is  the unknown, suppose that  $\tilde\eta(x,t) = \eta(x-ct),\,c\neq 0$, so that $\widetilde \Omega_t = \Omega + (ct,0)$,
$\widetilde {S}_t = S + (ct,0)$ where $\Omega=\{(x,y):y< \eta(x)\}$ and  $S=\{(x,\eta(x)): x \in \mathbb{R}\}$.
\begin{subequations}\label{pens}Then
$\tilde\phi(x,y;t) := c\big(\phi(x-ct, y)+(x-ct)\big)$ is a solution of
\eqref{fens} if
\begin{align}
 &\Delta  \phi (x,y) = 0 \text{ on } \Omega;
\\& \nabla  \phi(x,y) \text{ is $2\pi$-periodic in $x$};\\&\nabla  \phi
(x,y) \rightarrow (-1,0) \text{ as } y \rightarrow -\infty;
\\&  \eta'(x)\partial_x \phi(x,\eta(x)) -\partial_y\phi(x,\eta(x))  = 0;
\\& {\scriptstyle{\frac{c^2}2}} |\nabla  \phi(x,\eta(x))|^2 +g \eta(x) = {\scriptstyle{\frac{c^2}2}}.
\end{align}
\end{subequations}
Thus Stokes-waves may be regarded as steady, since \eqref{pens}, which    is time-independent, describes the flow  relative to a frame moving with the wave velocity.
\begin{subequations}\label{qens}
 It follows from \eqref{pens} and the Cauchy-Riemann equations that a harmonic conjugate $\psi$ of $\phi$ (that is, the relative stream function) satisfies
\begin{align}\label{qa}
 &\Delta  \psi (x,y) = 0 \text{ on } \Omega;
\\& \nabla  \psi(x,y) \text{ is $2\pi$-periodic in $x$}; \label{qb}
\\&\nabla  \psi(x,y) \rightarrow (0,-1) \text{ as } y \rightarrow -\infty; \label{qc}
\\&  \eta'(x)\partial_y \psi(x,\eta(x)) +\partial_x\psi(x,\eta(x)  = 0;\label{qd}
\\& {\scriptstyle{\frac{c^2}2}} |\nabla  \psi(x,\eta(x))|^2 +g \eta(x) = {\scriptstyle{\frac{c^2}2}}.\label{qe}
\\
\intertext{Since \eqref{qd} implies that $\psi$ is  constant on $S$, there is no loss in assuming that} &\psi \equiv 0 \text{ on $S$, and hence $\psi >0$ on $\Omega$}\label{qf}
\end{align}
\end{subequations}
by \eqref{qc} and the maximum principle. By \eqref{qe}, the pressure in the flow at the surface is constant (atmospheric pressure).
Let $Q$ be the speed of the steady flow at  the wave crest, $(0, \eta(0))$ say, in a frame  moving with the wave speed $c$:
\begin{equation} \label{Q} Q = |\nabla  \phi(0,\eta(0))| = |\nabla  \psi(0,\eta(0))|, \text{ and let }  \nu = \frac{Q^3}{3g c}.
\end{equation} \qed

\subsection{Nekrasov's Equation} If the periodic steady wave is symmetric about the crest, the flow is horizontal at $x = \pm \pi$,  $\phi(\pi, y)- \phi(-\pi, y) = 2\pi$, $y \leq \eta(\pm\pi)$. By \eqref{qens},  $\exp(-\psi + i\phi)$ maps one period of the flow domain  $\Omega$ conformally
onto the unit disc in $\mathbb C$, and   one period of the surface  $S$ onto the unit circle.
Let $\theta(s)$, where $s = \phi(x,\eta(x))$, denote the slope of the surface at the point $(x, \eta(x)) \in S$.
Then  Nekrasov showed that $\nu$ in   \eqref{Q} and
the odd function $\theta:[-\pi,\pi] \to \mathbb{R}$ satisfies the integral equation
\begin{gather}\label{nek}
\theta(s) = \frac{1}{3\pi}\int_{-\pi } ^\pi      \left( \sum_{k=1}^\infty
\frac{ \sin k s \sin k t } {k} \right) \frac{\sin \theta (t)}{\nu + \int_0 ^t
\sin \theta (\nu) d \nu },~  \theta (\pm \pi)=0,
\end{gather}
and from it obtained the first proof
 of
existence \cite{kuznetsov,Nekrasov,stokes} of non-trivial  steady periodic  travelling waves.  Although the equation is valid for waves of all amplitudes, his proof is for small amplitude waves bifurcating from $\theta = 0$ at  $\nu = 1$.
Subsequently, the  geometer Tullio Levi-Civita  found a slightly  different proof of the same result.
Then in 1961, Krasovskii \cite{krasovskii}, using  subtle  estimates in harmonic analysis with Leray-Schauder degree theory,
showed that
for all $\alpha \in [0, \pi/6)$, there exists a smooth solution $\theta$ of \eqref{nek} with $\|\theta\|_\infty = \alpha$, and heuristic arguments from hydrodynamics and  harmonic analysis led him to conjecture  that the constraint  $\alpha \in [0,\pi/6)$ is sharp.
However, although
 McLeod \cite{McLeod}  showed that
 $\|\theta\|_\infty  > \pi/6$ does happen,
  Amick \cite{amick} showed that always $\|\theta\|_\infty  < (1.098)(\pi/6)$. \qed

Since then, familiar  tools such as
the contraction mapping and  implicit function theorems, bifurcation theory,
 Leray-Schauder degree theory,
real-analytic function theory,  and Nash-Moser theory (for standing waves, see the closing remark)  have been applied successfully in  nonlinear water-wave theory.
However variational methods, saddle-point and mini-max principles, the mountain-pass lemma, Morse index theory,
Lyusternik-Schnirelman category etc have been less successful.
Although \eqref{fens}, \eqref{pens} and \eqref{nek} do not have obvious variational structure, \eqref{pens} is a conservative
system for which there should
be a formulation in terms of  Hamiltonians or Lagrangians. How this can be done was first realised by Zakharov \cite{zak}.
\qed

\subsection{Hamiltonian theory of 2D-irrotational water waves.  }

The wave energy,  kinetic + potential, in one period at time $t$, of a
(not necessarily steady) solution of \eqref{fens} is,
\begin{equation}\label{eny} {\frac{1}{2}} \int_0^{2\pi}\int_{-\infty}^{\tilde\eta(x,t)} |\nabla \tilde\phi(x,y;t)|^2|dy
dx + \frac{g}{2}\int_0^{2\pi}\tilde\eta^2(x;t) dx.\end{equation}
Now for any pair of
${2\pi}$-periodic functions $\eta,\Phi:\mathbb{R} \to \mathbb{R}$ let
$$ \Omega =\{(x,y):y<\eta(x)\}, \qquad S = \{(x,\eta(x)): x \in \mathbb{R}\},$$
and  let $\phi$ be the ${2\pi}$-periodic solution of
the  Dirichlet boundary-value problem
\begin{align*}
&\Delta \phi (x,y) = 0 \text{ on $\Omega$}; ~~\phi(x,\eta(x)) = \Phi (x) \text{ on $S$};
\\&\phi \text{ is $2\pi$-periodic in $x$}; ~~\phi(x,y) \rightarrow 0 \text{ as } y \rightarrow -\infty
\end{align*}
and define
 \begin{equation*}\mathcal{E}(\eta,\Phi) : = {\frac{1}{2}}
\int_0^{2\pi}\int_{-\infty}^{\eta(x)} |\nabla \phi(x,y)|^2dy dx +
{\frac {g}{2}\int_0^{2\pi}\eta^2(x) dx}.
\end{equation*}
So formally  $\mathcal{E}$
gives the energy \eqref{eny} of a ${2\pi}$-periodic solution of \eqref{fens} as a functional of $\eta$,  the surface elevation, and $\Phi$, the trace of the velocity potential on  $S$.
Zakharov's insight \cite{zak} was  to consider the infinite-dimensional Hamiltonian system
\begin{equation}\boxed{ \frac{\partial\eta}{\partial t} =
\frac{\partial\mathcal{E}}{\partial \Phi}(\eta,\Phi),\qquad
\frac{\partial\Phi}{\partial t} = -\frac{\partial\mathcal{E}}{\partial
\eta}(\eta,\Phi).}\label{HSB}\end{equation}
This system is in standard canonical form
 \begin{equation}\label{HS}  \dot{x} = J \nabla \mathcal{E} (x),
~~J
=\left(\begin{array}{c}0,I\\-I,0\end{array}\right),\quad x=(\eta,\Phi),\end{equation}
and he showed    how its solutions
yield solutions of     system \eqref{fens}.
However, its rigorous analysis \cite{craigsulem,lannes,nicholls}  is  technically difficult because Dirichlet-Neumann operators on time-dependant domains are involved in the definition  of $\mathcal{E}$.
Nevertheless, as noted in \cite{DKSZ1996}, conformal mapping methods transform   system \eqref{HSB}  into a tidy set of two equations for two real-valued functions,  $\omega$ and $\varphi$, of variables $(\xi,t)$:
\begin{equation}\label{boxed}
\boxed{\begin{aligned}
\dot \chi (1+\mathcal{C} \chi') -\mathcal{C} \varphi' -\chi'\mathcal{C} \dot \omega = 0\\\mathcal{C}\big(\chi'\dot \varphi
-\dot \chi\varphi' +g \omega\chi'\big) +(\dot \varphi+g \chi) (1+\mathcal{C}
\chi')-\varphi'\mathcal{C}
\dot \chi  = 0\\  \varphi(\xi,t) =  \varphi(\xi+2\pi,t),~~
\chi(\xi,t) = \chi(\xi+2\pi,t),\quad \xi,~t \in \mathbb{R}
\end{aligned}}
\end{equation}
Here  $\dot~ = \partial/\partial t$ and $'= \partial/\partial \xi$, and  $\mathcal{C} $ is the Hilbert transform with respect to $\xi \in [-\pi,\pi]$,
\begin{equation}\label{conjugate}
\mathcal{C} \chi (\xi ) =\,pv\frac{1}{2\pi}\int_{-\pi}^{\pi}
\frac{\chi(\sigma)\,d\sigma}{\tan{\scriptstyle{\frac12}} (\xi-\sigma)}.
\end{equation}
According to  \cite[eqns.\,(2.11),\,(2.12)]{{DKSZ1996}} solutions of \eqref{boxed} yield solutions of \eqref{fens} for which the  wave profile is  the periodic curve
\begin{equation}\label{profile} S_t= \{(\xi +\mathcal{C} \chi(\xi,t),\chi(\xi,t)):\xi,\,t \in \mathbb{R}\},
\end{equation}  and $\varphi(\xi,t)$ is the potential on the surface at the point $(\xi+\mathcal{C} \chi(\xi,t),\chi(\xi,t))\in S_t$.
A starting point of this analysis can be the observation that a rectifiable  curve in the plane $\{(x,\eta(x)):x \in \mathbb{R}\}$
where $\eta$ is $2\pi$-periodic
can be parametrized (see  \cite{bernoulli}) as in \eqref{profile}:
$$\{(x,\eta (x)): x \in \mathbb{R} \}= \{ (\xi+\mathcal{C} w(\xi),w(\xi)):\xi \in \mathbb{R}\}, \text{ where $w $ is $2 \pi$ periodic}.
$$

To see that \eqref{boxed} is a Hamiltonian system, although clearly not  in the form \eqref{HS}, let
\begin{multline*}
\varrho_{(\chi,\varphi)}\big((\chi_1,\varphi_1),(\chi_2,\varphi_2)\big) =\int_{-\pi}^{\pi} (1+\mathcal{C}
\chi') ( \varphi_2\chi_1-\varphi_1\chi_2)\\+\chi'\big(\varphi_1\mathcal{C} \chi_2-\varphi_2\mathcal{C} \chi_1\big)
- \varphi'\big(\chi_1\mathcal{C} \chi_2 -\chi_2\mathcal{C} \chi_1\big)\,d\xi.
\end{multline*}
 This bilinear form  is  obviously skew-symmetric, and  exact
because
 $$ \varrho = d\rho
~~\text{ where }~~ \rho_{\varphi,\chi}(\hat \chi,\hat \varphi) = \int_{-\pi}^{\pi}
\big\{\varphi(1+\mathcal{C} \chi') + \mathcal{C}\big(\varphi \chi'\big)\big\}\,\hat \chi d\xi,$$
and it is non-degenerate  by Riemann-Hilbert theory (see \cite{fareast} for details). Hence \eqref{boxed}, which is the Hamiltonian system corresponding  to the skew form $\varrho$ and the  Hamiltonian
$$\mathcal{E} (\chi, \varphi) =
{\frac{1}{2}}\int_{-\pi}^{\pi} \varphi\mathcal{C} \varphi' + g \chi^2(1+\mathcal{C} \chi') \,d\xi,
$$
for  two real functions $(\varphi(\xi,t),\chi(\xi,t))$ of two real variables which are $2\pi$-periodic in $\xi$,
yields solutions of  \eqref{fens}, and in particular leads to the existence of standing-waves \cite{Ioossp1,Ioossp2,ioosspt}, but
not by methods of the calculus of variations (see closing remark).\qed

\subsection {Variational formulation of Stokes waves.}
When $\varphi (\xi,t) = \varPhi(\xi-ct)$ and $\chi(\xi,t) = w(\xi-ct)$, system \eqref{boxed}
 simplifies  because $\varPhi' = c\,\mathcal{C} w'$ by the first equation, whence the second equation involves $w$ only:
\begin{equation}\label{P}\tag{$\mathcal{B}$}\mathcal{C} w'= \lambda\big(w+w\mathcal{C} w' +
\mathcal{C}(ww')\big),
\end{equation}
where $\lambda = g/c^2$, and by \eqref{profile} the wave profile which travels speed $(c,0)$
is
\begin{multline}\label{prof}
 S_t= \{\big(\xi+\mathcal{C} w(\xi-ct),w(\xi-ct)\big):\xi\in \mathbb{R}\}\\ =  (ct,0)+ \{\big(\xi+\mathcal{C} w(\xi),w(\xi)\big):\xi\in \mathbb{R}\}.
 \end{multline}
 It follows from \eqref{prof} that $ \{\big(\xi+\mathcal{C} w(\xi),w(\xi)\big):\xi\in \mathbb{R}\}$ is the fixed profile of a surface wave which travels without change of shape with  velocity $(c,0)$. Note that if $\eta$ is defined
  by $\eta (\xi +\mathcal{C} w(\xi))  = w(\xi)$, it follows that
 $$
 S_t= \{(x, \eta(x-ct): x \in \mathbb{R}\} \text{ where } x = \xi+\mathcal{C} w(\xi-ct).
 $$

 Since  $\lambda = g/c^2$ does not determine the sign of $c$, a solution $w$ of \eqref{P} yields two traveling waves with the same profile and speed, but travelling in opposite directions.

Now  let $v(\xi) =w(-\xi)$. Then from \eqref{conjugate} and  differentiation,
\begin{equation*} \mathcal{C} v(\xi) = -(\mathcal{C} w)(-\xi), \quad \mathcal{C} v'(\xi) =  (\mathcal{C} v)'(\xi) = (\mathcal{C} w')(-\xi),
\end{equation*}
and hence $w$ satisfies \eqref{P} if and only if $v$ satisfies \eqref{P}. So,  by \eqref{prof}, if $w$ satisfies \eqref{P},
\begin{align*}\notag
  \{\big(\xi+\mathcal{C} v(\xi&-ct),v(\xi-ct)\big):\xi\in \mathbb{R}\}
 = \{\big(\xi-\mathcal{C} w(-\xi+ct),w(-\xi+ct)\big):\xi\in \mathbb{R}\}\\\notag &= (ct,0)+\big\{\big(-(-\xi+ct)-\mathcal{C} w(-\xi+ct),w(-\xi+ct)\big):\xi\in \mathbb{R}\big\}\\ &=  (ct,0)+ \{\big(-\xi-\mathcal{C} w(\xi),w(\xi)\big):\xi\in \mathbb{R}\}.
 \end{align*}
 Therefore $ \{\big(-\xi-\mathcal{C} w(\xi),w(\xi)\big):\xi\in \mathbb{R}\}$ is the fixed profile of a another surface wave (the same one if $w$ is even) which travels without change of shape with  velocity $(c,0)$.
 Of course this observation, that the reflection $x \leftrightarrow -x$ of a steady travelling wave profile yields a steady wave travelling with the same speed  in the same direction, is a simple consequence of the reversibility of flows
 governed by system \eqref{fens}.

Since  $(\varphi,w)$ in \eqref{boxed} is related to $(\tilde \phi,\eta)$ in \eqref{fens} by
\begin{gather*} \omega(\xi,t) = \eta\big(\xi+\mathcal{C} w(\xi,t),t\big),~\varphi(\xi,t) = \tilde \phi\big(\xi+ \mathcal{C} w(\xi,t), w(\xi,t), t\big),\\
\intertext{and in the Stokes wave problem $\varphi(\xi,t) = \varPhi(\xi-ct)$  it follows  that  $\phi$ in \eqref{pens} satisfies}
\varPhi(\xi-ct) = c\phi\big(\xi-ct+\mathcal{C} w(\xi-ct),w(\xi-ct)\big) + c\big( \xi-ct+\mathcal{C} w (\xi-ct)\big),\\
\end{gather*}
Therefore, since $\varPhi' - c\,\mathcal{C} w' =0$, the independent variable $\xi$ is \eqref{P} is related to the dependent variable $\phi$ in \eqref{pens} by
$$\phi\big(\xi+\mathcal{C} w(\xi),w(\xi)\big)~ +~ \xi \text{ is constant}.$$  Recall similarly that   in Nekrasov's equation \eqref{nek} the independent variable $s= \phi(x,\eta(x))
$, also for $\phi$ in \eqref{pens}.
\qed

Equation \eqref{P} was first derived by Babenko \cite{Babenko}  for small amplitude steady travelling waves. However an equivalent identity  involving Fourier coefficients of Stokes waves had been uncovered earlier \cite{LongH:78}, but not written as  equation \eqref{P}. Babenko's equation was independently  re-discovered \cite{Balk, DKSZ1996,plotnikov:turning} and by now has been studied  extensively
for waves of all amplitudes, see  \cite{BuDaTo:98c, BDT2000b,BDT2000a,{BuTo:00},ShargTol:01,ShargTol:03, bernoulli,Toland:00}.\qed

{\it Dual Stokes waves.}   With  $\hat w= \lambda w^2-w$,  \eqref{P} can be written in symmetric form
$$\hat w'(1+\mathcal{C} w') + w'(1+\mathcal{C} \hat w') =0,
$$
which implies, see  \cite{BuTo:00,bernoulli}, that $\hat w$ also satisfies an equation similar to  \eqref{P} and so corresponds to another  free boundary problem on infinite depth, namely
\begin{subequations}\label{rens}\begin{align}\label{ra}
 &\Delta  \hat\psi (x,y) = 0 \text{ on } \hat\Omega:=\{(x,y):(y < \hat \eta(x))\};
\\& \nabla  \hat\psi(x,y) \text{ is $2\pi$-periodic in $x$}; \label{rb}
\\&\nabla  \hat\psi(x,y) \rightarrow (0,-1) \text{ as } y \rightarrow -\infty; \label{rc}
\\&  \hat\eta'(x)\partial_y \hat\psi(x,\hat \eta(x)) +\partial_x\hat\psi(x,\hat \eta(x)  = 0;\label{rd}
\\& (4\lambda+1) |\nabla  \hat\psi(x,\hat\eta(x))|^4  = 1,\quad \lambda = g/c^2.\label{re}
\end{align}
\end{subequations}
Although  \eqref{qens} and \eqref{rens} are equivalent, $\Omega \neq \hat \Omega$, $S \neq \hat{S}$, and by \eqref{re} the stream function $\hat \psi$ does not satisfy the constant
pressure boundary condition \eqref{qe}.
Thus any steady solution of  Babenko's equation \eqref{P}, or equivalently Nekrasov's equation, yields simultaneously solutions to two distinct Bernoulli free boundary problems.\qed

\subsection{Variational structure of Babenko's equation}

\begin{figure}
  \includegraphics[scale=1]{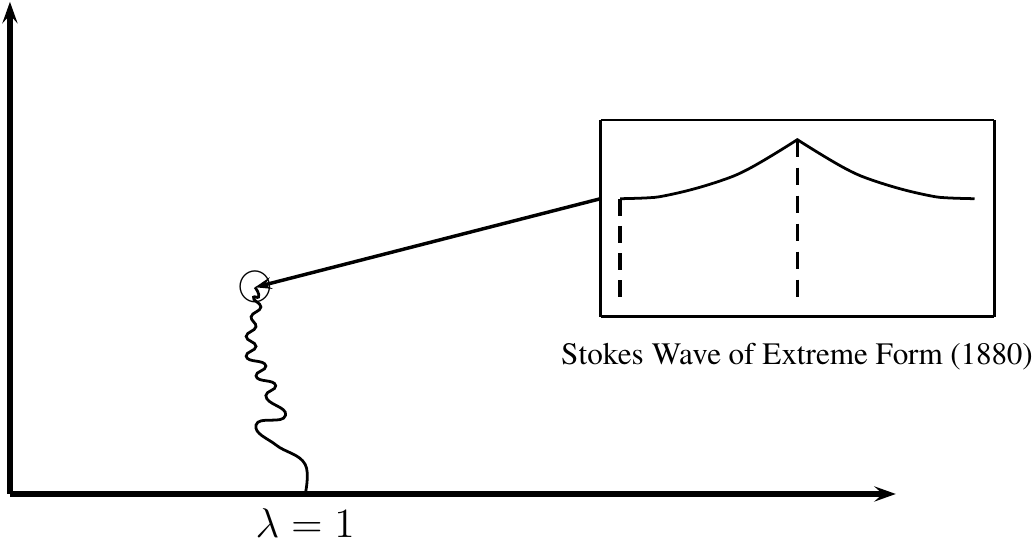}
  \caption{The primary branch}
  \label{Figure_Primary_Branch_Toland}
\end{figure}

The operator  $w \mapsto \mathcal{C} w'$  in equation \eqref{P} is  first-order, non-negative-definite, self-adjoint and densely defined on ${L^2_{2\pi}}$ by
$$\mathcal{C}(e^{ik})' = |k| e^{ik}, ~~ k\in \mathbb{Z}, \text{ whence }~ \mathcal{C} w' = \sqrt{-\frac {\partial^2~}{\partial \xi^2}~w}.$$
Therefore it
behaves somewhat like an elliptic differential operator
but lacks  a maximum principle. It is related to the $H^{1/2}(\mathbb S^1)$ norm of $w$  by the formula
$${\|w\|^2_{H^{1/2}(\mathbb S^1)} =\|w\|^2_{L_2(\mathbb S^1)} +\int_0^{2\pi}w\mathcal{C} w'\,dt}. $$
Let
\begin{equation}\label{EL} \mathcal{J}(\lambda,w)=\int_{-\pi}^{\pi} \big\{ w\mathcal{C} w' -
\lambda w^2(1+\mathcal{C} w')\big\}\,d\xi,\quad w \in W^{1,2}(\mathbb S^1).\end{equation}
Then $\mathcal{J}(\lambda,\cdot)$ is smooth on $W^{1,2}(\mathbb S^1)$ and, as is easily seen,   \eqref{P} is the Euler-Lagrange equation of $\mathcal{J}(\lambda,\cdot)$ in that setting.
Unfortunately $\mathcal{J}(\lambda,\cdot)$ is not well-defined on  $H^{1/2}(\mathbb S^1)$ (see \cite{Toland:00,To:01}) and
 is not bounded above or below, because of the cubic term.
These difficulties mean that in spite of some small-amplitude variational theory \cite{BSTa,BSTb}, so far there is
no self-contained variational proof of existence of critical points of
\eqref{EL} that yield large-amplitude solutions of \eqref{P}.
However, the variational formulation is  effective in alliance with other methods.
For example, equation \eqref{P}  can be written
\begin{equation}\label{magiceqn} (1-2\lambda w) \mathcal{C} w' = \lambda(w - (w\mathcal{C} w' -\mathcal{C}(ww')  )= \lambda(w-\mathcal{Q}(w))
\end{equation}
where, for $x \in [-\pi,\pi]$,
$$\mathcal{Q}(w)(x):= w(x)\mathcal{C} w'(x) -\mathcal{C}(ww')(x) = \frac{1}{8\pi} \int _0^{2\pi} \left\{\frac{ w(x)-w(y)}{\sin {\frac{1}{2}}(x-y)}\right\}^2 dy \geq 0.$$
Note $\mathcal{Q}(w) \in L_q$ for some $q>p$ when $w'\in L_p,\,p>1$, and
\eqref{magiceqn} implies
$$ \mathcal{C} w' = \lambda\left(\frac{w -\mathcal{Q} (w) }{1-2\lambda w}\right) \text{ when }  1-2\lambda w >0.
$$
 So difficulties  arises only when $1-2\lambda w$ has zeros, in which case $|\nabla \psi(x,\eta(x)| =0$ in \eqref{qe}
which corresponds to  Stokes waves of extreme form.
 Nowadays a very elementary application of bifurcation from a simple eigenvalue  yields their small-amplitude waves with $\lambda$ close to 1 (first proved by Nekrasov \cite{Nekrasov} and Levi-Civita
in the 1920s.)
Global bifurcation theory has more to say, but many questions remain unanswered, while
numerical evidence and \cite{PT2004} suggest the global bifurcation picture is as in Figure \ref{Figure_Primary_Branch_Toland}.

As the Stokes extreme wave, $1-2\lambda w(0) =0$, is approached, the Morse Index of solutions increases without bound
and, from numerical observation, $\|\theta\|_\infty$ in \eqref{nek} oscillates and the number of inflection points on the surfaces increases.
{\it The Morse Index  
${\mathcal M (\lambda,w)}$}
is the number of negative eigenvalues  of  $D^2\mathcal{J}[\lambda,w]$,
\begin{equation*}D^2 \mathcal{J}[\lambda,w]\,\phi\,= \mu \phi,\quad \mu <0,~~ \phi \neq 0,\end{equation*}
where $D^2\mathcal{J}[\lambda,w]$ is the linearisation of \eqref{P} with respect to $w$ at a solution $(\lambda, w)$:
\begin{equation*}\mathcal{C}\phi' -\lambda \big(\phi +\phi\mathcal{C} w' + w \mathcal{C} \phi ' +\mathcal{C} (w\phi)'\big) =D^2 \mathcal{J}[\lambda,w]\,\phi = \mu \phi,~~ (\lambda,w) \text{ satisfies \eqref{P}}. \end{equation*}

\subsection{Plotnikov's Theorem \cite{plotnikov:turning}}
Suppose
 a sequence $\{(\lambda_k,w_k)\}$ of solutions of \eqref{P} has
  $1-2\lambda_k w_k > 0$ and the Morse indices $\{\mathcal M(\lambda_k,w_k)\}$ are
 bounded. Then for some
 $\alpha > 0$
 $$\boxed{1-2\lambda_k w_k(x) \geq \alpha,\quad x \in \mathbb{R}, \quad k \in \mathbb{N}.}$$
Thus, by Plotnikov's theorem, Stokes waves  approaching extreme waves  becomes more and more unstable in the sense that the Morse indices
become unbounded.
Shargorodsky \cite{sharg} has quantified the relation between   Morse index and $\alpha$.\qed

We now mention the
implications of Plotnikov's Theorem for Stokes Waves \cite{BDT2000a}.
In abstract terms Babenko's equation for travelling waves is
\begin{equation}\label{dog} \mathcal{C} w' = \lambda \nabla \mathcal{L} (w)
 \text{ where }
\mathcal{L} ={\frac{1}{2}} \int_{-\pi}^{\pi} \left\{ w^2(1+\mathcal{C} w')\right\} \,dx.\end{equation}
By real-analytic function theory \cite{dancer,dancer:global,bt:analytic} there is a parameterized real-analytic curve of solutions $\{(\lambda_s,w_s):s \in [0, \infty)\}$ with
 $\mathcal M(\lambda_s,w_s) \to \infty$  as $s \to \infty$.

The variational structure of \eqref{dog}, with  a Lagrange multiplier $\lambda$ on the right,
implies that if the Morse Index changes as $s$ passes through $s^*$
then one of two
 things must happen: either there is  a crossing or a turning point. A point where one or other of these possibilities occurs will be called a bifurcation point on the primary branch. In the Stokes wave problem it is not known which of these possibilities occurs (this seems to be a very hard problem) but the numerical evidence points to all the bifurcation points being turning points.

\begin{figure}[h!]
\center
\includegraphics[scale=1]{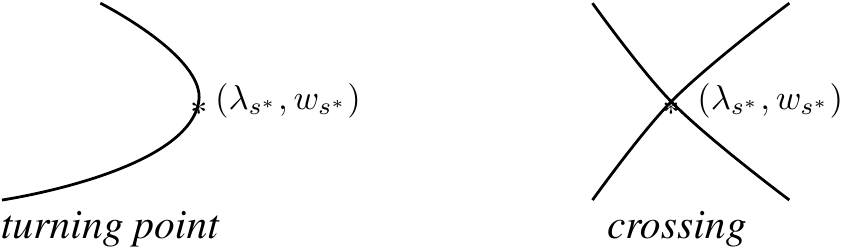}

\end{figure}

\subsection{Period-multiplying bifurcation on the primary branch.}
Motivated by numerical investigations of Chen \& Saffman \cite{chensaff} a lot  more can be extracted from the same point of view if the variational structure is developed in a slightly different setting.
 The solutions $\{(\lambda_s,w_s):s \in [0, \infty)\}$ with
 $\mathcal M(\lambda_s,w_s) \to \infty$  as $s \to \infty$ are $2\pi$-periodic, and hence  $2p\pi$-periodic for any prime number $p$. So consider $(\lambda_s,w_s)$ in a Banach space of $2p\pi$-periodic solutions for a prime number $p$, and let $\mathcal M_p(\lambda_s,w)$ denote their Morse index in this new setting.
Then,
for $p$  sufficiently large and prime, it can be shown that
$\mathcal M_p (\lambda_s,w_s)$ changes  as $s$ passes through $s_p^*$ where $s_p^*\neq s^*$ is close to $s^*$.
Since, for  primes $p$ sufficiently large there is, near $s^*$ on the primary branch, bifurcation point $s_p^*$ for solutions of  minimal period $2p\pi$, these
bifurcations must be crossings.

\begin{figure}[h!]
\center
\includegraphics[scale=1]{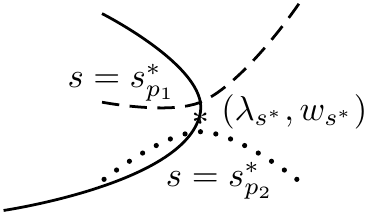}

\end{figure}

{\it Open question.}
It is known from methods of topological degree and real-analytic function theory that the global branch of Stokes waves ``terminates" at Stokes extreme wave
and Plotnikov's result guarantees the existence of solutions with arbitrarily large Morse index.
However, despite the apparently  simple form of $\mathcal{J}$, a satisfactory \emph{global
variational approach} to the  existence of Stokes waves capable of answering the question
\begin{center}
``for large $n \in \mathbb{N}$ does there exist a wave  with Morse index $n$?''
\end{center}
remains to be discovered. Nevertheless the variational approach and Morse index theory led to the
period-multiplying (sub-harmonic) bifurcation points near  turning points on the primary branch, as had been
 observed numerically by Chen \& Saffman \cite{chensaff}.\qed

\subsection{Summary of variational approach to Babenko's Equation.}
Equation \eqref{P} was introduced for small-amplitude steady  waves by Babenko \cite{Babenko}  in 1987, the year in which he died \cite{aabab}.
Independently, in 1992 Plotnikov \cite{plotnikov:turning},  and in 1996 {Balk \cite{Balk} } rediscovered equation \eqref{P}. In 1996
 Dyachenko, Kuznetsov, Spector \& Zakharov \cite{DKSZ1996} used conformal mappings to transform Zakharov's \cite{zak} Hamiltonian system \eqref{HSB} into \eqref{boxed},
but did not comment on its Hamilton structure. Apparently unaware
of \cite{Babenko}, they  derived \eqref{P} as the special case of \eqref{boxed} for travelling wave.

In his ground-breaking study of   non-uniqueness questions for
solitary waves, Plotnikov \cite{plotnikov:turning} introduced  Morse index
calculations  for the solitary-wave analogue of \eqref{P}. When Morse index considerations are taken into account with real-analytic  global bifurcation theory, this  led, in 2000,  to a proof \cite{BDT2000b, BDT2000a} that sub-harmonic bifurcation occurs when   periodic waves approach the Stokes wave of extreme form.
Then, in 2001, free boundaries  dual to Stokes waves were discovered  \cite{BuTo:00}.

In 2008,  the theory of \eqref{P} was generalised  to cover a wide class of free boundary problems
\cite{bernoulli},
but nothing emerged to make Stokes Waves seem special in that  wider class. See \cite{craik,craik2} for an historical account of water-wave theory.\qed

\subsection{Standing waves.}
When a fluid moves with temporal period $T$ when confined between verticals walls at which its surface remains horizontal, the result
is called a standing wave which is given by a solution  \eqref{boxed} which is periodic in space and time. Unlike the Stokes-wave problem,
the  standing-wave problem cannot can be treated as  time-independent in a moving frame. In fact, their theory, first considered by Sim\'eon Denis Poisson (1781--1840) who called them ``le clapotis" (lapping or sloshing waves), is much more complicated, as is clear from the observation that \eqref{boxed} linearised about the zero solution $\phi = w = 0$ is, for $x,t \in \mathbb{R}$,
\begin{gather*}
\dot w (x,t) =\mathcal{C} \phi'(x,t) ; \quad \dot \phi(x,t) + g w(x,t) = 0,~~ x \in [-\pi,\pi];~~ w'(\pm \pi, 0) = 0;\\
\big(w(x+2\pi, t),\phi(x+2\pi, t)\big) = \big(w(x,t),\phi(x,t)\big)= \big(w(x, t+T),\phi(x, t+T)\big).
\end{gather*}
Here $T$, the time period of $(w,\phi)$, is a parameter to be determined. From separation of variables it follows that there are no non-zero solutions $(w,\phi)$ if $\mu:=g(T/2\pi)^2$ is irrational, and if $\mu$ is rational,
 $\mu = n^2/m$ for infinitely many pairs $(m,n),\,m>0$. Hence the linearised problem has infinitely many solutions for every period $T$ for which  $\mu$ is positive rational:
\begin{align*}\phi(x,t) &= \sin \Big(\frac{2\pi nt}{T}\Big) \cos (mx), ~~ w(x,t) = -\frac{mT}{2\pi n} \cos \Big(\frac{2\pi nt}{T}\Big) \cos (mx).
\end{align*}
Hence, for a dense set of $T$ the linearised problem has infinitely many linearly independent solutions, and otherwise there are no solutions.
Nevertheless the Hamiltonian formulation \eqref{boxed}, combined with the Nash-Moser method for  small-divisor problems
from celestial mechanics
has led to the existence   \cite{ioosspt} of non-trivial   small amplitude  standing waves for  $\mu= g(T/2\pi)^2$
in a measurable set of positive density at $\mu = 1$.\qed

\section{Steady rotational periodic water waves in stratified media (by S.~V.~Haziot)}\label{sec:susanna}

\subsection{Introduction}\label{sect:euler_susanna}
The density of water, denoted by $\rho$, is a function of temperature, salinity and depth. It is strongly resistant to change for small variations of these quantities, and hence in many cases it is reasonable to assume $\rho$ to be a constant. However, when studying certain ocean waves, such as waves along the equator or internal waves at great depth, this assumption is no longer applicable. Any fluid with non-constant density is referred to as \emph{stratified} and is \emph{stably stratified} if $\rho$ is non-decreasing with depth.

We fix the wave speed $c>0$ and work in the moving frame coordinate system, $(X,Y):=(X-c,Y)$. We denote the fluid domain by $\Omega$, bounded below by the flat ocean bed $Y=0$ and above by the free surface $S$. In the first part, for simplicity, we will assume the free surface to be the graph of a function, that is $S=\{Y=H(X)\}$. In Section~\ref{subsect:conformal map_susanna}, we will parameterize $S$ to allow for waves with an overhanging profile.

We now recall Euler's incompressible equations of motion in the moving frame for stratified waves. These are expressed in terms of the velocity field $\textbf{u}=(u-c,v)$, the pressure distribution $P$, the density function $\rho>0$ and are coupled with the kinematic, \eqref{kinematic bottom_susanna}--\eqref{kinematic top_susanna}, and dynamic, \eqref{dynamic_susanna}, boundary conditions. For stratified waves, we additionally have the \emph{continuity equation}, see \eqref{continuity equation_susanna} below, expressing the fact that the density of a fluid particle remains constant throughout the flow.
\begin{subequations}\label{euler 2d_susanna}
	\begin{alignat}{2}	
	\label{incompressibility_susanna}
	u_X+v_Y&=0&\qquad&\text{in }\Omega,\\
	\label{continuity equation_susanna}
	(u-c)\rho_X+v\rho_Y&=0&\qquad&\text{in }\Omega,\\
	\label{momentum horizontal_susanna}
	\rho((u-c)u_X+vu_Y)&=-P_X&\qquad&\text{in }\Omega,\\
	\label{momentum vertical_susanna}
	\rho((u-c)v_X+vv_Y)&=-P_Y-g\rho&\qquad&\text{in }\Omega,\\
	\label{kinematic bottom_susanna}
	v&=0&\qquad&\text{on }Y=0,\\
	\label{kinematic top_susanna}
	(u-c)H_X-v&=0&\qquad&\text{on }S,\\
	\label{dynamic_susanna}
	P&=P_\text{atm}&\qquad&\text{on }S.
	\end{alignat}
\end{subequations}
Here $g$ denotes the gravitational constant and $P_\text{atm}$ the constant atmospheric pressure.

Finally, this paper will only deal with periodic symmetric waves, meaning that any solution to \eqref{euler 2d_susanna} will be $L$-periodic and even in the horizontal variable $X$. As is customary, we set $L:=2\pi/k$, where $k$ denotes the wave number. When studying solutions to \eqref{euler 2d_susanna}, we only need to consider the problem in a single period of the fluid domain. In Section~\ref{sec:sam:strat}, an in-depth review of the literature on stratified waves will be presented. In addition, different types of waves in stratified media, such as solitary waves and waves propagating along the interface of a multi-layered fluid, will be discussed.

\subsubsection{The pseudo stream function}\label{sect:stream function_susanna}
As in the homogeneous setting (when $\rho$ is constant), the system \eqref{euler 2d_susanna} can be reduced to a scalar elliptic problem. The incompressibility condition \eqref{incompressibility_susanna} ensures the existence of a stream function. However we cannot define it in the traditional way as the skew-gradient of the velocity field. Indeed, in that form it would not take into account the effects of stratification. Instead, we introduce the \emph{pseudo stream function} $\psi$ which contains an extra factor of $\rho$. It is defined by
\begin{equation*}
\psi_y=\sqrt{\rho}(u-c),\qquad\psi_x=-\sqrt{\rho}v.
\end{equation*}
\emph{Streamlines}, or the integral curves of the transformed vector field $\big(\sqrt{\rho}(u-c),\sqrt{\rho}v\big)$, are then level curves of $\psi$. Specifically, any particle on a given streamline will stay trapped there as the flow develops. In particular, from the kinematic boundary conditions \eqref{kinematic bottom_susanna}--\eqref{kinematic top_susanna}, we can see that the top and the bottom of the domain are streamlines. This enables us to normalize $\psi$ so that it vanishes along the free surface. Using the kinematic boundary conditions, we can then calculate the constant value $-m$ of $\psi$ along the bottom. Physically, $m$ represents the \emph{pseudo mass flux}.

From the continuity equation \eqref{continuity equation_susanna}, we see that $\rho$ is transported along particle trajectories. Hence it must be constant along streamlines and can therefore be expressed in terms of the pseudo stream function $\psi$ as
\begin{equation*}
\rho(X,Y)=:\varrho(\psi).
\end{equation*}
By reformulating \eqref{momentum horizontal_susanna}--\eqref{momentum vertical_susanna} in terms of $\psi$ and eliminating the pressure term, we see that $\Delta\psi+gy\varrho(\psi)$ and $\psi$ have the same level sets. This enables us to introduce the \emph{Yih--Long equation}:
\begin{equation}\label{yih long_susanna}
\Delta\psi+gY\varrho'(\psi)=\beta(\psi).
\end{equation}
Physically $\beta$ represents the variation of energy along streamlines.

Often when working with steady water waves, in both the homogeneous and heterogeneous setting, the additional assumption
\begin{equation}\label{no stagnation_susanna}
	\psi_Y>0\qquad\text{in }\Omega,
\end{equation}
is placed on the pseudo stream function. This assumption precludes the formation of internal stagnation points (points at which the gradient of the stream function vanishes), and forces $\beta$ to be a single-valued function $\mathbb{R}\to\mathbb{R}$. However, \eqref{no stagnation_susanna} is merely a sufficient condition for $\beta$ to exist and we are still free to assume the validity of \eqref{yih long_susanna} in the case when \eqref{no stagnation_susanna} fails: if we can prove that solutions to our problem exist, then the existence of $\beta$ is confirmed. Moreover, notice that in the case of constant density, the second term in the Yih--Long equation \eqref{yih long_susanna} vanishes and we recover an analogue of the vorticity equation for waves in a homogeneous fluid, with $\beta$ acting as the vorticity distribution.

By fixing $\beta$ and integrating \eqref{yih long_susanna}, we find that the energy in the system
\begin{equation}\label{energy_susanna}
	E=P+\tfrac{1}{2}|\nabla\psi|^2+g\varrho Y,
\end{equation}
is constant along streamlines. This is Bernoulli's law for stratified waves. We recover the dynamic boundary condition \eqref{dynamic_susanna} by evaluating \eqref{energy_susanna} along the free surface; see \eqref{fp dynamic_susanna} below.

Combining all these considerations, we obtain the following elliptic problem for stratified waves
\begin{subequations}\label{full problem streamfunction_susanna}
	\begin{alignat}{2}
		\label{fp yih long_susanna}
		\Delta\psi&=\beta(\psi)-gY\varrho'(\psi)&\qquad&\text{in }\Omega,\\
		\label{fp kinematic top_susanna}
		\psi&=0&\qquad&\text{on }S,\\
		\label{fp kinematic bottom_susanna}
		\psi&=-m&\qquad&\text{on }Y=0,\\
		\label{fp dynamic_susanna}
		|\nabla\psi|^2&=Q-2g\varrho(0)Y&\qquad&\text{on }S.
	\end{alignat}
\end{subequations}
Here the constant $Q$ is referred to as \emph{Bernoulli's constant}.

By general elliptic theory, \eqref{fp yih long_susanna}--\eqref{fp kinematic bottom_susanna} determine $\psi$ given the domain, while \eqref{fp dynamic_susanna} yields the domain, via $H(X)$, given $\psi$. This explains the necessity for the additional nonlinear dynamic boundary condition \eqref{fp dynamic_susanna}. The fact that the free surface $S$ is an unknown which needs to be determined as part of any solution is what makes the water wave problem so difficult. Indeed, the first step we must take is find a suitable change of variables with fixes the domain.

\subsubsection{Notation} Before we proceed, we gather some notation. Let $\Omega$ be an open connected subset of $\mathbb{R}^2$. For $\alpha\in(0,1)$ and $k\in\mathbb{N}$ we denote by $C^{k+\alpha}(\overline{\Omega})$ the space of functions whose partial derivatives up to order $k$ are Hölder continuous with exponent $\alpha$ over $\overline{\Omega}$. We denote by $C_{2\pi,\textup{e}}^{k+\alpha}(\overline{\Omega})$ the space of functions of class $C^{k+\alpha}$ that are $2\pi$-periodic and even in the $X$-variable.

\subsection{Transforming the domain}
\subsubsection{Dubreil-Jacotin transformation}
\begin{figure}
	\centering
	\includegraphics[scale=0.75]{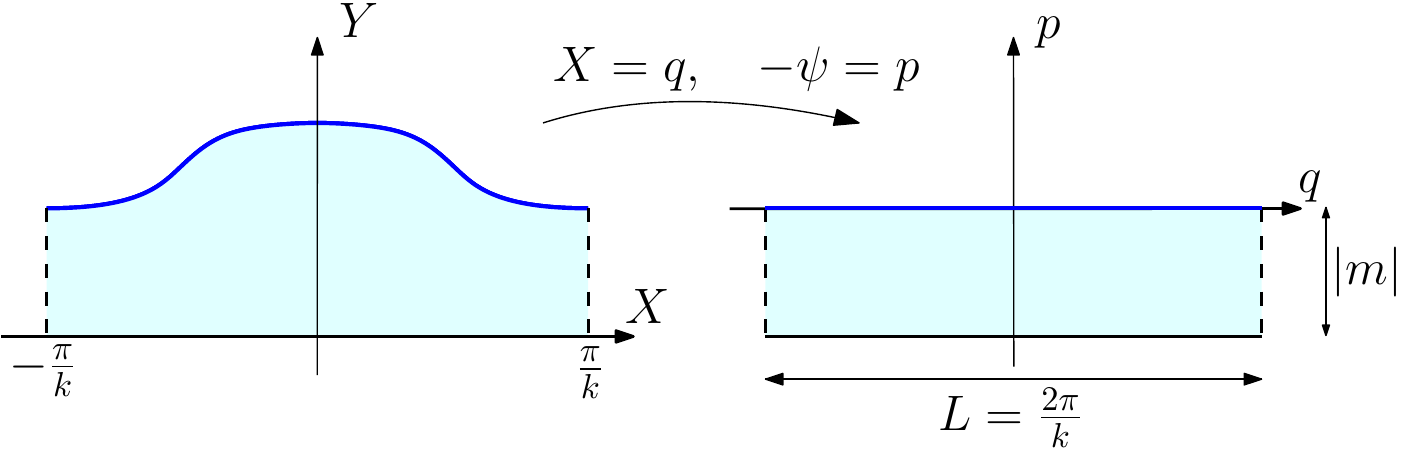}
	\caption{The Dubreil-Jacotin flattening transformation.
		\label{fig:dubreil_jacotin_susanna}
	}
\end{figure}

Typically, when studying rotational water waves in either a homogeneous or a stratified media for which \eqref{no stagnation_susanna} holds, we fix the domain using an elegant semi-hodograph change of variables due to Dubreil-Jacotin; see \cite{dubreil}. This transformation consists of straightening all the streamlines, including the free surface. In effect, one period in the fluid domain is mapped into a known rectangle of height $|m|$ and length $L$; see Figure~\ref{fig:dubreil_jacotin_susanna}. More precisely, we introduce the change of variables
\begin{equation}\label{DJ_susanna}
	(q,p)\mapsto(X,-\psi)
\end{equation}
along with the \emph{height function} defined by $h(q,p)=Y.$ This function indicates the height above the flat bottom on a streamline corresponding to $p$ at $X=q$. When mapped into the rectangle, the physical problem \eqref{full problem streamfunction_susanna} becomes a quasi-linear elliptic problem in terms of $h$.

In 2004, Constantin and Strauss, see \cite{CS2004}, used this transformation to prove the existence of large-amplitude rotational steady periodic water waves using global bifurcation theory. Specifically, they constructed a global curve of symmetric periodic water wave solutions which limits to stagnation. This result was the first of its kind and the approach was later extended to the setting of large-amplitude stratified waves \cite{walsh:stratified}, solitary waves \cite{wheeler:solitary}, stratified solitary waves \cite{cww:strat}, and many more!

As the change of variables \eqref{DJ_susanna} indicates, the assumption \eqref{no stagnation_susanna} is crucial for solutions in terms of $h$ in the rectangle to represent physical solutions to \eqref{full problem streamfunction_susanna}. As a result, this transformation is only applicable in settings for which all streamlines, including the free surface, are graphical. In other words, water wave solutions constructed using this approach cannot admit internal stagnation points or have overhanging wave profiles.

\subsubsection{Naive flattening transformation}
Another important transformation consists of using the change of variables
\begin{equation*}
	(X,Y)\mapsto\bigg(X,\frac{Y}{H(X)} \bigg),
\end{equation*}
which only flattens the free surface. This approach was used by Wahl\'en in \cite{wahlen:crit} to construct small-amplitude periodic waves with constant vorticity for which the flow contains closed streamlines in the form of \emph{Kelvin cat eyes}. Recently, Varholm, see \cite{varholm:global}, used this transformation to prove the existence of large-amplitude waves with a general vorticity distribution. Here, the global curve of solutions either limits to stagnation or loops back and reconnects. Similarly as for the Dubreil-Jacotin approach, this change of variables requires that the free surface be the graph of a function, and hence does not allow for overturning waves.

\subsubsection{Conformal map}\label{subsect:conformal map_susanna}
A rigorous existence proof for steady water waves with an overhanging wave profile remains to this day one of the biggest open problems in the field. However, both field data and numerical evidence strongly suggest that such waves do exist. For example, for periodic waves, the study in \cite{sp:steep} indicates that overhanging waves develop as the amplitude of a small-amplitude wave with closed streamlines drastically increases.

\begin{figure}
	\centering
	\includegraphics[scale=0.9]{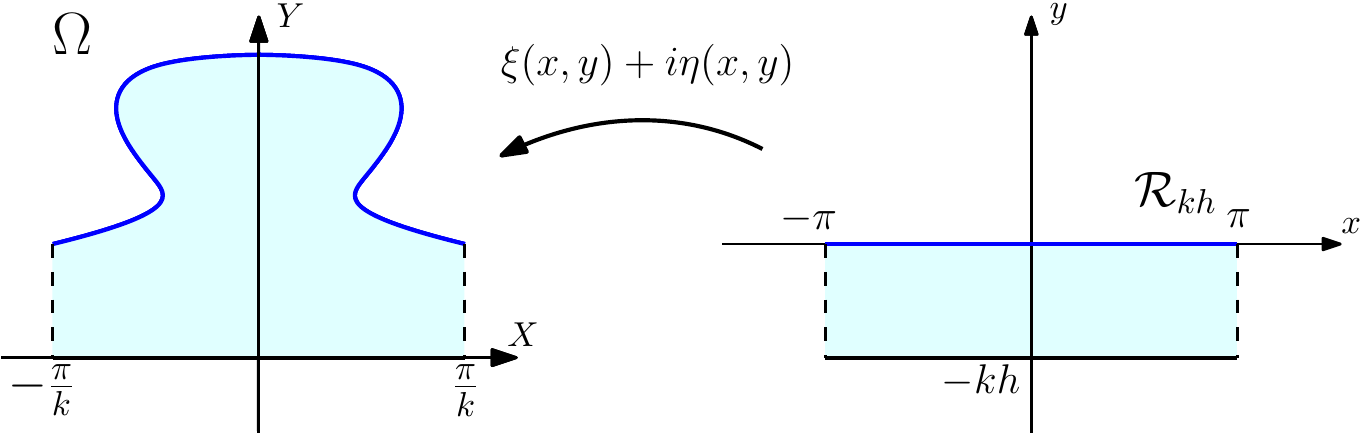}
	\caption{The conformal parametrization of the fluid domain $\Omega$.
		\label{fig:mapping_susanna}
	}
\end{figure}

In 2016, Constantin, Strauss and Varvaruca (see \cite{csv:critical}) proved the existence of large-amplitude periodic waves with constant vorticity which may have both an overhanging profile and internal stagnation points. In order to construct these waves, they viewed the physical fluid domain as the image via a conformal map of a rectangular strip; see Figure~\ref{fig:mapping_susanna}.

Let us outline the idea of this transformation. To begin with, in order to construct waves which may potentially overturn, one can no longer view the free surface $S$ as the graph of a function. Therefore, we parameterize it by two functions $\xi$ and $\eta$ such that
\begin{equation}\label{surface_susanna}
	S=\{\big(\xi(s,0),\eta(s,0)\big):s\in\mathbb{R} \},
\end{equation}
where the map $s\mapsto\big(\xi(s,0)-s,\eta(s,0)\big)$ is periodic of period $L$.

In \cite{cv:constant}, it is shown that for such a domain there exists a unique positive constant $h$, called the \emph{conformal mean depth}, such that one can find a conformal map $X+iY=\xi+i\eta$ from the strip
\begin{equation}\label{strip_susanna}
	\mathcal{R}_{kh}:=\{(x,y)\in\mathbb{R}^2:-kh<y<0 \}
\end{equation}
to $\Omega$. The wave number $k$ appears in \eqref{strip_susanna} due to rescaling the $x$ variable so that the image of a rectangle of length $2\pi$ is mapped to one period  of length $L$ in $\Omega$; see Figure~\ref{fig:mapping_susanna}. Specifically,
\begin{equation*}
\xi(x+2\pi,y)=\xi(x,y)+L\qquad\text{and}\qquad\eta(x+2\pi,y)=\eta(x,y),
\end{equation*}
for all $(x,y)\in\mathcal{R}_{kh}$. Moreover, the map is constructed such that the top of $\mathcal{R}_{kh}$ is mapped to $S$ and the bottom to $Y=0$. Finally, assuming that $\Omega$ is of class $C^{1+\alpha}$, the Kellogg--Warschawski Theorem (see \cite[Theorem 3.6]{pommerenke:book}) yields that $\xi,\eta\in C^{1+\alpha}(\overline{\mathcal{R}_{kh}})$ and that
\begin{equation*}
	\xi_x^2(x,0)+\eta_x^2(x,0)\neq0\qquad\text{for all }x\in\mathbb{R}.
\end{equation*}

It now suffices to determine $\eta$. Indeed, since $\eta$ is harmonic in $\mathcal{R}_{kh}$, using the Cauchy--Riemann equations, we can, up to a constant, determine its harmonic conjugate $\xi$. Provided the map $\xi+i\eta$ is injective on the surface of the domain, the Darboux--Picard Theorem, (see \cite[Corollary 9.16]{burckel:book}) ensures that $\xi+i\eta$ defines a conformal map, thus enabling us to determine the free surface $S$ from \eqref{surface_susanna}. In order to solve for $\eta$, the water wave problem is then reformulated in the conformal strip as a quasilinear pseudo differential system on the surface $y=0$. This formulation involves the periodic Hilbert transform and only depends on the function $\eta(x,0)$ and on several parameters. The specifics on how to achieve this are presented in the next section.

\subsection{Some large-amplitude stratified waves}

\subsubsection{Formulation of the problem}

\begin{figure}
	\centering
	\includegraphics[scale=0.9]{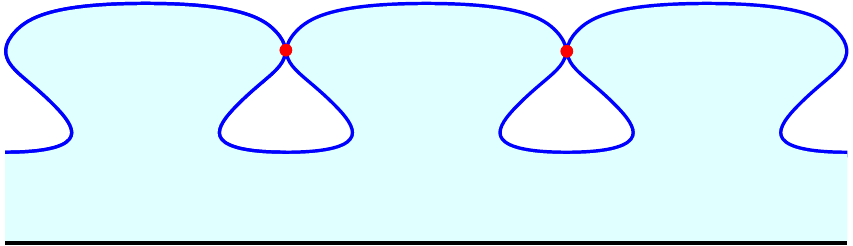}
	\caption{The free surface overturns and self-intersects along the trough line.
		\label{fig:self-intersection_susanna}
	}
\end{figure}

We describe the approach used in \cite{csv:critical} and extend it to a certain class of stratified waves satisfying the assumptions
\begin{equation}\label{assumptions_susanna}
	\beta(\psi)=\gamma\qquad\text{and}\qquad\varrho(\psi)=A\psi+B\qquad\text{for }\gamma,A\in\mathbb{R}\text{ and }B\in\mathbb{R}^+.
\end{equation}
From the considerations in Section~\ref{sect:euler_susanna}, for the fluid to be stably stratified we would expect a sign requirement on $\varrho'=A$. However, we refrain from doing this since an overhanging wave will necessarily be unstably stratified when it overturns. Implementing \eqref{assumptions_susanna} into \eqref{full problem streamfunction_susanna}, we obtain
\begin{subequations}\label{op stream_susanna}
	\begin{alignat}{2}
	\Delta\psi&=\gamma+aY&\qquad&\text{in }\Omega,\\
	\psi&=0&\qquad&\text{on }S,\\
	\psi&=-m&\qquad&\text{on }Y=0,\\
	|\nabla\psi|^2+bY&=Q&\qquad&\text{on }S,
	\end{alignat}
\end{subequations}
where for simplicity of notation, we denote $a:=-Ag$ and $b=2Bg$. Here $S$ is defined as in \eqref{surface_susanna}. We will now outline the proof of the following informally stated theorem.
\begin{theorem}[\cite{haziot2021stratified}]
	Fix a Hölder exponent $\alpha\in(0,1)$, a wave speed $c>0$, a wave number $k>0$, the gravitational constant $g>0$, and $\beta$ and $\varrho$ as in \eqref{assumptions_susanna}. Then there exists a family of global curves $\mathscr{C}$ of periodic, symmetric, monotone water wave solutions to \eqref{op stream_susanna}, parameterized by $s$, $0<s<\infty$. As $s\to\infty$, one of the following alternatives occurs:
	\begin{enumerate}[label=\rm(\roman*)]
		\item \label{greatest height_susanna} $\eta$ approaches a wave of greatest height with a stagnation point at the crest;
		\item \label{self-intersection_susanna} the wave profile overturns and self-intersects above the trough line (see Figure~\ref{fig:self-intersection_susanna});
		\item \label{parameter blowup_susanna} the parameters $m$ and $Q$ blow up;
		\item \label{eta blowup_susanna} the wave elevation $\eta$ blows up in the Hölder norm $C^{2+\alpha}$.
	\end{enumerate}
\end{theorem}

To begin with, the advantage of assumption \eqref{assumptions_susanna} is that we can introduce the harmonic function
\begin{equation*}
	\zeta(x,y):=\psi(\eta(x,y),\xi(x,y))+m-\tfrac{1}{2}\gamma\eta^2(x,y)-\tfrac{1}{6}a\eta(x,y)^3.
\end{equation*}
In terms of the conformal variables, \eqref{op stream_susanna} becomes
\begin{subequations}\label{full problem_susanna}
	\begin{alignat}{2}
		\label{zeta harmonic_susanna}
		\Delta\zeta&=0&\qquad&\text{in }\mathcal{R}_{kh},\\
		\label{zeta kinematic_susanna}
		\zeta&=m-\tfrac{1}{2}\gamma\eta^2-\tfrac{1}{6}a\eta^3&\qquad&\text{on }S,\\
		\label{zeta bottom_susanna}
		\zeta&=0&\qquad&\text{on }y=-kh,\\
		\label{zeta dynamic_susanna}
		\big(\zeta_y+\gamma\eta\eta_y+\tfrac{1}{2}a\eta^2\eta_y\big)&=(Q-b\eta)\big(\eta_x^2+\eta_y^2\big)&\qquad&\text{on }S,		
	\end{alignat}
which we couple with
	\begin{alignat}{2}
		\label{eta harmonic_susanna}
		\Delta\eta&=0&\qquad&\text{in }\mathcal{R}_{kh},\\
		\eta&=:\overline{\eta}(x)&\qquad&\text{on }S,\\
		\label{eta bottom_susanna}
		\eta&=0&\qquad&\text{on }y=0.
	\end{alignat}
It can be checked that
\begin{equation}\label{eta average_susanna}
	\big[\overline{\eta}\big]=h,
\end{equation}
 where the $[\,\cdot\,]$ denotes the mean over one period, and we will additionally require that
	\begin{equation}\label{regularity_susanna}
		\zeta,\eta\in C_{2\pi,\textup{e}}^{2+\alpha}(\mathcal{R}_{kh}).
	\end{equation}	
\end{subequations}
Since $\eta$ is harmonic, its harmonic conjugate $\xi$ can be determined uniquely up to a constant, so that $\xi+i\eta$ is holomorphic in $\mathcal{R}_{kh}$. Specifically, we have
\begin{equation}\label{xi_susanna}
	\xi(x,y)=\frac{\big[\overline{\eta}\big]}{kh}x+\xi_0(x,y),
\end{equation}
where the function $\xi_0$, $2\pi$-periodic in $x$, has been normalized so that it has zero mean over one period.

Let us consider the problem \eqref{eta harmonic_susanna}--\eqref{eta bottom_susanna} more closely. The Dirichlet--Neumann operator for the strip $\mathcal{R}_{kh}$ is defined by
\begin{equation*}
	\mathcal{G}_{kh}(\overline{\eta})(x)=\eta_y(x,0),
\end{equation*}
for all $x\in\mathbb{R}$. This is a bounded linear operator from $C_{2\pi}^{2+\alpha}(\mathbb{R})$ to $C_{2\pi}^{1+\alpha}(\mathbb{R})$. The Dirichlet--Neumann operator admits a conjugation operator $\mathcal{C}_{kh}$, called the periodic Hilbert transform. It is defined for $2\pi$-periodic functions of zero mean and is a bounded invertible linear operator from $C_{2\pi}^{2+\alpha}(\mathbb{R})$ into itself. It is defined by
\begin{equation*}
	\mathcal{C}_{kh}\big(\overline{\eta}-[\overline{\eta}]\big):=\xi_0(x,0),
\end{equation*}
which, combined with \eqref{xi_susanna}, yields
\begin{equation*}
	\xi(x,0)=\frac{\big[\overline{\eta}\big]}{kh}x+\big(\mathcal{C}_{kh}(\overline{\eta}-[\overline{\eta}])\big)(x),
\end{equation*}
for all $x\in\mathbb{R}$. Using the Cauchy--Riemann equations and \eqref{eta average_susanna}, we now see that the operators $\mathcal{G}_{kh}$ and $\mathcal{C}_{kh}$ are related by the identity
\begin{equation}\label{etay_susanna}
	\mathcal{G}_{kh}(\overline{\eta})=\eta_y(x,0)=\xi_x(x,0)=\frac{\big[\overline{\eta}\big]}{kh}+\big(\mathcal{C}_{kh}(\overline{\eta}-[\overline{\eta}]) \big)'=\frac{1}{k}+\mathcal{C}_{kh}(\overline{\eta}')(x).
\end{equation}
Similarly, from \eqref{zeta harmonic_susanna}--\eqref{zeta bottom_susanna}, we get
\begin{equation}\label{zetay_susanna}
	\zeta_y(x,0)=\frac{m}{kh}-\frac{\gamma}{2kh}\big[\overline{\eta}^2\big]-\gamma\mathcal{C}_{kh}(\overline{\eta}\overline{\eta}')-\frac{a}{6kh}\big[\overline{\eta}^3\big]-\frac{a}{2}\mathcal{C}_{kh}(\overline{\eta}^2\overline{\eta}')(x).
\end{equation}
Substituting \eqref{etay_susanna} and \eqref{zetay_susanna} in the dynamic boundary condition \eqref{zeta dynamic_susanna}, we effectively reduce the elliptic system problem \eqref{full problem_susanna} to the following nonlinear pseudo differential scalar equation
\begin{subequations}\label{first reformulation_susanna}
\begin{alignat}{2}\label{nonlinear pseudodiff_susanna}
&\Bigg(\frac{m}{kh}-\frac{\gamma}{2kh}[\overline{\eta}^2]-\gamma\mathcal{C}_{kh}(\overline{\eta}\overline{\eta}')-\frac{a}{6kh}[\overline{\eta}^3]-\frac{a}{2}\mathcal{C}_{kh}(\overline{\eta}^2\overline{\eta}')\\&\quad+\gamma \overline{\eta}\bigg(\frac{1}{k}+\mathcal{C}_{kh}(\overline{\eta}') \bigg) +\frac{a}{2}\overline{\eta}^2\bigg(\frac{1}{k}+\mathcal{C}_{kh}(\overline{\eta}') \bigg) \Bigg)^2\nonumber=(Q-b\overline{\eta})\Bigg((\overline{\eta}')^2+\bigg(\frac{1}{k}+\mathcal{C}_{kh}(\overline{\eta}')\bigg)^2 \Bigg),
\end{alignat}
for all $x\in\mathbb{R}$, for
\begin{equation}\label{regularity first formulation_susanna}
	(\overline{\eta},h)\in C_{2\pi,e}^{2+\alpha}(\mathbb{R})\times\mathbb{R}^+\qquad\text{and}\qquad\big[\overline{\eta}\big]=h.
\end{equation}
Moreover, any solution $(\overline{\eta},h)$ to \eqref{nonlinear pseudodiff_susanna}--\eqref{regularity first formulation_susanna} leads to a solution of \eqref{op stream_susanna} provided the following addition requirements are met:
\begin{alignat}{2}
\label{positive free surface_susanna}
&\overline{\eta}(x)>0&\qquad&\text{for all }x\in\mathbb{R},\\
\label{injectivity_susanna}
&x\mapsto\bigg(\frac{x}{k}+\mathcal{C}_{kh}(\overline{\eta}-h)(x),\overline{\eta}(x) \bigg)&\qquad&\text{is injective on }\mathbb{R},\\
\label{singularity_susanna}
&(\overline{\eta}'(x))^2+\bigg(\frac{1}{k}+\mathcal{C}_{kh}(\overline{\eta}')\bigg)^2\neq0&\qquad&\text{for all } x\in\mathbb{R}.
\end{alignat}
Condition \eqref{positive free surface_susanna} ensures that the free surface remains positive, and \eqref{injectivity_susanna}, that it does not self-intersect. Finally, \eqref{singularity_susanna} guarantees that we have no singularities in the conformal map.

The problem admits a family of \emph{trivial solutions}. These consist of the simplest solutions to our problem: they are independent of the $x$-variable and have a flat free surface over parallel streamlines. Clearly, any such solution would be of the form $\overline{\eta}\equiv h$. Hence, it will be convenient to work with the function
\begin{equation}\label{definition w_susanna}
w:=\overline{\eta}-h,
\end{equation}
since, by \eqref{regularity first formulation_susanna}, we must necessarily have $\big[w\big]=0$.
\end{subequations}

For irrotational water waves, via Riemann--Hilbert theory the analogue of \eqref{nonlinear pseudodiff_susanna} can be reformulated as a quasi-linear pseudo differential equation known as \emph{Babenko's equation}. In \cite{csv:critical}, the authors showed that, with substantially more work, a Babenko-type formulation also exists for waves with constant density. We extend their approach to our stratified setting. Specifically, provided \eqref{singularity_susanna} holds, we can show using Riemann--Hilbert theory, that \eqref{nonlinear pseudodiff_susanna} expressed in terms of $w$ is equivalent to the following system:
\begin{subequations}\label{babenko_susanna}
\begin{alignat}{2}\label{babenko equation_susanna}
&\mathcal{C}_{kh}\Bigg(\bigg(Q-b(w+h)-\gamma^2(w+h)^2-\frac{a^2}{4}(w+h)^4-\gamma a(w+h)^3 \bigg)w'\Bigg)\nonumber\\&+ \Bigg(\bigg(Q-b(w+h)-\gamma^2(w+h)^2-\frac{a^2}{4}(w+h)^4-\gamma a(w+h)^3 \bigg)\bigg(\frac{1}{k}+\mathcal{C}_{kh}(w') \bigg)\Bigg)\nonumber\\&-\big(2\gamma(w+h)+a(w+h)^2	 \big)\bigg(\frac{m}{kh}-\frac{\gamma}{2kh}([w^2]+h^2)\nonumber\\&-\gamma\mathcal{C}_{kh}\big((w+h)w'\big)-\frac{a}{6kh}\big([w^3]+3h[w^2]+h^3\big)-\frac{a}{2}\mathcal{C}_{kh}\big((w+h)^2w'\big) \bigg)-\kappa=0,
\end{alignat}
for $(m,Q,w)\in\mathbb{R}\times\mathbb{R}\times C_{2\pi}^{2+\alpha}(\mathbb{R})$, coupled with the scalar constraint
\begin{alignat}{2}\label{babenko scalar_susanna}
\Bigg[\Bigg(\frac{m}{kh}-\frac{\gamma}{2kh}([w^2]+h^2)-\mathcal{C}_{kh}((w+h)w')-\frac{a}{6kh}([w^3]+3h[w^2]+h^3)\nonumber\\-\frac{a}{2}\mathcal{C}_{kh}((w+h)^2w')+(\gamma(w+h)+\frac{a}{2}(w+h)^2)\left(\frac{1}{k}+\mathcal{C}_{kh}(w')\right) \Bigg) \Bigg]\nonumber\\-\left[(Q-b(w+h))\left((w')^2+\left(\frac{1}{k}+\mathcal{C}_{kh}(w') \right)^2 \right) \right]=0.
\end{alignat}
\end{subequations}
In \eqref{babenko equation_susanna}, $\kappa$ consists of scalar terms whose exact formulation is not important to our analysis.

By removing all effects of stratification in \eqref{babenko_susanna}, that is, by deleting all terms involving $a$, we recover the Babenko-type system derived for constant vorticity in \cite{csv:critical}, with $\gamma$ acting as the vorticity distribution. If we additionally set $\gamma$ equal to zero, we obtain an analogue of Babenko's equation for irrotational waves in a fluid of infinite depth.  With our reformulation \eqref{babenko_susanna} now in hand, we prove the existence of large-amplitude water wave solutions to our problem \eqref{op stream_susanna}. This will be done via bifurcation theory.

\subsubsection{Functional analytic formulation}
To begin with, we see that when plugging the trivial solution $w=0$ into \eqref{babenko_susanna}, we find that the parameters $Q$ and $m$ are related by
\begin{equation*}
Q=bh+\bigg(\frac{m}{h}+\frac{\gamma h}{2}+\frac{ah^2}{3}\bigg)^2.
\end{equation*}
Since this identity is not necessarily true for nontrivial solutions, we introduce new parameters $(\mu,\lambda)$, such that
\begin{equation*}
	\mu=Q-bh-\lambda^2,\qquad\text{where}\qquad\lambda:=\frac{m}{h}+\frac{\gamma h}{2}+\frac{ah^2}{3},
\end{equation*}
denotes the speed of the particles of the trivial solutions at the free surface. The parameter $\lambda$ will serve as a bifurcation parameter for small-amplitude waves.

In order to apply any kind of bifurcation theory, we need to put our problem in a suitable functional analytic setting. We first define the Banach spaces
\begin{equation*}
\mathcal{X}=\mathbb{R}\times C_{2\pi,e}^{2+\alpha}\qquad\text{and}\qquad\mathcal{Y}=C_{2\pi,e}^{1+\alpha}\times\mathbb{R},
\end{equation*}
along with the open subset
\begin{equation*}
\mathcal{O}=\{(\lambda,(\mu,w))\in\mathbb{R}\times\mathcal{X}:\mu+\lambda^2-bw(x)>0 \text{ for all }x\in\mathbb{R} \}\subset\mathbb{R}\times\mathcal{X}.
\end{equation*}
In terms of the new parameters $(\mu,\lambda)$, we can now express \eqref{babenko_susanna} as the operator equation
\begin{equation}\label{operator equation_susanna}
	\mathscr{F}\big(\lambda,(\mu,w)\big)=0,
\end{equation}
where $\mathscr{F}\colon\mathbb{R}\times\mathcal{X}\to\mathcal{Y}$ is given such that $\mathscr{F}:=(\mathscr{F}_1,\mathscr{F}_2)$. Here, $\mathscr{F}_1$ is defined by the right hand side of \eqref{babenko equation_susanna} and $\mathscr{F}_2$ by the right hand side of \eqref{babenko scalar_susanna}.

\subsubsection{Local bifurcation theory}
The idea is now the following. We begin with the family of trivial solutions $w=0$ to \eqref{operator equation_susanna}. They represent a curve which is parameterized by $\lambda$. Since the trivial solutions are not particularly interesting, we hope to find values of $\lambda$ for which there is a change (specifically, an increase) in the number of solutions to \eqref{operator equation_susanna}. In other words, when $\lambda$ passes a certain threshold value, the flat surface would be \emph{perturbed} leading to the emergence of small-amplitude waves. Any such value of $\lambda$ is referred to as a \emph{bifurcation point} and clearly the implicit function theorem applied to $\mathscr{F}$ must fail at these points. As a result, we can find all potential bifurcation points by computing the Fréchet derivative of \eqref{operator equation_susanna} and identifying for which values of $\lambda$ this derivative vanishes.

We hence find that any bifurcation point must satisfy
\begin{equation}\label{bifurcation points_susanna}
\lambda(\gamma+ah)+\frac{b}{2}=\lambda^2nk\coth(nkh),\qquad\forall n\in\mathbb{N},
\end{equation}
and denote the set of all these points by $\{\lambda_{n,\pm}^*:n\in\mathbb{N} \}$. By studying $\partial_{(\mu,w)}\mathscr{F}\big(\lambda^*(0,0)\big)$, the linearized operator at the bifurcation parameter and the trivial solution, we find that it is Fredholm of index zero with a one-dimensional kernel, generated by $(0,w^*)\in\mathcal{X}$ where $w^*(x)=\cos(nx)$. In addition, using \eqref{bifurcation points_susanna}, we see that the transversality condition
\begin{equation*}
	\partial_\lambda\partial_{(\mu,w)}\mathscr{F}\big(\lambda^*(0,0)\big)(1,(0,w^*))\notin\operatorname{ran}\big(\partial_{(\mu,w)}\mathscr{F}\big(\lambda^*(0,0)\big)\big),
\end{equation*}
is satisfied (see \cite{haziot2021stratified} for details). Consequently, the Crandall--Rabinowitz local bifurcation theorem (see \cite{rabinowitz:simple}) yields the existence of a whole family of continuous local curves, bifurcating from the curve of trivial solutions at every $\lambda^*\in\{\lambda_{n,\pm}^*:n\in\mathbb{N} \}$. The solutions along these local curves are referred to as \emph{small-amplitude} solutions.

\subsubsection{Global bifurcation theory}
The next step is to extend each local curve to a global one. Solutions on the global curve will no longer be mere perturbations of the flat free surface and are hence referred to as \emph{large-amplitude} solutions. This continuation is achieved via analytic global bifurcation theory. Specifically, we use a theorem originally due to Dancer (see \cite{dancer}) which was later improved to be applicable to the study of water waves by Buffoni and Toland; see \cite{bt:analytic}. The theorem yields the existence of a global continuation of the local curve, which either limits to a blow-up scenario, or reconnects with the curve of trivial solutions.

Certain requirements need to be satisfied in order to apply this result. Provided $\mathscr{F}\colon\mathcal{O}\to\mathcal{Y}$ is real-analytic, we require the existence of a local curve of solution (as the one obtained from the local bifurcation theorem), and the following additional assumptions
\begin{enumerate}
	\item \label{global fredholm_susanna} the linearized operator $\partial_{(\mu,w)}\mathscr{F}\big(\lambda(\mu,w)\big)$ must be Fredholm of index zero for any $\big(\lambda,(\mu,w)\big)\in\mathcal{O}$ such that \eqref{operator equation_susanna} holds;
	\item \label{global compactness_susanna} for some sequence $(\mathcal{Q}_j)_{j\in\mathbb{N}}$ of bounded closed subsets of $\mathcal{O}$ with $\mathcal{O}=\cup_{j\in\mathbb{N}}\mathcal{Q}_j$, the set $\{(\lambda,(\mu,w))\in\mathcal{O}:\mathscr{F}(\lambda,(\mu,w))=0 \}\cap\mathcal{Q}_j,$ is compact for each $j\in\mathbb{N}$.
\end{enumerate}

In order to show that the operator $\mathscr{F}$ satisfies \eqref{global fredholm_susanna}--\eqref{global compactness_susanna}, we reorganize $\mathscr{F}_1$ as follows:
\begin{equation}\label{F1_susanna}
\mathscr{F}_1(\lambda,(\mu,w))=2(\mu+\lambda^2-bw)\mathcal{C}_{kh}(w')+4[w\mathcal{C}_{kh}(w')]+J(\lambda,w).
\end{equation}
Here $J(\lambda,w):=J_1(w)+J_2(\lambda,w)$, where $J_1(w)$ consists of a sum of commutators of the Hilbert transform of the form $f\mathcal{C}_{kh}(g)-\mathcal{C}_{kh}(fg)$ and $J_2(\lambda,w)$ gathers up all the remaining terms (a polynomial in $w$ and scalar terms which depend on $w$). The heart of the argument now relies on the following theorem concerning the regularity of commutators of the periodic Hilbert transform.
\begin{theorem}[see \cite{cv:constant}]\label{thm: hilbert reg_susanna}
	If $f\in C_{2\pi}^{j,\alpha}(\mathbb{R})$ and $g\in C_{2\pi}^{j-1,\alpha}(\mathbb{R})$, with $j\in\mathbb{N}$ and $\alpha\in(0,1)$, then $f\mathcal{C}_{kh}(g)-\mathcal{C}_{kh}(fg)\in C_{2\pi}^{j,\delta}(\mathbb{R})$ for all $\delta\in(0,\alpha)$, and there exists a constant $C:=C(j,\alpha,\delta)$ such that $$\|f\mathcal{C}_{kh}(g)-\mathcal{C}_{kh}(fg)\|_{j,\delta}\leq C\|f\|_{j,\alpha}\|g\|_{j-1,\alpha}.$$
\end{theorem}
Indeed, Theorem~\ref{thm: hilbert reg_susanna} yields that $J_1(w)$  maps bounded sets of $C_{2\pi}^{2,\alpha}(\mathbb{R})$ into bounded sets of $C_{2\pi}^{2,\delta}(\mathbb{R})$ and therefore into relatively compact subsets of $C_{2\pi}^{1,\alpha}(\mathbb{R})$. With this in hand, one can then show that the linearized operator $\partial_{(\mu,w)}\mathscr{F}\big(\lambda,(\mu,w)\big)$ is, for any $\big(\lambda,(\mu,w)\big)\in\mathcal{O}$, the sum of an invertible linear operator and a compact linear operator, and must hence be Fredholm of index zero.

For all $\big(\lambda,(\mu,w)\big)\in Q_j$ such that $\mathscr{F}\big(\lambda,(\mu,w)\big)=0$, by suitably choosing $\mathcal{Q}_j$ and by setting the left hand side of \eqref{F1_susanna} equal to zero, we can invert the linear operator $w\mapsto\mathcal{C}_{kh}(w')$ to get
\begin{equation*}
w=-(\mathcal{C}_{kh}\partial_x)^{-1}\left(\frac{1}{\mu+\lambda^2-bw}\big(J(\lambda,w)+4[w\mathcal{C}_{kh}(w')]\big) \right).
\end{equation*}
By applying Theorem~\ref{thm: hilbert reg_susanna} once more it is now easy to show that \eqref{global compactness_susanna} also holds.

The analytic global bifurcation theorem now yields the existence of a family of global continuous curves $\mathscr{C}$ of solutions to \eqref{babenko_susanna}, parameterized by $s$, with $0<s<\infty$ along which either \ref{greatest height_susanna}, \ref{parameter blowup_susanna} or \ref{eta blowup_susanna} occur, or the global curve $\mathscr{C}$ loops back and reconnects with the curve of trivial solutions.

\subsubsection{Nodal Analysis}
It remains to rule out the undesirable looping alternative. This is achieved by means of \emph{nodal analysis}. Motivated by the fact that the small-amplitude solutions are strictly monotone between consecutive crests and trough, this analysis consists of proving that this monotonicity property is conserved along the entire global curve $\mathscr{C}$. Since $\mathscr{C}$ is continuous, and hence connected, we achieve this by showing that the nodal property
\begin{equation}\label{nodal property_susanna}
	w_x<0\quad\text{for }0<x<\pi\qquad\text{and}\qquad w_x>0\quad\text{for }-\pi<x<0,
\end{equation}
defines a relatively open and closed subset of the set of nontrivial solutions to \eqref{babenko_susanna}. As in \cite{csv:critical}, the proof consists of applying maximum principle type arguments to the vertical component of the velocity field in conformal variables $V(x,y):=v\big(\xi(x,y),\eta(x,y)\big)$. Indeed, using the chain rule to compute $V$ and applying the kinematic boundary condition \eqref{kinematic top_susanna} to it, one can see that the sign of $V$ on the surface $y=0$ depends only on that of $\overline{\eta}_x$ and hence by \eqref{definition w_susanna}, on that of $w_x$. This argument relies very heavily on $v$ being harmonic in $\mathcal{R}_{kh}$.

\subsubsection{Water wave solutions}
In the last step, we must show that the solutions we have found do in fact provide us with physically relevant solutions in the fluid domain $\Omega$. Solving \eqref{babenko_susanna} for $w$ yields, via \eqref{definition w_susanna} and \eqref{eta harmonic_susanna}--\eqref{eta bottom_susanna}, $\eta$ in $\mathcal{R}_{kh}$. From the considerations in Section~\ref{subsect:conformal map_susanna}, it now suffices to check that the holomorphic function $\xi+i\eta$ is injective along the surface $y=0$. Injectivity is lost when the free surface self-intersects, scenario which, by symmetry, can only take place along the crest or the trough line. Using the Hopf boundary point lemma and the Cauchy--Riemann equations, we see that self-intersection along the crest line would violate the nodal property \eqref{nodal property_susanna}, and can hence only occur along the trough line, as stated in \ref{self-intersection_susanna}; see Figure~\ref{fig:self-intersection_susanna}.

\section{Stokes waves in constant vorticity flows (by V.~M.~Hur)}\label{sec:vera}
\subsection{Introduction}
Stokes \cite{Stokes1847, Stokes1880} made many contributions about periodic waves at the free surface of an incompressible inviscid fluid in two dimensions, under the influence of gravity, traveling at a constant velocity without change of shape. Particularly, he observed that crests become sharper and troughs flatter as the amplitude increases, and the so-called wave of greatest height ({\em extreme} wave) exhibits a $120^\circ$ corner at the crest. Most of the existing mathematical treatments of Stokes waves assumed that the flow is irrotational, whereby the stream function is harmonic inside the fluid. For instance, based on the reformulation of the problem via conformal mapping as the nonlinear pseudo-differential equation \cite{Babenko} (see also \cite{DKSZ1996} among others)
\begin{equation}\label{eqn4:Babenko0}
c^2\mathcal{C}w'-gw=g(w\mathcal{C}w'+\mathcal{C}(ww')),
\end{equation}
impressive progress was achieved analytically \cite{BDT2000a,BDT2000b} and numerically \cite{DLK2013,DLK2016,Lushnikov2016,LDS2017}. See also Section~\ref{sec:john} and references therein for more discussion.

On the other hand, vorticity has profound effects in many circumstances, for instance, for wind waves or waves in a shear flow, for which vorticity affects the wave-current interactions significantly. Stokes waves in rotational flows have had a major renewal of interest during the past two decades \cite{CS2004,Hur2006,Hur2011, CS2011,csv:critical}. See also Sections~\ref{sec:walter} and \ref{sec:susanna}, and references therein. {\em Constant vorticity} is of particular interest because the fluid flow can be written as the sum of a linear shear flow and an irrotational flow, whereby one can adapt the approaches for zero vorticity.

For instance, Simmen and Saffman \cite{SS1985} (see also \cite{TdSP1988} for finite depth) numerically solved a boundary integral equation for the problem and discovered that overhanging waves appear for a large value of positive constant vorticity, in marked contrast to the fact that the wave profile must be the graph of a single-valued function in an irrotational flow.  A  `limiting' configuration is either an extreme wave, like in an irrotational flow, or a {\em touching} wave whose profile overturns to intersect itself tangentially along the trough line, enclosing a bubble of air. See Figure~\ref{fig4:touching1} for examples. Here we distinguish positive vorticity for waves propagating upstream versus negative vorticity for downstream.

Recently, Dyachenko and Hur \cite{DH2019b} (see also \cite{DH2019c}) followed the same line of argument as in \cite{DKSZ1996} and others, to derive the equation
\begin{equation}\label{eqn4:Babenko}
c^2\mathcal{C}w'-(g+\omega c)w=g(w\mathcal{C}w'+\mathcal{C}(ww'))
+\tfrac12\omega^2(w^2+\mathcal{C}(w^2w')+w^2\mathcal{C}w'-2w\mathcal{C}(ww')),
\end{equation}
where $\omega$ denotes the constant vorticity. Throughout this section we use the notation of Sections~\ref{sec:john} and \ref{sec:susanna}. When $\omega=0$, \eqref{eqn4:Babenko} becomes \eqref{eqn4:Babenko0}. A solution of \eqref{eqn4:Babenko} gives rise to a Stokes wave in a constant vorticity flow, provided that
\begin{align}
&\text{$t\mapsto (t+\mathcal{C}w(t),w(t))$, $t\in\mathbb{R}$, is injective} \label{eqn4:touching}
\intertext{and}
&(1+\mathcal{C}w'(t))^2+w'(t)^2\neq0\quad\text{for all $t\in\mathbb{R}$}, \label{eqn4:extreme}
\end{align}
in which case the fluid surface is given parametrically as $\{(t+\mathcal{C}w(t),w(t)):t\in\mathbb{R}\}$.

We remark that numerical solutions of \eqref{eqn4:Babenko} can be found even if \eqref{eqn4:touching} fails to hold, although such solutions are `physically unrealistic' because the fluid surface intersects itself and the fluid flow becomes multi-valued. See Figures~\ref{fig4:w=2} and \ref{fig4:w=2.5} for examples. If  \eqref{eqn4:extreme} fails to hold, on the other hand, the conformal mapping from the lower-half plane to the fluid region is not well-defined at the boundary of the half plane and there is a stagnation point at the fluid surface, namely an extreme wave.

Constantin, Strauss and Varvaruca \cite{csv:critical} (see also Section~\ref{sec:susanna} for more discussion, and \cite{DH2019a,DH2019b,DH2019c} for infinite depth in the notation herein) derived
\begin{subequations}\label{eqn4:Babenko(b)}
\begin{equation}
\begin{aligned}
(c^2+b)\mathcal{C}w'-&(g+\omega c)w \\=&g(w\mathcal{C}w'+\mathcal{C}(ww'))
+\tfrac12\omega^2(w^2+\mathcal{C}(w^2w')+w^2\mathcal{C}w'-2w\mathcal{C}(ww')),
\end{aligned}
\end{equation}
where $b$ is part of the solution, subject to
\begin{equation}
[(c+\omega w(1+\mathcal{C}w')-w\mathcal{C}(ww'))^2]
=[(c^2+2b-2gw)((1+\mathcal{C}w')^2+(w')^2)],
\end{equation}
\end{subequations}
and they developed the methodology of \cite{BDT2000a,BDT2000b} further for the global bifurcation of Stokes waves in a constant vorticity flow, permitting overhanging profiles and critical layers. Also \eqref{eqn4:Babenko(b)} is equivalent to
\begin{equation}\label{eqn4:Bernoulli}
\frac{(c+\omega w(1+\mathcal{C}w')-\omega\mathcal{C}(ww'))^2}{(1+\mathcal{C}w')^2+(w')^2}=c^2+2b-2gy.
\end{equation}
See \cite{DH2019a,DH2019b,DH2019c} for more discussion.

\subsection{Numerical method}
The linearized operator of \eqref{eqn4:Babenko} is self-adjoint, whereby one can employ, for instance, the conjugate gradient (CG) method for numerical treatments \cite{DH2019b}. On the other hand, the linearized operator of \eqref{eqn4:Babenko(b)} is not self-adjoint, whence the CG method is no longer applicable. One can use, for instance, the GMRES method instead \cite{DH2019a}, but the method turns out to be impracticable, if not futile, for large values of vorticity for large amplitude. See \cite{DH2019b} for more discussion.

The linearized operator of \eqref{eqn4:Babenko} is not positive definite, though, and the CG method turns out to diverge for `almost' extreme waves. The method converges for any value of vorticity for sufficiently large amplitude nonetheless. One can resort to other Krylov subspace methods, particularly, MINRES methods. In fact, replacing the CG method by the conjugate residual method works well. For instance, Dyachenko, Hur and Silantyev \cite{DHS2021} required the residuals $\lesssim10^{-43}$ and the Fourier coefficients $\lesssim 10^{-38}$ to approximate more than $200$ decimal digits of extreme waves, after concentrating grid points near the crest via an auxiliary conformal mapping. Alternatively, one can solve the linearized equation of \eqref{eqn4:Babenko} directly. For instance, for $\omega=0$, a $4096\times4096$ linear system suffices to approximate the extreme wave with high precision. We refer the reader to \cite{DH2019a,DH2019b,DH2019c,DHS2021} for details of how to solve \eqref{eqn4:Babenko} numerically.

\subsection{Folds and gaps}
When $\omega=0$, Longuet-Higgins and Fox \cite{LHF1977a, LHF1977b}, among others, predicted that the wave speed oscillates infinitely many times as the steepness increases monotonically toward the extreme wave. Here the steepness $s$, say, measures the crest-to-trough vertical distance divided by the period. See also Section~\ref{sec:john} for more discussion. In fact, Dyachenko, Hur and Silantyev \cite{DHS2021} resolved at least six oscillations of the wave speed numerically, approximating the steepness of the extreme wave with high precision. For instance, for the Nekrasov parameter \cite{Nekrasov} $\nu\approx2.2\times 10^{-26}$ (see Section~\ref{sec:john}, $\nu=0$ for the extreme wave), $s\approx 0.141063483979936081$.

\begin{figure}
\centerline{
\includegraphics{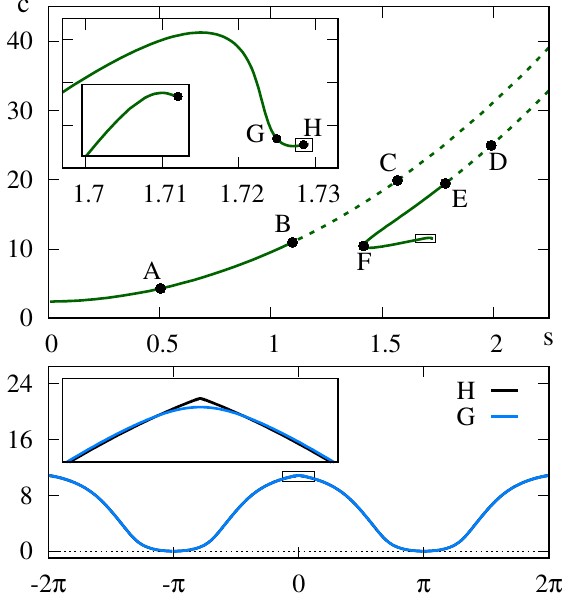}
\includegraphics{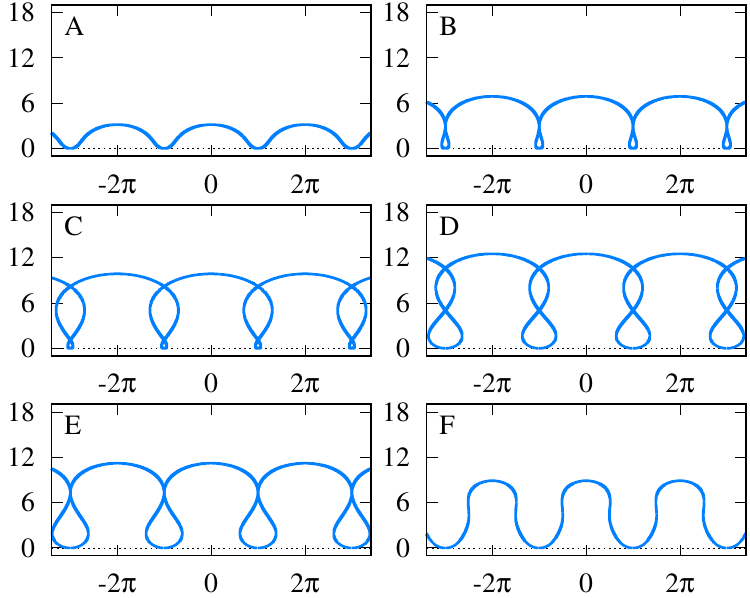}}
\caption{$\omega=2$: (upper left) Wave speed vs. steepness. Solid curves are for physical solutions and dashed curves for unphysical. The insets are closeups near the extreme wave, the endpoint of numerical continuation. (right) Several wave profiles in the physical space, labelled by $A$ through $H$ along the solution curve. (lower left) Almost extreme waves. The inset is a closeup near the crests.}
\label{fig4:w=2}
\end{figure}

For a large value of positive constant vorticity, for instance for $\omega=2$, on the other hand, the upper left panel of Figure~\ref{fig4:w=2} shows the wave speed versus the steepness for numerical solutions of \eqref{eqn4:Babenko} \cite{DH2019b}, illustrating vast differences from zero vorticity. In the right and lower left panels are the free surfaces in the physical space at the indicated points along the $c=c(s)$ curve.

To begin, $s$ increases and then decreases from $s=0$ to wave $F$, namely a {\em fold} in the $c=c(s)$ curve. Consequently, there is more than one solution for some values of $s$. When $s$ is small, for instance for wave~$A$, the profile is not overhanging. As $s$ increases along the fold, on the other hand, the wave profile becomes more rounded, so that overhanging waves appear and ultimately a touching wave, whose profile intersects itself tangentially along the trough line, enclosing a bubble of air. For instance, wave~$B$ is an almost touching wave. See also Figure~\ref{fig4:touching1}. Past such a touching wave, numerical solutions are unphysical because \eqref{eqn4:touching} no longer holds true. For instance, for wave~$C$, the fluid surface intersects itself transversely along the trough line and the fluid flow becomes multi-valued. As $s$ decreases along the fold, the wave profile becomes less rounded, so that another touching wave appears, past which numerical solutions become physical. For instance, wave~$D$ is unphysical, wave~$E$ is an almost touching wave, and wave~$F$ is physical. Therefore, there is a {\em gap} in the $c=c(s)$ curve, consisting of unphysical solutions, bounded by two touching waves, waves~$B$ and $E$. By the way, wave~$E$ encloses a larger bubble of air than wave~$B$.

\begin{figure}
\centerline{
\includegraphics{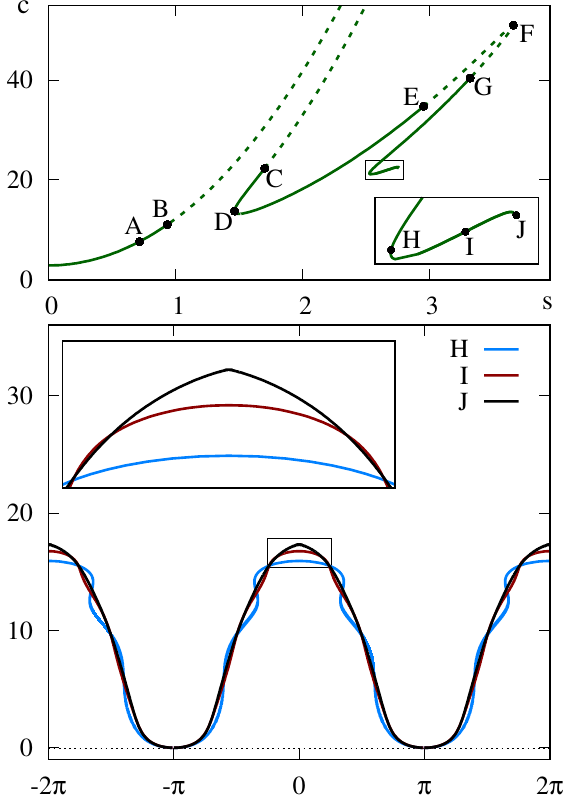}
\includegraphics{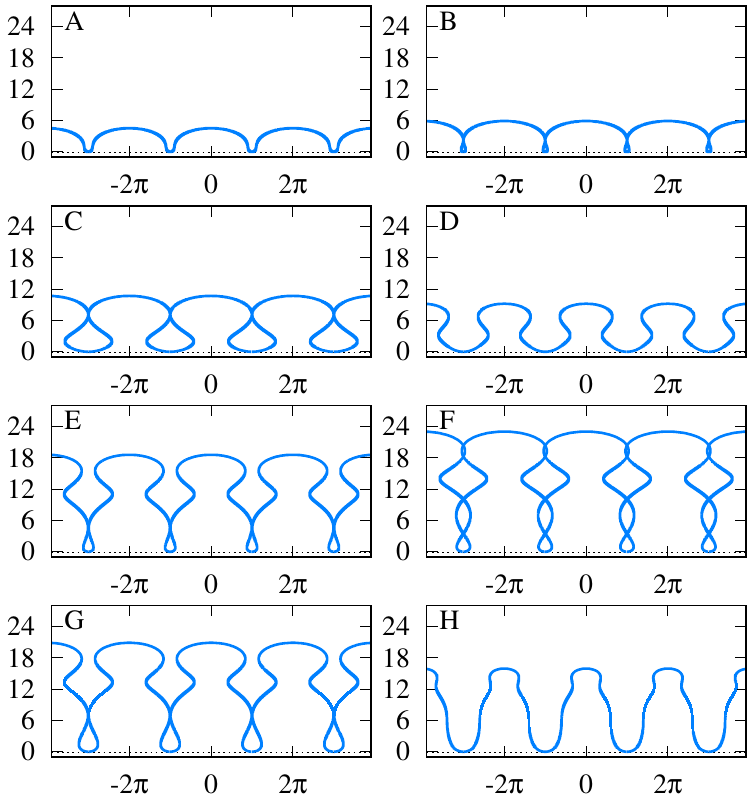}}
\caption{\small $\omega=2.5$: (upper left) Wave speed vs. steepness. Solid and dashed curves are for physical and unphysical solutions. The inset is a closeup past the second fold. (right) Several wave profiles in the physical space, labelled by A through H along the solution curve. (lower left) Almost extreme waves. The inset is a closeup near the crests.}
\label{fig4:w=2.5}
\end{figure}

We point out that the numerical method of \cite{SS1985,sp:steep} and others, based on a boundary integral equation for the problem, diverges in a gap. Based on \eqref{eqn4:Babenko}, involving the periodic Hilbert transform, on the other hand, the numerical method of \cite{DH2019a,DH2019b,DH2019c} converges throughout and explains what is going on in the gap.

The upper left panel of Figure~\ref{fig4:w=2} shows that beyond wave~$F$ at the endpoint of the fold, $s$ increases monotonically whereas $c$ oscillates, like in an irrotational flow. Also the lower left panel suggests that as $s$ increases past the fold, the crest becomes sharper and the wave profile is no longer overhanging, as in an irrotational flow. In fact, numerical evidence \cite{DH2019b} supports that the solution curve is ultimately limited by an extreme wave whose profiles have a sharp corner at the crest. For instance, wave~$H$ is an almost extreme wave. One may not be able to continue numerical solutions past such an extreme wave where \eqref{eqn4:extreme} no longer holds true.

We point out that this extreme wave is not the wave of greatest height. In fact, for wave~$E$, $s=1.78406376$ and for the extreme wave, $s\approx 1.72849105$.

For a larger value of positive constant vorticity, for instance, for $\omega=2.5$, the upper left panel of Figure~\ref{fig4:w=2.5} shows a fold from $s=0$ to wave~$D$, and a gap bounded by waves~$B$ and $C$, like Figure~\ref{fig4:w=2}. An important difference is another fold from wave~$D$ to $H$. As $s$ increases along the second fold, the wave profile becomes more rounded, so that a touching wave, wave~$E$, appears, like along the lowest fold. Past the touching wave, numerical solutions are unphysical. See, for instance, wave~$F$. As $s$ decreases along the second fold, the wave profile becomes less rounded, so that another touching wave appears, wave~$G$, past which numerical solutions are physical, like along the lowest fold. See, for instance, wave~$H$. Together, there is another gap bounded by two touching waves, waves~$E$ and $G$. See \cite{DH2019b,DH2019c} for more discussion.

Recall that the numerical method of \cite{SS1985,sp:steep} and others diverges in a gap, and in order to seek physical solutions past the gap, the authors solved the boundary integral equation assuming that the solution would have a stagnation point at the crest, in order to approximate an almost extreme wave numerically.  They continued the solution until an almost touching wave appears, past which numerical solutions would become unphysical. Such a strategy would work when there is only one gap but  it is incapable of finding a second gap, which does exist.

The upper left panel of Figure~\ref{fig4:w=2.5} shows that past the second fold, $s$ increases monotonically although $c$ oscillates and the crest becomes sharper, so that one expects that the solution curve is ultimately limited by an extreme wave, like Figure~\ref{fig4:w=2}. For instance, wave~$J$ is an extreme wave whose profile is not overhanging. See \cite{DH2019b,DH2019c} for more discussion.

\begin{figure}
\centerline{
\includegraphics[scale=1.0]{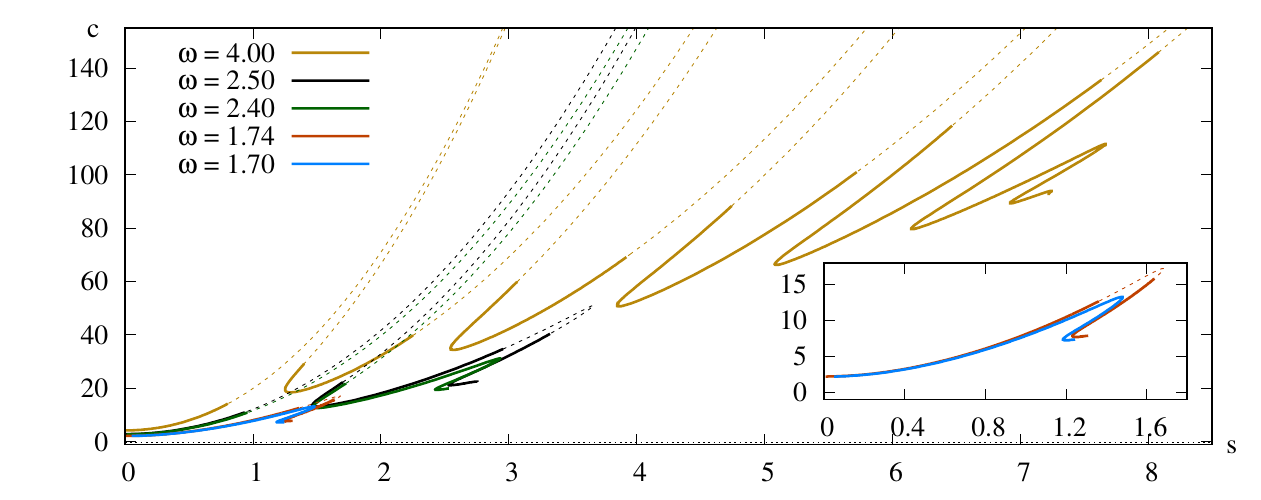}}
\caption{Wave speed vs. steepness for five values of positive constant vorticity. Solid and dashed curves are for physical and unphysical solutions. The inset is a closeup for $\omega=1.7$ and $1.74$.}
\label{fig4:c=c(s)}
\end{figure}
Figure~\ref{fig4:c=c(s)} shows the wave speed versus steepness for several values of positive constant vorticity. For $\omega=0$ (not shown because the solution curve is much smaller than the others), there is analytical and numerical evidence \cite{LHF1977a,LHF1977b} that $c$ oscillates infinitely many times whereas $s$ increases monotonically. For negative constant vorticity, crests become sharper and lower. See \cite{DH2019a,DH2019c} for more discussion. For instance, for $\omega=1.7$, on the other hand, $s$ no longer increases monotonically but, rather, the lowest oscillation of $c$ turns into a fold. 
The fold becomes larger as $\omega$ increases so that, for instance for $\omega=1.74$, part of the fold turns into a gap of unphysical solutions. Moreover, for instance for $\omega=2.4$, the second oscillation of $c$ turns into another fold, and the fold becomes larger as $\omega$ increases so that, for instance for $\omega=2.5$, part of the second fold turns into another gap.

More folds and gaps appear in a similar manner as $\omega$ increases. In fact, when $\omega=4$, Figure~\ref{fig4:c=c(s)} shows five folds and five gaps. There is numerical evidence \cite{DH2019a,DH2019b,DH2019c} that past all the folds, $s$ increases monotonically toward an extreme wave whose profile has a sharp corner at the crest. Also there is analytical evidence \cite{DHS2021} that the extreme wave has a $120^\circ$ corner at the crest, regardless of the value of the vorticity.

\subsection{Almost extreme waves}

In what follows, by the {\em angle} we mean the angle that the fluid surface of a Stokes wave makes with the horizontal.

\begin{figure}
\centerline{
\includegraphics[scale=0.9]{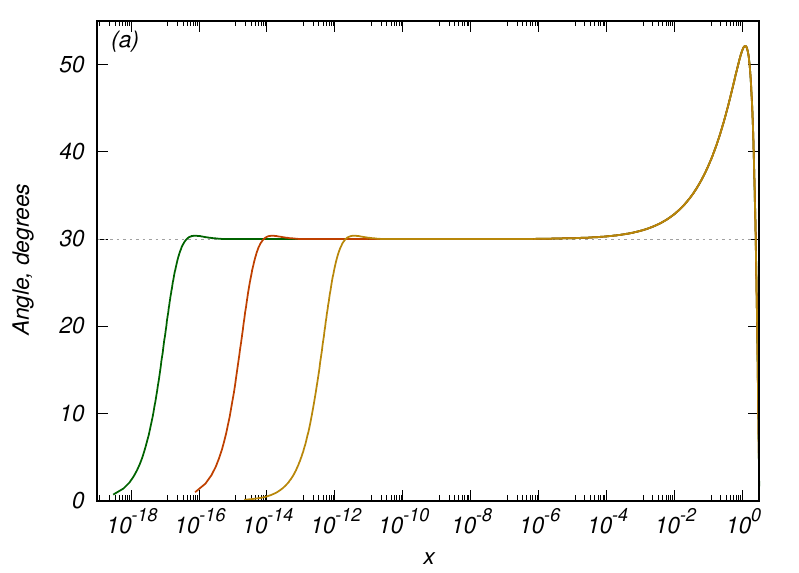}
\includegraphics[scale=0.9]{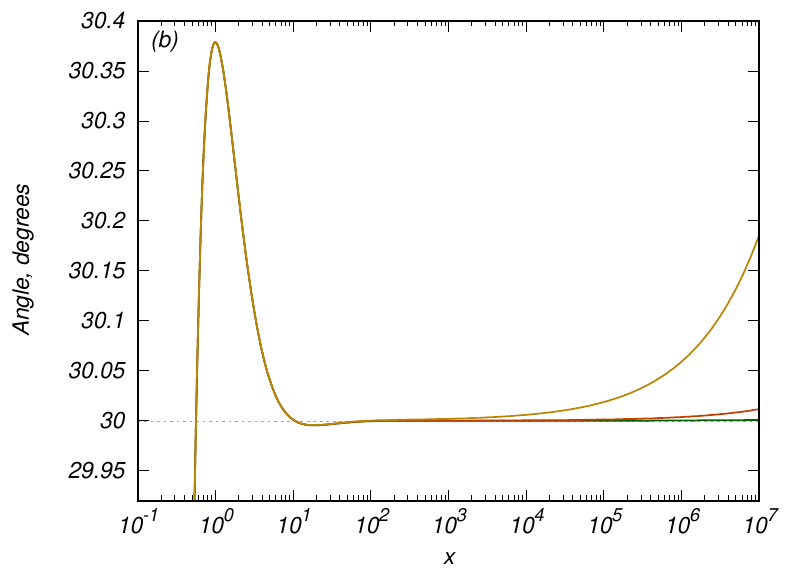}}
\caption{$\omega=1$: (left) Angle vs. $\log(x)$ for almost extreme waves. (right) A closeup near the first local maximum angle.}
\label{fig4:angle}
\end{figure}

When $\omega=0$, the maximum angle of the extreme wave $=30^\circ$ at the crest \cite{AFT1982,Plotnikov1982,PT2004}. So it came as a surprise when McLeod \cite{McLeod} proved that the angle of an almost extreme wave $>30^\circ$ near the crest. Longuet-Higgins and Fox \cite{LHF1977a,LHF1977b} combined asymptotic analysis and numerical computation to predict the maximum angle of an almost extreme wave $\approx 30.37^\circ$. See also Section~\ref{sec:john} for more discussion. Chandler and Graham \cite{CG1993} solved Nekrasov's equation \cite{Nekrasov} (see also Section~\ref{sec:john}) numerically to find that:
\begin{itemize}
\item there is a `boundary layer' where the angle increases sharply from $0^\circ$ at the crest to the maximum~$\approx 30.3787032466^\circ$ (see also Figure~\ref{fig4:angle}); this agrees well with the result of \cite{LHF1977a,LHF1977b};
\item past the boundary layer, there is a region where the angle oscillates about $30^\circ$, resembling the Gibbs phenomenon (see also Figure~\ref{fig4:angle}); the authors resolved two and a half oscillations; and
\item past the oscillations, the angle decreases monotonically to $0^\circ$ at the trough.
\end{itemize}
Recently, Dyachenko, Hur and Silantyev \cite{DHS2021} improved the result, calculating the maximum angle in the boundary layer $\approx 30.3787032466519114^\circ$ and resolving at least three and a half oscillations of the angle with high precision, and they took the matters further to nonzero constant vorticity.

For instance for $\omega=1$, Figure~\ref{fig4:angle} shows the angles of almost extreme waves versus the horizontal coordinate in the logarithmic scale for numerical solutions of \eqref{eqn4:Babenko}. The steepness $\approx 0.443104987824811270$. This is noticeably higher than the steepness of the extreme wave for $\omega=0$. Also the extreme steepness is $\approx 0.04929917508893317$ for $\omega=-1$. In fact, there is numerical evidence \cite{DH2019b,DHS2021} that the extreme steepness increases monotonically as $\omega$ increases.

The left panel of Figure~\ref{fig4:angle} shows a boundary layer where the angle increases sharply from $0^\circ$ at the crest to a local maximum, like in an irrotational flow, where the maximum angle $\approx 30.3787035187621690^\circ$. This is very close to the maximum angle for $\omega=0$. Also the maximum angle $\approx 30.3787021268445382^\circ$ for $\omega=-1$. In fact, numerical evidence supports the first local maximum angle being independent of the value of the vorticity. It would be interesting to give a rigorous proof of this. The right panel shows the oscillations of the angle about $30^\circ$, like when $\omega=0$. The first local minimum angle $\approx 29.995397640954266194^\circ$, which is very close to the first local minimum angle $\approx 29.995396467370841183^\circ$  for $\omega=0$. Numerical evidence clarifies that the first local minimum angle is also independent of the value of the vorticity. See \cite{DHS2021} for more discussion.

But an important difference is that past the oscillations the angle increases toward the global maximum $\approx 52.1426155193^\circ$ and then decreases sharply to $0^\circ$ at the trough. There is numerical evidence \cite{DH2019b,DH2019c} that the global maximum angle increases as $\omega$ increases but remains $<90^\circ$. That means, an extreme wave profile is not overhanging. It would be interesting to give a rigorous proof of this.

\subsection{Touching waves}

\begin{figure}
\centerline{
\includegraphics[scale=1]{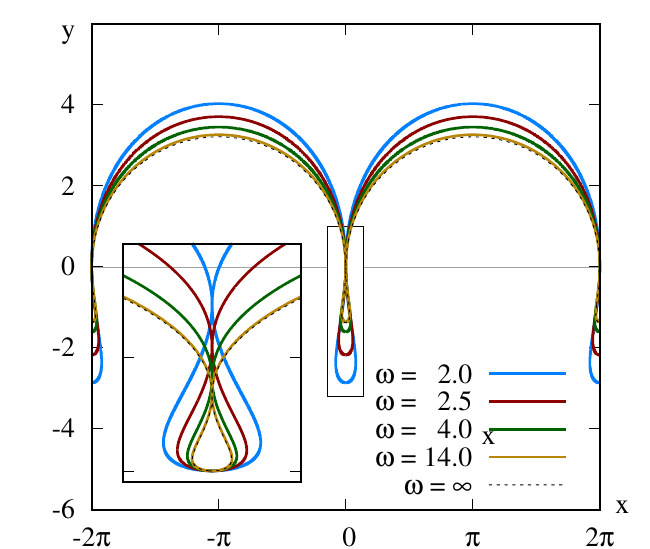}
\includegraphics[scale=1]{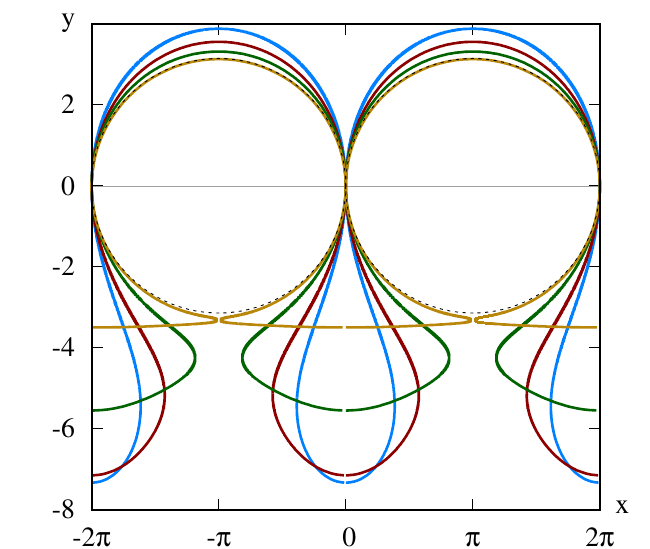}}
\caption{\small (left) Touching waves at the beginnings of the lowest gaps. The dashed curve is the limiting Crapper wave. The inset is a closeup near the troughs. (right) At the ends of the gaps. The dashed curve encloses a circular fluid region.}
\label{fig4:touching1}
\end{figure}

The left panel of Figure~\ref{fig4:touching1} shows the profiles of almost touching waves near the beginnings of the lowest gaps, and the right panel near the ends of the gaps, for four values of positive constant vorticity. In the left panel, at the beginnings of the lowest gaps, $s$ decreases monotonically toward $\approx 0.7297$ as $\omega\to \infty$. In the right panel, $s\to1$ as $\omega\to\infty$. See \cite{DH2019a} for details.

Crapper \cite{Crapper} produced an exact solution for capillary waves---nonzero surface tension and zero gravity---in an irrotational flow, whereby he deduced that crests become flatter and troughs more curved as the amplitude increases, which is opposite to gravity waves, toward a touching wave whose steepness is $0.7297\dots$. This fortuitous coincidence is no accidence. In fact, numerical evidence \cite{DH2019a,DH2019b,DH2019c} is clear that touching waves at the beginnings of the lowest gaps approach the limiting Crapper wave as the value of positive constant vorticity increases unboundedly, or equivalently, gravitational acceleration vanishes----which is an unexpected and remarkable connection between rotational and capillary effects.

Recently, Hur and Vanden-Broeck \cite{HVB2020} offered analytical and numerical evidence that any Crapper wave, not necessarily the limiting form, gives the profile of a periodic traveling wave in a constant vorticity flow for zero gravitational acceleration. Also Hur and Wheeler \cite{HW2020} proved that Crapper's formula
\begin{equation}\label{eqn4:Crapper}
(x+iy)(t)=t-\frac{4iAe^{-it}}{1+Ae^{-it}},
\end{equation}
in the complex form, makes an exact solution of \eqref{eqn4:Bernoulli}, where $g=0$, for appropriate values of $\omega$, $c$ and $b$. The wave profile of \eqref{eqn4:Crapper} is not overhanging as long as $A<\sqrt{2}-1$, and it does not intersect itself so long as $A<A_{\max}\approx 0.4546700164520109$. When $A>A_{\max}$ the profile intersects itself transversely at two points along the trough line and the fluid flow becomes multi-valued, whence it makes an unphysical solution. See \cite{HW2021} for more discussion. More recently, Hur and Wheeler \cite{HW2021} used the implicit function theorem to construct overhanging and touching waves for small values of the gravity constant. Therefore, overhanging and touching waves do exist in constant vorticity flows under gravity. We remark that there is a global bifurcation result \cite{csv:critical} which permits overhanging profiles, but it is incapable of determining whether overhanging profiles actually exist.

The right panel of Figure~\ref{fig4:touching1} suggests that touching waves at the ends of the lowest gaps approach a fluid disk in rigid body rotation as the value of the vorticity increases unboundedly or, equivalently, gravitational acceleration vanishes. See \cite{DH2019a,DH2019b,DH2019c} for more discussion. Teles da Silva and Peregrine \cite{sp:steep}, among others, computed periodic waves in constant vorticity for zero gravity and argued that a limiting configuration is such a `circular vortex wave'.

There is numerical evidence \cite{DH2019b,DH2019c} that touching waves at the beginnings of the second gaps approach the circular vortex wave on top of the limiting Crapper wave as the value of the vorticity increases unboundedly, whereas the circular vortex wave on top of itself at the ends of the gaps. Touching waves at the boundaries of higher gaps accommodate additional circular vortex waves in like manner. For instance, see \cite{DH2019b} for a profile nearly enclosing five fluid disks.

\section{Solitary waves and fronts (by M. H. Wheeler)}\label{sec:miles}

 \subsection{Introduction}\label{sec:miles:intro}
In this section we turn our attention to solitary water waves, that is,
traveling waves whose surfaces approach some asymptotic height at
infinity; see Figure~\ref{fig:miles:basic}(a). In many respects, the
theory for solitary waves is more difficult and more subtle than that
for periodic waves considered in the preceding
Sections~\ref{sec:john}--\ref{sec:vera}. Yet in other ways the problem
is much simpler, and indeed many results for solitary waves are
stronger than their periodic counterparts.

Rather than attempt a comprehensive survey of this broad topic, we
will focus on the most classical case of two-dimensional irrotational
waves, without surface tension or density stratification. Stratified
solitary waves and fronts will be discussed in Section~\ref{sec:sam:strat},
and three-dimensional solitary waves with surface tension will be
discussed in Section~\ref{sec:erik}. This still leaves out, among
other things, the theory for non-stratified waves with vorticity, as
well as the extremely rich theory of small-amplitude two-dimensional
solitary waves with surface tension. For an introduction to the
literature for the latter problem, we refer the reader to the surveys
\cite{di:handbook,groves:survey}.

After defining some basic terminology and notation in
Section~\ref{sec:miles:definitions}, in Section~\ref{sec:miles:small}
we discuss the existence of small-amplitude solitary
waves. The linear theory is significantly less informative than in the
periodic case, and instead we must use weakly nonlinear
models such as the Korteweg--de Vries equation for intuition. We
outline how this intuition can be made into a rigorous proof
using ``spatial dynamics'' and center manifold reduction techniques.
Section~\ref{sec:miles:qual} collects a series of
results which together provide a surprisingly detailed qualitative
picture of an arbitrary solitary wave, independent of its
amplitude or construction. In Section~\ref{sec:miles:nobores}, we
explain why front-type solutions or \emph{bores}, where the free surface
limits to different values as $x \to \pm\infty$
(Figure~\ref{fig:miles:basic}(b)), cannot exist in this
classical model. Perhaps surprisingly, this nonexistence result has
important implications for the existence of large-amplitude solitary
waves, which is the topic of Section~\ref{sec:miles:large}. Finally,
Section~\ref{sec:miles:open} lists some open problems, many of
which appear to be quite difficult.

\subsection{Basic definitions}\label{sec:miles:definitions}

\begin{figure}
  \centering
  \includegraphics[page=1]{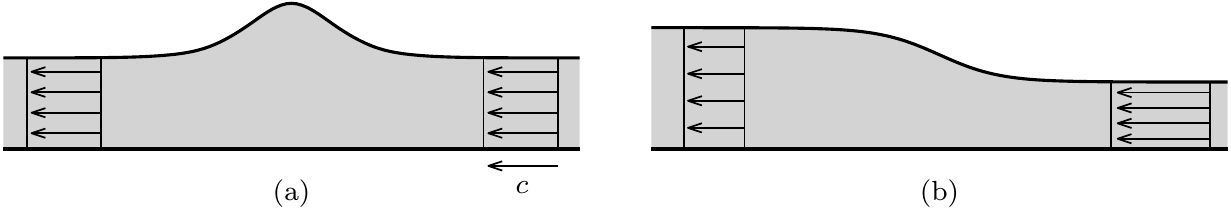}
  \caption{(a) A solitary wave. The arrows in the far field represent
  the asymptotic relative fluid velocity $(u-c,v) \to (-c,0)$.
  (b) A bore, with distinct asymptotic velocities and depths at $x =
  \pm\infty$. }
  \label{fig:miles:basic}
\end{figure}

The stream function formulation \eqref{stream eq} of the two-dimensional
steady water wave problem, specialized to irrotational waves with
$\omega\equiv 0$, reads
\begin{equation}\label{eqn:miles:stream}
  \left\{
    \begin{alignedat}{2}
      \Delta \psi &= 0 &\qquad& \text{ in } {-d} < y < \eta(x),\\
      \psi &= 0 &\qquad& \text { on }  y = \eta(x),\\
      \psi &= -m &\qquad& \text { on }  y = -d,\\
      \tfrac 12 |\nabla\psi|^2 + g\eta &= \tfrac 12 Q &\qquad& \text { on }  y = \eta(x).
    \end{alignedat}
  \right.
\end{equation}
For any wave speed $c > 0$ and depth $d > 0$, \eqref{eqn:miles:stream} has
the explicit ``trivial'' solution
\begin{align}
  \label{eqn:miles:trivial}
  \psi = -cy,
  \qquad
  \eta \equiv 0,
  \qquad
  Q = c^2,
  \qquad
  m = cd.
\end{align}
corresponding to body of a fluid at rest with a flat free surface. In
a reference frame moving with the wave, the fluid instead has
uniform horizontal velocity $-c$. By a \emph{solitary} wave, we mean a
solution to \eqref{eqn:miles:stream} which is non-constant in $x$ but
converges to \eqref{eqn:miles:trivial} as $x \to \pm\infty$ in the
sense that
\begin{align}
  \label{eqn:miles:asym}
  \psi \to -cy,
  \quad \eta \to 0
  \qquad
  \text{ as }
  x \to \pm\infty.
\end{align}
See Figure~\ref{fig:miles:basic}(a).
A \emph{front} or \emph{bore}, by contrast, is a
solution where $\psi$ and $\eta$ have different limits as $x \to
+\infty$ and $x \to -\infty$; see Figure~\ref{fig:miles:basic}(b).
As we will show in Section~\ref{sec:miles:nobores}, such solutions to
\eqref{eqn:miles:stream} do not exist, but they do exist in the
stratified setting considered in Section~\ref{sec:sam:strat}.

\begin{figure}
  \centering
  \includegraphics[page=2]{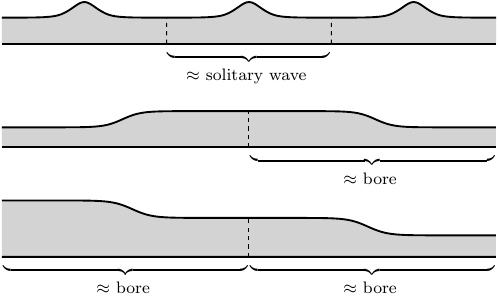}
  \caption{One can imagine a solitary waves arising as broadening limits of
  periodic waves, and bores arising as broadening
  limits of solitary waves or indeed other bores.}
  \label{fig:miles:broadening}
\end{figure}

In addition to their own independent interest, solitary waves can
occur as broadening limits of periodic waves with increasingly wide
and flat troughs; see Figure~\ref{fig:miles:broadening}. Similarly,
bores might occur as broadening limits of solitary waves, or indeed as broadening
limits of other fonts. It is by no means clear, however, that
\emph{all} solitary waves or bores can be realized in this way.

Taking $x \to \pm\infty$ in
\eqref{eqn:miles:stream} and \eqref{eqn:miles:asym}, we discover that
solitary waves necessarily have $Q=c^2$ and $m=cd$. Thus there are
only three dimensional parameters: $c,g,d$. Because we are free to
choose both a length scale and a velocity scale, there is in fact just
a single dimensionless parameter
\begin{align}
  \label{eqn:miles:froude}
  F = \frac c{\sqrt{gd}}
\end{align}
called the \emph{Froude number}. As we will see in the next
subsection, $F=1$ is a critical threshold for the linearized problem.
Waves with $F>1$ are traditionally called \emph{supercritical}, while
those with $F<1$ are called \emph{subcritical}.

\subsection{Existence of small-amplitude solitary waves}\label{sec:miles:small}
\begin{figure}
  \centering
  \includegraphics[page=3]{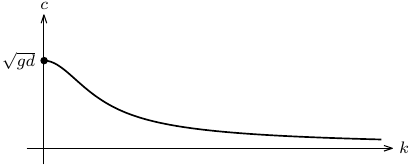}
  \caption{The dispersion relation \eqref{eqn:miles:dispersion}.}
  \label{fig:miles:dispersion}
\end{figure}

Formally linearizing \eqref{eqn:miles:stream} about the trivial
solution \eqref{eqn:miles:trivial} leads to the problem
\begin{equation}\label{eqn:miles:linear}
  \left\{
    \begin{alignedat}{2}
      \Delta \dot\psi &= 0 &\qquad& \text{ in } {-d} < y < 0,\\
      \dot\psi - c\dot\eta &= 0 &\qquad& \text { on }  y = 0,\\
      \dot\psi &= 0 &\qquad& \text { on }  y = -d,\\
      -c\dot \psi_y + g \dot \eta  &= 0 &\qquad& \text { on }  y = 0,
    \end{alignedat}
  \right.
\end{equation}
where here dots denote linearized variables. Note that, unlike the
original system \eqref{eqn:miles:stream}, \eqref{eqn:miles:linear} is
posed on a fixed domain. Separating variables, we discover that
\eqref{eqn:miles:linear} has a solution with $\dot \eta = e^{ikx}$ if
and only if the \emph{dispersion relation}
\begin{align}
  \label{eqn:miles:dispersion}
  c^2 = g \frac{\tanh kd}k
\end{align}
holds. Looking at the graph in Figure~\ref{fig:miles:dispersion}, we
see that periodic solutions with real wavenumbers $k$ exist if and only
if the Froude number $F=c/\sqrt{gd} < 1$ is subcritical.
These are indeed first-order approximations of periodic solutions to
the full problem \cite{Nekrasov,lc,struik}, a phenomenon that can now
be understood using the local bifurcation theory discussed in
Sections~\ref{sec:john} and \ref{sec:susanna}.

Such arguments \emph{fail completely} for solitary waves! Thinking of
periodic waves in the $k \to 0$ limit, and glancing again at
Figure~\ref{fig:miles:dispersion}, we might expect solitary waves to
bifurcate from the critical Froude number $F=c/\sqrt{gd}=1$. But for
this value of $F$ the only bounded solutions of
\eqref{eqn:miles:linear} are constant in $x$, not localized. From a
more rigorous point of view, the linearized operator on the left hand
of \eqref{eqn:miles:linear} is non-Fredholm at $F=1$, and so local
bifurcation techniques cannot be applied.

Indeed, Russell's initial observations of solitary waves \cite{russell}
were famously controversial at the time precisely because they did
not fit into the model provided by \eqref{eqn:miles:linear}. This
issue was not resolved until decades later, when Boussinesq and
Rayleigh independently developed weakly nonlinear models for long
waves. One of these models was later rediscovered by Korteweg and
de Vries, and is now referred to as the KdV equation. We refer the
reader to \cite{miles:survey,darrigol:horse,craik} for more on this
interesting history.

Unlike \eqref{eqn:miles:linear}, the KdV equation is nonlinear,
derived by carrying out an asymptotic expansion to second order.
Crucially, it is also a long-wave model, meaning that the expansion
involves rescaling the horizontal variable $x$. The free surface
profile $\eta$, for instance, is expanded as
\begin{align}
  \label{eqn:miles:kdvexp}
  \eta(x)
  = \varepsilon \eta_1(\sqrt{\varepsilon} x)
  + \varepsilon^2 \eta_2(\sqrt{\varepsilon} x)
  + \cdots,
\end{align}
where here $\varepsilon > 0$ is a small parameter roughly
corresponding to the amplitude of the wave.
Eventually one obtains the model equation
\begin{align}
  \label{eqn:miles:kdv}
  \eta'' + \frac3{F^2d^2} (1-F^2) \eta + \frac 9{2F^2d^3}  \eta^2 &=
  0,
\end{align}
which is a dimensional and steady version of the time-dependent equation
\eqref{KdV} given in Section~\ref{sec:walter}.
\begin{figure}
  \centering
  \includegraphics[page=4]{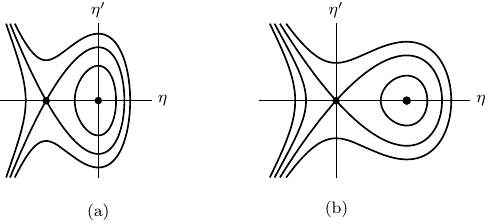}
  \caption{Phase portraits for \eqref{eqn:miles:kdv} when (a) $F<1$
  and (b) $F>1$.}
  \label{fig:miles:phase}
\end{figure}
The corresponding phase portrait is sketched in
Figure~\ref{fig:miles:phase}. When $F < 1$ is subcritical, the origin
is a center, and all small solutions are periodic. When $F>1$ is
supercritical, on the other hand, there is an explicit holomoclinic
orbit
\begin{align}
  \label{eqn:miles:sech2}
  \eta
  = d(F^2-1)  \operatorname{sech}^2
  \Big( \frac{\sqrt{3}}2 \frac{\sqrt{F^2-1}}{Fd} x \Big).
\end{align}
corresponding to a solitary wave. Comparing with
\eqref{eqn:miles:kdvexp}, we expect \eqref{eqn:miles:sech2} to be a
valid approximation only for small waves $F \approx 1$.

Now that we have some formal intuition about the existence of
small-amplitude solitary waves, we turn to the rigorous theory. As we
have seen, local bifurcation arguments based on an analysis of the
linearized problem \eqref{eqn:miles:linear} are doomed to failure, and
so a different strategy is required. The first proof is due to
Lavrentiev~\cite{lavrentiev}, who constructed solitary waves as
long-wavelength limits of (weakly nonlinear) periodic waves.
Friedrichs and Hyers subsequently gave a simpler
argument~\cite{fh:existence}.
By rescaling the horizontal variable $x$ as in
\eqref{eqn:miles:kdvexp}, subtracting off the approximate solution
coming from \eqref{eqn:miles:sech2}, and working in carefully chosen
exponentially weighted spaces, they are essentially able apply the
implicit function theorem. Later Beale gave another argument based on
the Nash--Moser implicit function theorem~\cite{beale:existence}, and
Mielke~\cite{mielke:reduction} gave a proof using dynamical systems
techniques. Perhaps the simplest proof in the literature is that of
Pego and Sun in \cite[Appendix~A]{ps:stability}, where they give a
more modern fixed-point argument in the spirit of \cite{fh:existence}.

The center manifold techniques employed by
Mielke~\cite{mielke:reduction} arguably lead to the strongest results.
Roughly speaking, Mielke shows that when $F \approx 1$, all
small-amplitude solutions of the full water wave problem
\eqref{eqn:miles:stream} solve an ODE
\begin{align}
  \label{eqn:miles:mielke}
  \eta'' + \frac3{F^2d^2} (1-F^2) \eta + \frac 9{2F^2d^3}  \eta^2 =
  f(\eta,\eta',F^2-1)
\end{align}
\emph{exactly}, where here the left hand side is nothing other than
the left hand side of \eqref{eqn:miles:kdv}. While the function $f$ on right
hand side is not explicit, it is higher-order (after an appropriate
rescaling),  has the expected symmetries, and can in principle be Taylor
expanded to any finite order. The reduction of the PDE
\eqref{eqn:miles:stream} to the ODE \eqref{eqn:miles:mielke} is truly
remarkable --- even if it is restricted to a certain perturbative regime ---
and hugely useful.
The existence of solitary waves, for instance, now becomes
a question about the persistence of the homoclinic orbit
\eqref{eqn:miles:sech2} to \eqref{eqn:miles:kdv} when the additional
nonlinear terms on the right hand side of \eqref{eqn:miles:mielke} are
included. Applying phase portrait arguments to
\eqref{eqn:miles:mielke} also yields strong uniqueness results, as
well as qualitative information, e.g.~about the monotonicity of solutions.

Center manifold techniques have been extremely successful in the study
of small-amplitude water waves of many kinds, both in two and three
dimensions; see Section~\ref{sec:erik}, as well as, for instance,
\cite{di:handbook,hi:book}. Interestingly, these techniques are not
restricted to local equations such as \eqref{eqn:miles:stream}, but
also apply to certain non-local problems \cite{fs:nophase}, including
some non-local models for traveling water waves \cite{jww:whitham,
jtw:capgrav}.

\subsection{Qualitative properties of solitary waves}\label{sec:miles:qual}

The explicit solution \eqref{eqn:miles:sech2} to \eqref{eqn:miles:kdv}
has several basic properties:
\begin{enumerate}
\item It has a supercritical Froude number $F>1$.
\item It is a \emph{wave of elevation} with $\eta > 0$.
  In other words, as in Figure~\ref{fig:miles:basic}(a), the free
  surface lies everywhere above its asymptotic level.
\item It is \emph{symmetric} and \emph{monotone} in that,
  after a translation, $\eta$ is even in $x$ with $\eta_x < 0$
  for $x > 0$. In particular, there is exactly one crest.
\end{enumerate}
By analyzing the reduced ODE \eqref{eqn:miles:mielke}, one can show
that all small-amplitude solitary waves also satisfy these properties.
In this section we will see that, in fact, these properties hold for
\emph{any} solitary wave, regardless of its amplitude or construction.

\begin{figure}
  \centering
  \includegraphics[page=5]{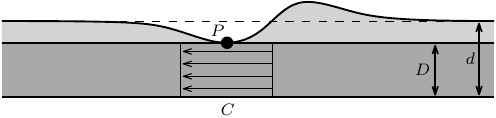}
  \caption{The auxiliary uniform flow in the proof of
  Theorem~\ref{thm:miles:elevation}.}
  \label{fig:miles:comp}
\end{figure}
\begin{theorem}[\cite{cs:sym,craig:nonexistence}]\label{thm:miles:elevation}
  All supercritical solitary waves are waves of elevation
  \begin{proof}
    We sketch the elegant argument given in \cite{craig:nonexistence},
    based on the maximum principle. Supposing that $\min \eta \le 0$
    is achieved,
    we consider the auxiliary uniform flow shown in
    Figure~\ref{fig:miles:comp}. The depth $D$ of this flow is chosen
    so that the two fluid surfaces touch tangentially at the point $P$
    where $\eta$ achieves its minimum, and the speed $C$ is chosen
    so that the mass flux $CD$ of the auxiliary flow matches the flux
    $cd$ of the solitary wave. Letting $\Psi$ denote the stream
    function for the auxiliary flow, we then argue that the
    difference $\psi-\Psi$ achieves its minimum value of $P$. Applying
    the Hopf lemma, and using the various boundary conditions, we
    discover that the Froude number $F < 1$.
  \end{proof}
\end{theorem}

The converse to Theorem~\ref{thm:miles:elevation} is also true. The
proof is based on an integral identity formally
derived by Starr~\cite{starr},
\begin{align}
  \label{eqn:miles:starr}
  F^2 = 1 + \frac 3{2d} \frac{ \int \eta^2 \, dx}{\int \eta\, dx}.
\end{align}
For waves of elevation with $\eta > 0$, the second term on the right
hand side of \eqref{eqn:miles:starr} is positive, immediately implying
that $F>1$. Recall that this inequality is sharp for
small-amplitude waves. Starr's
argument is formal in that he assumes that the integrals on the right
hand side are well-defined. For waves of elevation, however, this
assumption can be removed, as shown first by Amick and
Toland~\cite{at:finite} and then, much more simply, by
McLeod~\cite{mcleod:froude}.

\begin{theorem}[\cite{starr,at:finite,mcleod:froude}]\label{thm:miles:starr}
  All solitary waves of elevation are supercritical. Moreover, there
  are no solitary waves with critical Froude number $F=1$.
\end{theorem}

Combining Theorems~\ref{thm:miles:elevation} and
\ref{thm:miles:starr}, we see that supercriticality is equivalent to
being a wave of elevation. Craig and Sternberg~\cite{cs:sym} discovered that
such waves can be analyzed using a moving planes argument, leading to
the following celebrated result.

\begin{theorem}[\cite{cs:sym}]\label{thm:miles:sym}
  All supercritical waves of elevation are symmetric and monotone.
\end{theorem}

Theorems~\ref{thm:miles:elevation}--\ref{thm:miles:sym} leave open the
possibility of subcritical solitary waves with $F<1$, which have long
been conjectured not to exist. The analogous fact for the KdV model
\eqref{eqn:miles:kdv} is obvious; when $F<1$ the origin in the phase
portrait Figure~\ref{fig:miles:phase}(a) is a center, and there is
simply no room in phase space for a homoclinic orbit connecting the
origin to itself. For higher-order ODE models of water waves, however,
the situation is already quite subtle. There is still a
two-dimensional center manifold of periodic waves, but now also stable
and unstable manifolds which could in principle intersect to form a
large homoclinic orbit. A dimension-counting argument suggests that
this is unlikely, but that by no means proves that it is impossible.
For more on related issues in the context of ODEs, see
\cite{at:homoclinic, lombardi:book, champneys:persistence,
mho:saddlecenter}.

Very recently, this conjecture has been positively resolved by
by Kozlov, Lokharu, and the author.
\begin{theorem}[\cite{klw:nosubcrit}]\label{thm:miles:nosubcrit}
  There are no subcritical solitary waves.
\end{theorem}
The main ingredient in the proof of Theorem~\ref{thm:miles:nosubcrit}
is the introduction of an auxiliary function $\Phi$, called the
``flow force flux function'', which for irrotational solitary waves is
given by
\begin{align}
  \label{eqn:miles:ffff}
  \Phi(x,y) = \int_{-d}^y \Big( (\psi_y+c)^2 - \psi_x^2 \Big)\, dy'.
\end{align}
The definition of $\Phi$ is loosely motivated by the so-called ``flow
force'' \cite{benjamin:impulse}, an invariant for steady water waves
which is related to the balance of horizontal momentum. The surprising
and useful fact about $\Phi$ is that it solves the elliptic problem
\begin{equation}\label{eqn:miles:ffffprob}
  \left\{
    \begin{alignedat}{2}
      \Delta \Phi &= 0 &\qquad& \text{ in } {-d} < y < \eta(x),\\
      \Phi &= 0 &\qquad& \text { on }  y = -d,\\
      \Phi &\to 0 &\qquad& \text { as }  x \to \pm\infty,\\
      \Phi &= g\eta^2 \ge 0 &\qquad& \text { on }  y = \eta(x).
    \end{alignedat}
  \right.
\end{equation}
Applying the strong maximum principle to \eqref{eqn:miles:ffffprob},
for instance, we at once discover that $\Phi \ge 0$, which was not at
all obvious from the integral representation
\eqref{eqn:miles:ffff}.
More sophisticated versions of this function have appeared in
subsequent work of Lokharu \cite{lokharu:bounds}, and applied to
periodic as well as solitary waves.

\subsection{Non-existence of bores}\label{sec:miles:nobores}
We now show that \eqref{eqn:miles:stream} has no solutions which are
bores, that is, no solutions with distinct limits as $x \to
\pm\infty$. The argument, which is based on physical conservation
laws, goes back at least to Rayleigh; see \cite[Chapter
VIII, \S187]{lamb}.

\begin{figure}
  \centering
  \includegraphics[page=6]{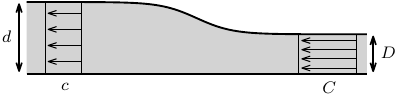}
  \caption{A bore limiting to distinct uniform flows as $x \to \pm\infty$.}
  \label{fig:miles:bore}
\end{figure}
Consider a hypothetical bore solution of \eqref{eqn:miles:stream}. One
can show that its limits as $x\to \pm\infty$ must be uniform flows,
similar to \eqref{eqn:miles:trivial}, but with differing speeds and
depths. Let $c,d$ denote the speed and depth of the flow at $x
=-\infty$, and $C,D$ the speed and depth at $x=+\infty$; see
Figure~\ref{fig:miles:bore}.
Using the third and fourth equations in \eqref{eqn:miles:stream}, we
see that these depths and speeds must be related by
\begin{align}
  \label{eqn:miles:conj:mass}
  cd &= CD, \\
  \label{eqn:miles:conj:energy}
  \tfrac 12 c^2 + gd &= \tfrac 12 C^2 + gD,
\end{align}
corresponding to mass and energy conservation, respectively. Another
calculation involving horizontal momentum and
the ``flow force'' from the previous subsection yields a third
condition,
\begin{align}
  \label{eqn:miles:conj:momentum}
  dc^2 + \tfrac 12 gd^2 &= DC^2 + \tfrac 12 gD^2.
\end{align}
Together,
\eqref{eqn:miles:conj:mass}--\eqref{eqn:miles:conj:momentum}
are called the ``conjugate flow equations'' for
\eqref{eqn:miles:stream}; see \cite{benjamin:conjugate}.

It is an easy algebraic exercise to show that
\eqref{eqn:miles:conj:mass}--\eqref{eqn:miles:conj:momentum} have no
solutions with distinct nonzero depths $d,D$. Let us briefly
outline a version of this argument which generalizes
to waves with vorticity~\cite[Lemma~3.8]{wheeler:pressure}.
Using \eqref{eqn:miles:conj:mass} to eliminate $C$, we view the
right hand sides of \eqref{eqn:miles:conj:energy} and
\eqref{eqn:miles:conj:momentum} as functions $\mathcal Q(D)$ and
$\mathcal S(D)$ of the downstream depth $D$ alone, and define $\tilde
{\mathcal S}=\mathcal S+D(\mathcal Q(d)-\mathcal Q)$. The desired
result now follows from the strict convexity of $\mathcal Q$ and
the fact that $\tilde {\mathcal S}'=\mathcal Q(d)-\mathcal Q$.

As we will see in Section~\ref{sec:sam:strat:bores}, the above
nonexistence proof breaks down in the presence of density
stratification. One can still derive conjugate flow equations, but they
are much more complicated, and often have nontrivial solutions, some
which have indeed been realized as the limiting states of bores.

\subsection{Existence of large-amplitude solitary waves}\label{sec:miles:large}

Having constructed solitary waves of small amplitude, and studied the
properties of an arbitrary solitary wave, we now wish to construct
solitary waves whose amplitudes are not small. For periodic waves, we
saw in Sections~\ref{sec:john} and \ref{sec:susanna} that
topological~\cite{rabinowitz:global} and
real-analytic~\cite{dancer:global,bt:analytic} global bifurcation
theory have been extremely successful. These theories cannot be
directly applied to the solitary wave problem, however, because of
many of the same issues we encountered in
Section~\ref{sec:miles:small}. Nevertheless, Amick and Toland were
able to circumvent these substantial difficulties in \cite{at:finite}
by approximating \eqref{eqn:miles:stream} with a family of problems
with better compactness properties. In \cite{at:periodic}, the same
authors were able to arrange for the approximate problem to be nothing
other than the periodic wave problem, thus constructing large solitary
waves as broadening limits of periodic waves
(Figure~\ref{fig:miles:broadening}).

\begin{figure}
  \centering
  \includegraphics[page=7]{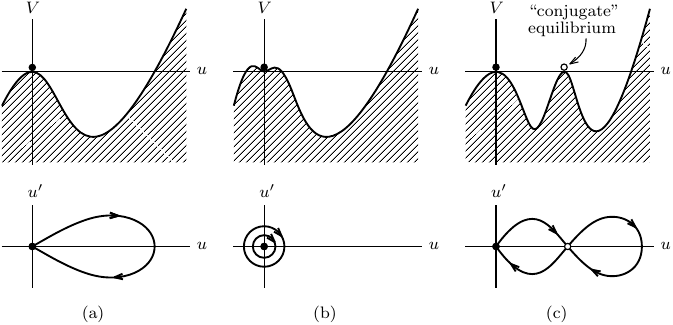}
  \caption{Potential functions $V$ and corresponding phase portraits.
  (a) Potential with a homoclinic orbit. (b) Perturbed potential for
  which the origin is now a center. (c) Perturbed potential for which
  the homoclinic orbit degenerates into two heteroclinic orbits.}
  \label{fig:miles:ball}
\end{figure}
It this subsection we will outline a different approach to solitary
wave problems, developed by Chen, Walsh, and the author in the series
of papers~\cite{wheeler:solitary, wheeler:pressure, cww:strat,
chen2020global}. To illustrate the method,  consider a family of
second order ODEs of the form
\begin{align}
  \label{eqn:miles:potential}
  u''+V_u(u;\lambda)=0,
\end{align}
modeling a small ball rolling along the graph of the potential function
$V$ under the influence of gravity. Suppose that, for some value of
$\lambda$, \eqref{eqn:miles:potential} has a homoclinic orbit as shown
in Figure~\ref{fig:miles:ball}(a). As we vary $\lambda$, what can
happen to this homoclinic? One possibility, shown in
Figure~\ref{fig:miles:ball}(b), is that the origin will flip from
a saddle to a center, in which case all small orbits will be periodic
and not homoclinic. We call this scenario ``spectral degeneracy''. Another
possibility is that $V$ develops a second critical point at the same
height as the origin, as shown in Figure~\ref{fig:miles:ball}(c). In
this case, the homoclinic loop from Figure~\ref{fig:miles:ball}(a)
degenerates into a pair of heteroclinic orbits. We call this scenario
``heteroclinic degeneracy''.

Very roughly speaking, the results in \cite{cww:strat,chen2020global}
state that, for a wide class of solitary-wave-type problems, such
``spectral'' and ``heteroclinic'' degeneracies are the \emph{only}
obstructions to applying global bifurcation theory. The results in
particular apply to the present problem \eqref{eqn:miles:stream}. In
this case, ``heteroclinic degeneracy'' is nothing other than the
broadening of a family of solitary waves into a bore
(Figure~\ref{fig:miles:broadening}), which was ruled out in
Section~\ref{sec:miles:nobores}. ``Spectral degeneracy'', on the other
hand, corresponds to the Froude number $F$ approaching the critical
value $F=1$. This can be analyzed by combining
Theorem~\ref{thm:miles:nosubcrit} on the nonexistence of solitary
waves with $F \le 1$ and the uniqueness results for small waves with
$F \approx 1$ alluded to in Section~\ref{sec:miles:small}.

\subsection{Some open problems}\label{sec:miles:open}
Many open questions remain, even in the classical setting we have
considered in this section. For instance:
\begin{enumerate}
\item[(1)] Does \cite{at:finite} construct \emph{all} solitary waves,
  or are there other connected components of solutions?
\item[(2)] Is there a set of parameters (mass, energy, amplitude,
  speed, etc.) which uniquely determine a solitary wave up to
  symmetries? Plotnikov~\cite{plotnikov:turning} has shown that the
  Froude number alone is not sufficient, even though it is sufficient
  for the small-amplitude waves in Section~\ref{sec:miles:small}.
\item[(3)] The solitary waves in \cite{at:periodic} are constructed as
  broadening limits of periodic waves. Can all solitary waves be
  realized in this way?
\item[(4)] Are small-amplitude solitary waves nonlinearly stable? Pego
  and Sun have shown that they are linearly
  stable~\cite{ps:stability}; also see \cite{lin:instability}.
\end{enumerate}
As seen in the other sections, the inclusion of additional
effects such as vorticity, stratification, surface tension, and
three-dimensionality significantly complicates the steady water wave
problem. This is certainly true for solitary waves and fronts, and
there are a great many open questions to be studied.

\section{Localized vorticity (by S. Walsh)}\label{sec:sam:localized}
\subsection{Introduction}
Water passing quickly over a blunt object can produce remarkably intricate spiral patterns, with distinct vortices appearing to shed from the object then propagating along in its wake.  The same phenomenon occurs at larger, even atmospheric scales in the form of von K\'arm\'an vortices swirling through clouds in the lee of mountains and islands.   Quite different from what we have encountered thus far, this type of vorticity is strikingly \emph{localized}: we observe isolated vortical regions moving through an expanse of irrotational flow.

To understand how this can happen, we return to the time-dependent (constant density) Euler equations \eqref{Euler} in $\mathbb{R}^2$.  Taking its curl leads to the $2$D \emph{vorticity equation},
\begin{equation}\label{sam:2-d vorticity equation}
	\frac{\partial \omega}{\partial t} + (\mathbf{u} \cdot \nabla) \omega = 0,
\end{equation}
which says, rather surprisingly, that $\omega$ is simply transported by the velocity field.  We have already seen one consequence of this fact, namely that $\omega$ is constant along streamlines for steady $2$D water waves.  But it also helps explain the vorticity distribution described above.  Intuitively, \eqref{sam:2-d vorticity equation} means that the vorticity of each fluid element is conserved.  So, if $\omega$ is initially localized --- say it is supported only on a collection of disjoint compact sets ---  it will remain localized for at least some period of time, perhaps a very long one.  Suppose that $\omega|_{t=0} = \chi_D$, where $\chi_D$ denotes the characteristic function of some bounded set $D \subset \mathbb{R}^2$.  From \eqref{sam:2-d vorticity equation} it follows that $\omega(t) = \chi_{D(t)}$, with $D(t)$ being the image of $D$ under the flow induced by $\mathbf{u}$.  This set may very well become extremely complicated, but the growth of its diameter is controlled by $\mathbf{u}$.  Incompressibility implies that the measure of $D(t)$ is conserved, so we know that vorticity can never completely suffuse the domain.

It is natural then to ask: can there be a traveling water wave for which $\omega$ is spatially localized?  While at first this might appear to be a simple variant of the literature discussed above, a little thought reveals that there are some serious obstructions to adapting earlier techniques.  Because $\omega$ is constant on the streamlines, in order for the vorticity to vanish at infinity, it must vanish identically along any unbounded streamline.  But, we have seen that the main strategies for constructing rotational water waves rely on perturbative methods (local bifurcation theory or spatial dynamics, for example) beginning at a shear flow.  Since \emph{every} streamline in a shear flow is unbounded, clearly some innovation will be needed.  Most critically, we must find exact (or approximate) solutions with many closed streamlines to use as the new starting point of our analysis.

Recently, there has been a lot of research in this direction, and some substantial progress has been made both in terms of constructing water waves with localized vorticity and in understanding their dynamical and qualitative properties.  We will present a number of these results, working from the most localized case (point vortices) to the least (spike vortices).

\subsection{Point vortices}
\begin{figure}
	\includegraphics[scale=0.8]{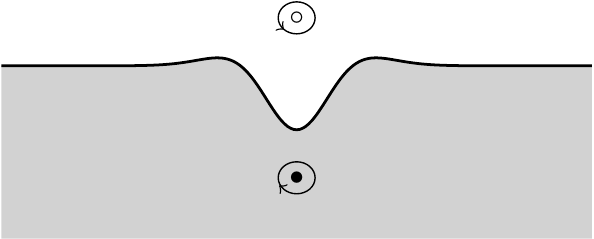}
	\includegraphics[scale=0.8]{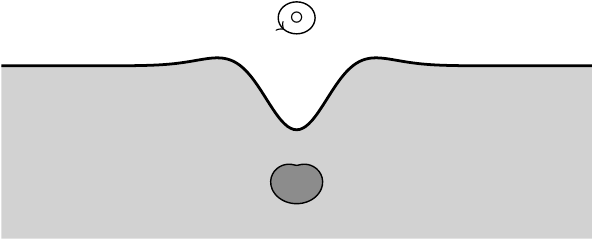}
	\caption{Two types of capillary-gravity waves with compactly supported vorticity from \cite{shatah2013travelling}.  On the left is a water wave with a submerged point vortex; on the right is a wave-borne vortex patch (the darker gray shaded region).  The white dot represents a counter rotating phantom vortex in the air.}
	\label{sam:point vortex figure}
\end{figure}

The simplest way to model localized vorticity is to shrink each vortical region down to a point. Doing so effectively ignores the small-scale structure of the flow there, freeing us to focus on the motion of the vortex as a whole.    Formally, a \emph{point vortex} corresponds to the case $\omega = \varepsilon \delta_{\mathbf{z}}$; the real number $\varepsilon \neq 0$ is its \emph{strength} and $\mathbf{z} = \mathbf{z}(t) \in \mathbb{R}^2$ its \emph{center}.  A collection of point vortices likewise means that $\omega$ is the sum of such Dirac measures.  Here we are following a long tradition in physics, where point masses and point charges are ubiquitous approximations.

Mathematically, though, this type of vorticity is too singular to give even a weak solution to the Euler equations.  Completely faithful adherence to \eqref{sam:2-d vorticity equation} would require that the center move according to $d\mathbf{z}/dt = \mathbf{u}(\mathbf{z})$.  But if $\omega$ is a Dirac measure, the magnitude of the velocity diverges as one approaches $\mathbf{z}$, which makes this ODE ill-defined.

Helmholtz \cite{helmholtz1858integrale} and Kirchhoff \cite{kirchhoff1876vorlesungen} independently arrived at a resolution to this issue.  The key physical reasoning is that the point vortex is not transported by the full velocity field because it does not self-advect; it is driven only by the \emph{irrotational} part of $\mathbf{u}$ at $\mathbf{z}$.   To make sense of this, observe that we can always split the velocity field as
\begin{equation}\label{sam:velocity splitting}
	\mathbf{u} = \nabla \Phi + \nabla^\perp \Psi
\end{equation}
where $\Phi$ is harmonic.  This decomposition is unique once appropriate boundary conditions have been prescribed.    As \eqref{sam:velocity splitting} implies $-\Delta\Psi = \omega = \varepsilon \delta_{\mathbf{z}}$, we must have that $\Psi = \varepsilon \Gamma(\placeholder - \mathbf{z}) + \Psi_{\mathcal{H}}$, with $\Gamma := -\frac{1}{2\pi} \log{|\placeholder|}$ being the fundamental solution of $-\Delta$ in $\mathbb{R}^2$, and $\Psi_{\mathcal{H}}$ a harmonic function.  The Helmholtz--Kirchhoff model states that the location of the vortex center satisfies
\begin{equation}\label{sam:helmholtz kirchhoff}
	\frac{d \mathbf{z}}{d t} = \left( \nabla \Phi + \nabla^\perp \Psi_{\mathcal{H}} \right)(\mathbf{z}),
\end{equation}
and $\mathbf{u}$ is a (distributional) solution of the Euler equations away from $\mathbf{z}$.  Having removed the singular self-interaction term $\varepsilon\nabla^\perp \Gamma$, the vector field above is analytic.  There are multiple ways to rigorously justify \eqref{sam:helmholtz kirchhoff} by taking limits of initial data with increasingly concentrated but smooth vorticity; see, for example, \cite{marchioro1993vortices,marchioro1994book,gallay2011interaction,glass2018point}.

Point vortices in fixed domains are a classical subject in fluid mechanics.  Our interest here is point vortices inside water waves, however, and this topic is far less well-explored.  Early papers by Ter-Krikorov \cite{terkrkorov1958vortex} and Filippov \cite{filippov1960vortex,filippov1961motion} studied traveling waves with point vortex forcing (that is, they do not require that \eqref{sam:helmholtz kirchhoff} holds).  Beginning with \cite{shatah2013travelling}, a number of recent papers have addressed wave-borne point vortices in a variety of regimes.  Let us sketch the basic idea behind these results.

Imagine we have a single point vortex carried by a solitary capillary-gravity wave in infinite-depth water.  Supposing that the undisturbed fluid domain is the lower half-plane, we impose the ansatz
\begin{equation}\label{sam:Psi ansatz}
	\Psi = \varepsilon \Gamma{(\placeholder - \mathbf{z})} -  \varepsilon \Gamma{(\placeholder - \mathbf{z^\prime})}
\end{equation}
where $\mathbf{z^\prime}$ is the reflection of $\mathbf{z}$ over the $x$-axis.  The second term is just a specific choice of $\Psi_{\mathcal{H}}$; it is indeed harmonic inside the fluid domain so long as $\mathbf{z^\prime}$ remains in the exterior.  Notice that this means  $\nabla^\perp \Psi \in L^2$, so its contribution to the kinetic energy is finite.  More whimsically, we can think of the second term in \eqref{sam:Psi ansatz} as a \emph{phantom vortex} in the air that has the opposite strength $-\varepsilon$.

Because $\Psi$ is harmonic in a simply connected neighborhood of the free boundary, it has a single-valued harmonic conjugate $\Theta$ there.  Near the surface we can therefore write $\mathbf{u} = \nabla (\Phi + \Theta)$.  Using familiar ideas from irrotational theory, one can then formulate the kinematic and Bernoulli conditions \eqref{WWBC} as nonlocal equations for the unknowns $\eta$, the trace $\varphi = \varphi(t,x) := \Phi(t, x, \eta(t,x))$, and $\mathbf{z}$.  The trace of $\Theta$ is small if $|\varepsilon| \ll 1$, and it is explicit given $\mathbf{z}$, so it behaves like a perturbative forcing term in the boundary conditions.  On the other hand, $\mathbf{z}$ evolves according to the ODE \eqref{sam:helmholtz kirchhoff}, with $\Phi$ obtained by harmonically extending $\varphi$.  The resulting evolution equation for $(\eta, \varphi, \mathbf{z})$ is thus quite reasonable; it is even Hamiltonian \cite{rouhi1993hamiltonian,varholm2020stability}.  Long-time well-posedness of the Cauchy problem absent surface tension was established by Su \cite{su2020longtime}.

As in Section~\ref{sec:walter}, the traveling wave equation is found by shifting to a moving reference frame and asking that the resulting relative velocity field be time independent. The Helmholtz--Kirchhoff condition \eqref{sam:helmholtz kirchhoff} becomes
\[
	c = \nabla \Phi(\mathbf{z}) - \varepsilon \nabla^\perp \Gamma{(\mathbf{z} - \mathbf{z^\prime})},
\]
where now $\mathbf{z}$ and $\mathbf{z^\prime}$ are fixed.  For simplicity, we may as well take them to be $-\mathbf{e}_2$ and $\mathbf{e}_2$, respectively.  Then the following theorem from \cite{shatah2013travelling} establishes the existence of small-amplitude solitary gravity waves with a submerged point vortex.

\begin{theorem}[Solitary waves with a point vortex] \label{sam:point vortex theorem}
For any $s > 3/2$, there exists $\varepsilon_0 > 0$, and a smooth one-parameter family
\[
	\left\{( \eta(\varepsilon), \varphi(\varepsilon), c(\varepsilon)) : \varepsilon \in (-\varepsilon_0, \varepsilon_0) \right\} \subset H^s(\mathbb{R}) \times \left( \dot H^s(\mathbb{R}) \cap \dot H^{\frac{1}{2}}(\mathbb{R}) \right) \times \mathbb{R}
\]
of solitary capillary-gravity water waves in infinite-depth with a submerged point vortex with center $\mathbf{z} = -\mathbf{e}_2$, strength $\varepsilon$, and wave speed $c(\varepsilon)$.  It bifurcates from the trivial solution $(0,0,0)$ at $\varepsilon = 0$.
\end{theorem}

Using the same formulation, one can obtain many other types of waves with a point vortex.  To name a few: in \cite{shatah2013travelling}, the same authors also study the periodic case $\omega =  \varepsilon\sum_{k \in \mathbb{Z}}  \delta_{\mathbf{z} + 2\pi k \mathbf{e}_1}$ and prove a global bifurcation result;  Varholm \cite{varholm2016solitary} constructed solitary capillary-gravity waves in finite-depth water with one or more point vortices;  and Le \cite{le2019existence} treated solitary waves carrying a submerged finite dipole in infinite depth.  All of these works rely on the fact that the linearized problem at $(0,0,0)$ is invertible --- which is due to the presence of surface tension --- and an implicit function theorem argument.   Refined qualitative properties, including sharp decay rates for $(\eta(\varepsilon),\varphi(\varepsilon))$, were later obtained in \cite{chen2019existence}, and a higher-order expansion of $\eta(\varepsilon)$ is derived in \cite{varholm2020stability}.  See Figure~\ref{sam:point vortex figure} for an illustration.

It is worth mentioning that the waves given by Theorem~\ref{sam:point vortex theorem} are (conditionally) orbitally stable \cite{varholm2020stability}, while Le \cite{le2019existence} proved that those carrying a submerged dipole are in fact orbitally unstable.  Both papers treat the small-amplitude regime where it turns out the water wave dynamics are dominated by those of the point vortices.  These are among the very few examples where nonlinear stability/instability of rotational water waves is known.

Finally, we note that C\'ordoba and Di Iorio \cite{cordoba2021existence} have very recently constructed periodic capillary-gravity waves with a point vortex or vortex patch using an entirely different approach.  Rather than start with the trivial solution, they gravity and vorticity perturb a Crapper wave.

\subsection{Vortex patches}
A \emph{vortex patch} describes the situation where $\omega$ is a function with compact, connected support $D$ (the \emph{patch}).  Most commonly, the vorticity is assumed to be constant on $D$, as in the example discussed at the beginning of this section.  While that means $\omega$ will have a jump discontinuity along $\partial D$, it nevertheless corresponds to a weak solution of the full Euler equations.  There is an enormous literature regarding the existence and stability properties of vortex patches in the plane, but again we wish to maintain our focus on water waves, and in that case much less is known.

Suppose that $D$ is positively separated from the free boundary.  We can still avail ourselves of the splitting \eqref{sam:velocity splitting} and say $\Psi = \Gamma * \omega + \Psi_{\mathcal{H}}$, for some harmonic function $\Psi_{\mathcal{H}}$.  In the typical case where $\omega = \chi_D$, the first term is fairly explicit and so it can be treated using analysis similar to the point vortex setting.  Let us discuss instead a family of solitary capillary-gravity waves with submerged vortex patches that have $C^0$ vorticity of a quite general type.

The strategy is to once more divide the task of constructing the wave into two parts. On the free surface, the kinematic and Bernoulli conditions must be satisfied.    Just like for point vortices, the vorticity's contribution there is perturbative, and for capillary-gravity waves the corresponding components of the linearized operator at the trivial solution are isomorphisms.  Things are more interesting inside the patch.   Conceptually, we imagine taking a point vortex and ``opening it up'' so that it becomes a set $D$ that is approximately a ball $B_r$ of radius $0 < r \ll 1$.  The main trouble is then to decide what the streamlines look like and select the value of $\omega$ along each one.  As in Section~\ref{sec:walter}, this can be done by specifying a vorticity function $\gamma$ and solving the semi-linear elliptic equation \eqref{vorticity eq} for the relative stream function $\psi$.  Note that this is a free boundary problem since the exact shape of the patch is a priori unknown.

We have a great deal of freedom in picking $\gamma$, but we do make a few restrictions:
\begin{equation}\label{sam:patch gamma assumption}
	 \gamma(0) = 0, \quad \gamma^\prime(0) < 0, \quad \gamma > 0 \quad \textrm{on } (-\infty,0).
\end{equation}
The first of these is to ensures continuity of the vorticity.  Without loss of generality, we can assume that $\psi$ vanishes on $\partial D$ (which is a streamline), so if $\omega$ is to be $C^0$, we need $\omega|_{\partial D} = \gamma(0) = 0$.  The other two are slightly more technical but essentially tell us that we have single-signed vorticity in the interior of $D$ and this persists under perturbations of $\gamma$.

Next we need to determine the relative stream function on the patch.  This motivates us to assume that when the semi-linear elliptic problem $\Delta \tilde\psi_0 = \gamma(\tilde\psi_0)$ is posed on the unit ball $B_1$, it has a negative, radial solution $\tilde\psi_0 \in H^2(B_1) \cap H_0^1(B_1)$, and the corresponding linearized operator
\begin{equation}\label{sam:patch gamma nondegen}
	\Delta - \gamma^\prime(\tilde\psi_0) : H^2(B_1) \cap H_0^1(B_1) \to L^2(B_1) \quad \textrm{is invertible.}
\end{equation}
The non-degeneracy condition \eqref{sam:patch gamma nondegen} is essentially generic and allows us to say that, for any domain $\tilde D \approx B_1$, there exists a slightly perturbed vorticity function $\tilde\gamma \approx \gamma$ and a unique solution $\tilde \psi \approx \tilde \psi_0$ to $\Delta \tilde \psi = \tilde \gamma(\tilde \psi)$ on $\tilde D$.  Finally, we obtain our stream function $\psi$ on $D$ as a rescaling of $\tilde \psi$.

In total then, we have accounted for the kinematic and Bernoulli conditions at the upper surface (which work as before) and the stream function problem in $D$ (solvable for any small deformation of a ball).  The last step is to impose a matching condition so that we can glue the solutions together and get a globally defined velocity field.  This can be accomplished through a somewhat delicate Lyapunov--Schmidt argument.  The result is the following small-amplitude existence theorem; see \cite[Theorem 2.3]{shatah2013travelling}.  An artistic rendering of the wave and patch is given in Figure~\ref{sam:point vortex figure}.

\begin{theorem}[Solitary wave with a vortex patch] Suppose that the vorticity function $\gamma$ satisfies \eqref{sam:patch gamma assumption} and \eqref{sam:patch gamma nondegen}.  For every $s > 3/2$, there exists $\varepsilon_0, \delta_0, \tau_0 > 0$ and a smooth family
\begin{align*}
	& \left\{ \left(\eta(\varepsilon,r, \tau),\, \varphi(\varepsilon,r, \tau),\,  \omega(\varepsilon,r, \tau),\, c(\varepsilon,r, \tau) \right) : (\varepsilon, r, \tau) \in [0,\varepsilon_0) \times (0, r_0) \times [0, \tau_0) \right\} \\
	& \qquad \qquad  \subset H^s(\mathbb{R}) \times \left( \dot H^s(\mathbb{R}) \cap \dot H^{\frac{1}{2}}(\mathbb{R}) \right) \times H^1(\mathbb{R}^2) \times \mathbb{R}.
\end{align*}
of traveling solitary capillary-gravity waves in infinite depth with a submerged vortex patch.   At each parameter value, the patch is a compact $H^s$ domain whose boundary is given by
\[
	\partial D(\varepsilon, r, \tau) = \left\{ r \left( \cos{\theta} + \tau \sin{(2\theta)}, \, \sin{\theta} - \tau \cos{(2\theta)} \right)  + O(r^2(r+\varepsilon)) : ~ \theta \in [0,2\pi) \right\}.
\]
Moreover, each $\omega(\varepsilon,r,\tau)$ is globally Lipschitz continuous and $H^{s+\frac{1}{2}}(D(\varepsilon,r,\tau))$.
\end{theorem}

\subsection{Spike vortices}

\begin{figure}
\includegraphics{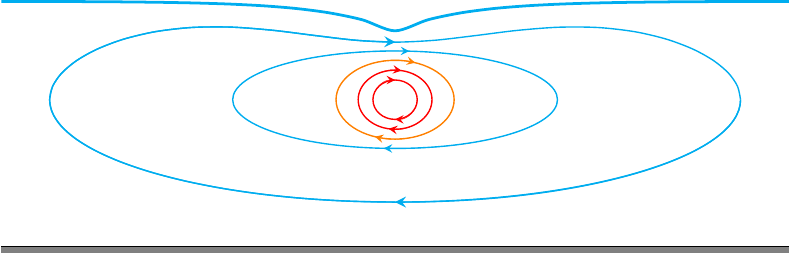}
\caption{
Streamline pattern for the stationary capillary-gravity wave with a spike vortex.  The arrows show the direction of the flow and the colors indicate the sign of the vorticity: red for negative vorticity; blue for positive; and orange for the unique interior streamline on which $\omega = 0$.
}
\label{sam:spike vortex figure}
\end{figure}

The vorticity localization in the previous two examples is rather severe. Indeed, $\omega$ is supported either at a point or on a small, compact set.  Can there be a wave for which the vorticity vanishes at infinity but is not completely confined to a bounded subdomain of the fluid?  This question was considered in the recent paper \cite{ehrnstrom2019smooth}, where a family of stationary waves in finite-depth water were constructed with exponentially decaying $\omega$.  While vortex patches are in some sense desingularized point vortices, these so-called \emph{spike vortex} waves are a different species entirely.

To obtain them we change tack yet again.  The separation between the support of $\omega$ and the free boundary was central to the analysis of the previous subsections as it allowed us to decouple the vortical considerations in the bulk from the determination of the free surface.  That is no longer possible, since the vorticity comes right up to the top.  Moreover, the streamline pattern for a wave carrying a spike vortex must be even more dramatically non-laminar: we cannot have \emph{any} unbounded streamlines besides the free surface and bed.  In particular, this requires there be a critical layer --- a line of stagnation points --- stretching from upstream to downstream infinity.  It is also more natural to consider stationary waves, that is waves with speed $c = 0$, since the far-field flow cannot be laminar.

Our new strategy is inspired by the study of spike solutions to singularly perturbed elliptic PDE. A prototypical example of such a problem is
\begin{equation}\label{sam:singularly perturbed PDE}
	r^2 \Delta u = \gamma (u) \quad \textrm{in } K, \qquad u = 0 \quad \textrm{on } \partial K,
\end{equation}
where $K \subset \mathbb{R}^n$ is fixed and bounded, and $0 < r \ll 1$ is a spatial scale parameter.  These arise as models for chemotaxis and certain concentration phenomena in biology and chemistry; see \cite{lin1988large,ni1991shape,ni1995location,
cao1996existence,ni1998location,ni1991shape,li1998dirichlet} for example.  The role of the parameter $r$ is the following.  Suppose that there exists a solution to the unscaled PDE when it is posed on the whole space
\[
	\Delta U = \gamma(U) \qquad \textrm{in } \mathbb{R}^n,
\]
such that $U$ is positive, radial, monotonically decreasing, and exponentially localized.  A sufficient (but not necessary) condition ensuring this is
\begin{equation}\label{sam:vortex spike gamma}
	\gamma(t) = t - |t|^p t,
\end{equation}
for some integer $p \geq 2$; see \cite{li1998dirichlet,ehrnstrom2019smooth}.

Then the rescaled solution $U(\placeholder/r)$ satisfies the PDE \eqref{sam:singularly perturbed PDE} in $K$, and it \emph{almost} satisfies the boundary conditions when $0 < r \ll 1$ due to the exponential decay.  One therefore searches for exact solutions with the ansatz
\[ u = U\left(\frac{\placeholder - \mathbf{z}}{r} \right) + v,\]
where $\mathbf{z} \in K$ is chosen to be as far from $\partial K$ as possible, and $v$ is a \emph{boundary correction}.  The point is that $v$ will be small, so we can hope to find it via perturbative methods.

A highly nontrivial adaptation of these ideas can be applied to the water wave problem.    We must contend with the unboundedness of the domain, a fully nonlinear boundary condition, and of course the free boundary itself.  Omitting all of the details, the ultimate result is given below.  Figure~\ref{sam:spike vortex figure} is a qualitative depiction of the shape of the wave and streamlines.

\begin{theorem}[Spike vortex]
Let $\gamma$ be given as in \eqref{sam:vortex spike gamma}. There exists $r_0 > 0$ such that, for all $r \in (0, r_0)$ there is a stationary capillary-gravity wave in finite-depth whose vorticity takes the approximate form
\[
	\omega \approx \frac{1}{r^2} \gamma\left( U\left(\frac{\placeholder -\tau \mathbf{e}_2}{r} \right)\right)
\]
where $\tau = \tau(r)$ is a vertical translation.
\end{theorem}

These are very small-amplitude waves with very small total vorticity.  However, the vorticity is spiked in that $\omega  \to 0$ in measure while $\| \omega \|_{L^\infty} \to \infty$ as $r \searrow 0$. Unlike the point vortices and vortex patches, the kinetic energy carried by these waves is non-perturbative:
\[
	\frac{1}{2} \| \nabla \psi \|_{L^2}^2 \approx \frac{1}{2} \| \nabla U(\placeholder - \tau e_2 )\|_{L^2(\tfrac{1}{r} \Omega)}^2 = O(1) \qquad \textrm{as } r \searrow 0.
\]
That is, we are not bifurcating from $0$.  Actually, this isn't a traditional local bifurcation at all: the existence of these waves is proved by first using a Lyapunov--Schmidt argument to reduce the system to a one-dimensional problem for $\tau$, then applying the intermediate value theorem!

\section{Stratified water waves (by S. Walsh)}\label{sec:sam:strat}
\subsection{Continuous stratification}
We have already discussed in Section~\ref{sec:susanna} a specific example of density stratified water waves and, in Section~\ref{sec:miles}, we hinted at the effects of stratification on conjugate flow analysis. Let us now delve a bit deeper into the literature with a decided bias towards new developments.

We start by considering the case of a single, finite-depth fluid domain $\Omega$ bounded above by a free surface $\{ y= \eta(x)\}$, and exhibiting non-constant density $\rho$ that is continuous.  Exact steady solutions to this system date back to Dubreil-Jacotin \cite{dubreil1937theoremes}, who used a power series approach to construct small-amplitude periodic gravity waves.  Through quite different methods,  Yanowitch \cite{yanowitch1962gravity}, Ter-Krikorov \cite{ter1963theorie}, and Turner \cite{turner1984variational} subsequently developed their own small-amplitude existence theories.  The latter two of these additionally treated solitary waves using the classical plan of first building periodic waves, then sending the period to infinity.  They both assume that the resulting solitary wave limits to constant velocity upstream and downstream, and, for that reason, only considered a specific class of Bernoulli functions.  Large-amplitude periodic waves with arbitrary (smooth) $\beta$ and $\varrho$ were obtained much later in \cite{walsh:stratified,walsh2014local,walsh2014global} via an adaptation of the Constantin--Strauss \cite{CS2004} technique to the heterogeneous regime.

Similarly, the existence of large-amplitude solitary waves with a general background current was only recently proved in \cite{cww:strat}.   Here we quote a simplified version of their main result.

\begin{theorem}[Stratified solitary waves] \label{sam:global solitary stratified theorem}
Fix $\alpha \in (0,1/2]$, and let a streamline density function $\varrho \in C^{2+\alpha}$, and asymptotic relative velocity $\mathbf{u}^* \in C^{2+\alpha}$ be given.  There exists a global curve of stratified solitary waves parameterized by the Froude number $F$ with the regularity $(\mathbf{u}, \eta) \in C^{2+\alpha} \times C^{3+\alpha}$.  It bifurcates from $(\mathbf{u},\eta,F) = (\mathbf{u}^*, 0, F_{\mathrm{cr}})$, and at its extreme the waves limit to (horizontal) stagnation.
\end{theorem}

This theorem is obtained using an abstract global bifurcation theory based on the characterization of ``loss of compactness'' discussed in Section~\ref{sec:miles:large}.  A key ingredient is the non-existence of monotone front-type solutions (which we call hydrodynamic \emph{bores}) for this problem.   Together with other qualitative results, this fact allows one to rule out heteroclinic degeneracy as a potential limiting behavior along the global curve.  The curve cannot be a closed loop (because the point of bifurcation is singular), so the only possibilities are that the curve is unbounded or that the problem degenerates along it.  Both of these alternatives are then shown to lead to horizontal stagnation.

A parallel set of results treats the case of waves in a stratified fluid bounded both above and below by rigid horizontal walls; we will call this system \emph{channel flow} and its steady solutions \emph{internal waves}.   A pleasant feature of continuously stratified channel flow is that the free boundary and the Bernoulli boundary condition are no longer an issue.   Small-amplitude waves were constructed by Ter-Krikorov \cite{ter1963theorie}, Turner \cite{turner1981internal}, Kirchg\"assner \cite{kirchgassner1982wavesolutions}, Kirchg\"assner and Lankers \cite{kirchgassner1993structure}, James \cite{james1997small}, and Sun \cite{sun2002solitary}. Large-amplitude existence theory for continuously stratified channel flows was provided by Bona, Bose, and Turner \cite{bona1983finite}, Amick \cite{amick1984semilinear}, and Lankers and Friesecke \cite{lankers1997fast}.

It is important to note that all of the work mentioned thus far assumes there are no stagnation points.  A major focus of current research is periodic stratified waves that have critical layers or very low regularity density; see, for example, \cite{escher2011stratified,henry2014global,wheeler2019stratified,escher2020stratified,haziot2021stratified}.

\subsection{Periodic and solitary waves in multilayer fluids}

\begin{figure}
	\includegraphics[scale=0.8]{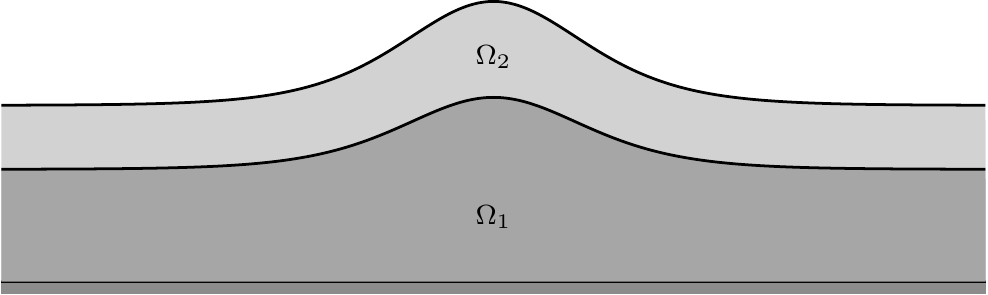}
	\caption{Two-layer solitary gravity waves studied in \cite{sinambela2020large}.    The upper fluid occupies the domain $\Omega_2$, while the the lower fluid lies in $\Omega_1$.  The interface dividing them is a free boundary along which the pressure and normal velocity are continuous but the tangential velocity and density may jump.}
	\label{sam:multilayer figure}
\end{figure}

We next consider systems comprised of multiple immiscible fluids arranged in vertical layers.  Modifying somewhat the model introduced in Section~\ref{sec:walter}, let us take the fluid domain to be given as the disjoint union
\[
	\Omega = \Omega_1 \cup \Omega_2 \cup \cdots \cup \Omega_M,
\]
for some $M \geq 2$.  Each $\Omega_k$ represents a single layer and the boundaries separating them are assumed to be free.  The restriction of the velocity and density to a layer will be smooth, but they can in general jump as one passes between phases.  In the absence of surface tension, the pressure is continuous throughout $\overline{\Omega}$; otherwise, it will jump across the interfaces according to the Young--Laplace law, resulting in a transmission-type boundary condition.  A two-fluid example is shown in Figure~\ref{sam:multilayer figure}.

Multilayer models are popular in oceanography as an idealization of observed stratification patterns.  In the field, one finds large regions of nearly constant density separated by much thinner regions where the density varies rapidly.  It is quite natural, then, to collapse the transition regions down to curves along which the density now jumps.  (In fact, this sort of sharp interface limit approach has been used to analytically construct multilayer waves \cite{turner1981internal,james2001internal}.)  Being even greedier, we might suppose that the density in each layer is constant and the flow is irrotational.  Approximating a continuously stratified fluid by a many-layered one in this way is a widespread practice; it has been rigorously justified for a certain class of periodic waves in \cite{chen2016continuous}.  Frequently, applied authors consider the case of exactly two-layers.  Three-layer models have been proposed \cite{grue2002solitary}, though,  and numerical evidence suggests that they support more complicated traveling waves such as breathers \cite{nakayama2020breathers}.

Multilayer gravity waves are the most extensively studied.  While these are completely mathematically valid as steady solutions of the Euler equations, the dynamical problem is actually ill-posed \cite{shatah2011interface} due to a Kelvin--Helmholtz instability that creates explosive growth at high frequencies.  Introducing any capillarity at all quashes the effect, however; see the discussion in \cite{lannes2013stability}.

Small-amplitude two-layered solitary and generalized solitary gravity waves with general stratification were first obtained by Wang \cite{wang2017small} using spatial dynamics techniques.  Recently, Sinambela \cite{sinambela2020large} proved a global bifurcation result that essentially extends Theorem~\ref{sam:global solitary stratified theorem} to the multilayer setting.  One of the main challenges in doing so is the absence of good bounds on the pressure.  This is a recurring theme in the study of multilayer waves: maximum principle arguments are deeply frustrated by a lack of quantities that are continuous across the interface.   In \cite{sinambela2020large}, Sinambela was able to control the pressure in part through a clever application of the almost monotonicity formula of Caffarelli--Jerison--Kenig \cite{caffarelli2002monotonicity} to Yih's equation \eqref{yih long_susanna}.  Gravity waves in a channel have been studied by many authors.  We note in particular that a large-amplitude existence theory was given by Amick and Turner \cite{amick1986global}.  As we will see shortly, monotone front-type solutions do exist in this regime, so it is an open question whether heteroclinic degeneracy occurs along the global curve \cite{amick1984limiting}.

There is a smaller but still quite well-developed literature on interfacial capillary-gravity waves.  These are interesting both in their own right and because of the well-posedness issues raised above.  Via global bifurcation, large-amplitude periodic capillary-gravity waves in a channel were obtained by Ambrose, Strauss, and Wright \cite{ambrose2016global}.  Solitary internal capillary-gravity waves in a channel have been analytically by Kirrmann \cite{kirrmann1991reduction} and Nilsson \cite{nilsson2017internal} using spatial dynamics techniques, and numerically by Laget and Dias \cite{Laget1997interfacial}.  Some of these waves have been shown to be orbitally stable \cite{chen2021orbital}.

Multilayer flows are also studied in connection to the wind-generation of water waves.  A fundamental problem in geophysics is to describe the primary mechanisms by which energy is transferred to the water to create persistent surface waves.  The air--sea system can be thought of as a two-layer stratified fluid where the lower layer is about $1000$ times denser than the upper one.  The quasi-laminar model of Miles \cite{miles1957windwaves1,miles1959windwaves2} attributes the creation of ocean waves by wind to an instability brought on by a critical layer in the air region; in particular, this says atmospheric vorticity plays a key role.  For mathematical work in this direction, see \cite{walsh2013wind,buhler2016wind,le2018transmission}.

\subsection{Bores} \label{sec:sam:strat:bores}
\begin{figure}
	\includegraphics[scale=0.8]{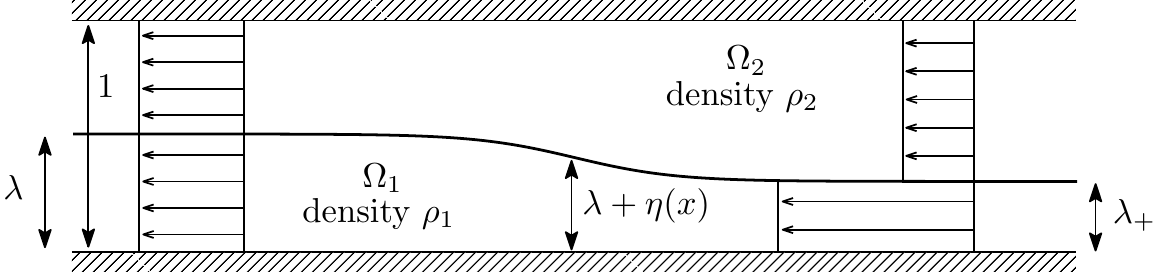}
	\caption{Bore system studied in \cite{chen2020global}.  The layers have constant but differing densities, and the flow in each is irrotational.  Upstream, the velocity tends to $(-1,0)$ and the lower layer has thickness $\lambda$.  }
\end{figure}

We conclude our tour of stratified waves with a discussion of bores, which have been the subject of recent work.  Specifically, we seek front-type solutions in a channel with two fluid layers.  For simplicity, suppose that the flow is irrotational flow in each region:
\begin{subequations} \label{sam:bore problem}
\begin{equation}
	\Delta \psi  = 0 \qquad \textrm{in } \Omega_1 \cup \Omega_2.
\end{equation}
The Bernoulli condition on the internal interface takes the form
\begin{equation}
	\frac{1}{2} \jump{\rho |\nabla \psi|^2} + \frac{\jump{\rho}}{F^2} y = \frac{\jump{\rho}}{2}   \qquad  \textrm{on } y = \eta(x)
\end{equation}
where $\jump{\placeholder} = (\placeholder)_2 - (\placeholder)_1$ denotes the jump of a quantity over the internal interface.  Suppose the velocity at infinity is uniform upstream, which after normalizing implies that
\begin{equation}\label{sam:bore upstream}
\left\{
\begin{aligned}
	\nabla \psi &\to (0,-1), & \eta &\to 0 & \qquad &\textrm{as } x \to -\infty \\
	\psi_x &\to 0, & \eta &\to \lambda_+ - \lambda \neq 0 & \qquad & \textrm{as } x \to +\infty.
\end{aligned}
\right.
\end{equation}
Notice that the velocity downstream is also purely horizontal and constant, but it need not be the same constant in each layer.  Here $\lambda_+$ denotes the limiting width of the lower layer downstream and $\lambda$ is the upstream width.  Based on conjugate flow analysis, it can be shown that for fixed (constant) densities $\rho_1$ and $\rho_2$, bores can only exist for a specific value of the Froude number and downstream layer depth  \cite{Laget1997interfacial}.  Thus the upstream layer depth $\lambda$ is the sole free parameter.  As the walls and internal interface are streamlines, we can infer from \eqref{sam:bore upstream} that
\begin{equation}
\left\{
\begin{aligned}
	\psi & = 0 & \qquad & \textrm{on } y = \eta(x) \\
	\psi & = \lambda &  \qquad & \textrm{on } y = -\lambda \\
	\psi & = \lambda - 1& \qquad & \textrm{on } y = 1-\lambda.
\end{aligned}
\right.
\end{equation}
\end{subequations}

Small-amplitude bores were first constructed by Amick and Turner \cite{amick1989small}; similar results using different methods were later obtained in \cite{mielke1995homoclinic,makarenko1992bore,makarenko1999conjugate,chen2019center}.  These were very recently extended to the large-amplitude regime \cite{chen2020large,chen2020global}:

\begin{figure}
	\includegraphics[scale=0.8]{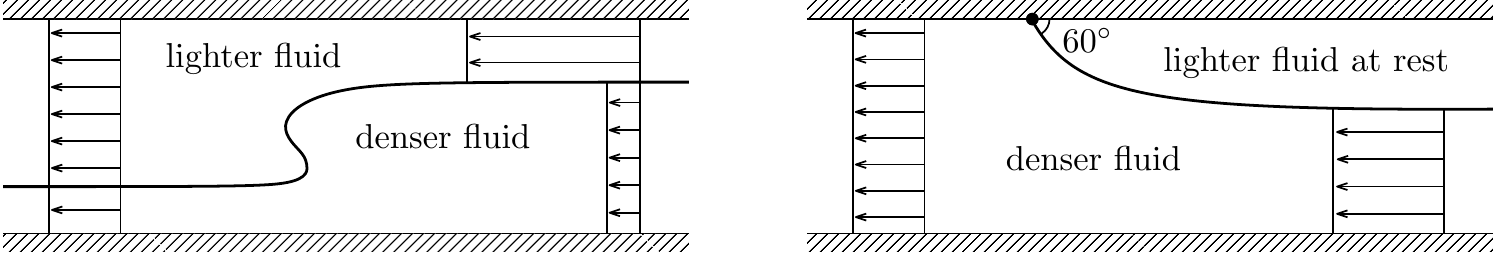}
\caption{Conjectured extremal configurations for the curves of bores from Theorem~\ref{sam:global bore theorem}.  On the left, the wave has overturned in that the interface has a vertical tangent.  On the right, the interface has touched the upper wall forming a \emph{gravity current}.  The contact angle of $60^\circ$ was computed formally by von K\'arm\'an \cite{vonkarman1940engineer,benjamin1968gravity}.  }
\label{sam:bore figure}
\end{figure}

 \begin{theorem}[Large bores] \label{sam:global bore theorem}
  For any fixed constant densities $0 < \rho_2 < \rho_1$, there exist two $C^0$ curves $\mathscr{C}^\pm$ of classical solutions to the internal wave problem \eqref{sam:bore problem}.  Each wave on $\mathscr{C}^\pm$ is a strictly monotone bore with $\pm  \eta_x < 0$.
  Following $\mathscr{C}^\pm$, we encounter waves that are arbitrarily close to having a horizontal stagnation point on the internal interface.
  \end{theorem}

This result is obtained through a new abstract global bifurcation theorem specialized to study heteroclinic solutions of elliptic PDE set on infinite cylinders.  We have already discussed the intuition underlying this machinery in Section~\ref{sec:miles:large}, namely that to ensure that the global curve is unbounded, we must be able to rule out spectral degeneracy and heteroclinic degeneracy.   For two irrotational layers, the conjugate flows are easily classified, and indeed both of these undesirable scenarios can be excluded.  It should be noted that adding constant vorticity, say, would considerably complicate the conjugate flow equations \cite{chen2019center}, as would having even a single additional layer \cite{lamb2000three}.

Theorem~\ref{sam:global bore theorem} guarantees that a stagnation point develops on the interface as one approaches the extreme of both $\mathscr{C}^\pm$.  For Stokes waves, this is known to coincide with the development of a corner, but the situation here is actually quite different.  Figure~\ref{sam:bore figure} shows the two scenarios observed in the numerical work \cite{dias2003internal}.   In terms of rigorous results, using free boundary regularity theory among other tools, the authors in \cite{chen2020global} are able to winnow the alternatives down to the following.  For the limit along $\mathscr{C}^-$, either the interface \emph{overturns} in that $\| \eta_x \|_{L^\infty}$ is unbounded or the flow in both layers limits to stagnation at a single point on the interface.  The second scenario has not been observed numerically and is conjectured never to occur.  On the other hand, following $\mathscr{C}^+$, it is proven that either the interface overturns or it comes into contact with the upper wall.

\section{Three-dimensional waves (by E. Wahl\'en)}\label{sec:erik}
\subsection{Introduction}
In this last section we consider steady water waves in three dimensions. They  have received considerably less attention than 2D waves, a strong reason being that the problem is mathematically much more challenging. Most of the existing theory is for irrotational flows  with surface tension. In fact, there is no clear analogue of the stream function formulation, which makes it very hard to even begin investigating solutions with vorticity. Already in the irrotational case, some 2D arguments based on the stream function formulation are hard to replace. It also turns out that surface tension plays a much greater role than in 2D; without surface tension, the combined problem for the velocity field and the surface profile is not elliptic, and one encounters problems with small divisors. Therefore, most existence results are for gravity-capillary waves. On the other hand, the 3D world offers many interesting geometrical possibilities in the sense that one can allow the waves to have different behaviours (for example, periodic, quasiperiodic or solitary) in different horizontal directions. In this section we will primarily focus on {\em doubly periodic waves}, that is, solutions which are periodic in two different horizontal directions, and  {\em fully localised solitary waves}, that is, solutions which decay in all horizontal directions. We will discuss both irrotational theory and recent results for waves with vorticity. The methods  will typically be related to local bifurcation theory or the calculus of variations. At this time there is a complete lack of large-amplitude theory. A very successful method for finding small-amplitude solutions with different behaviours in different directions is {\em spatial dynamics}. This involves  formulating the equations as an infinite-dimensional, typically ill-posed dynamical system, where one of the horizontal variables plays the role of time (see also Section \ref{sec:miles} for applications in 2D). In the interest of keeping the section to a reasonable length, we will only touch  briefly on the subject, even though it is important. For a good overview of spatial dynamics methods for 3D irrotational water waves until 2007 we refer to \cite{Groves07}. Some more recent additions are \cite{BagriGroves14,  DengSun09, DengSun10, GrovesSunWahlen16a, Haragus15, Nilsson19b}.

\subsubsection{Notation and problem formulation} We begin with some general comments and notation. The density $\rho$ is assumed constant throughout and we pick units so that it has value one\footnote{For a general value of $\rho$ one just has to replace $P$ and $\sigma$ by $P/\rho$ and $\sigma/\rho$ everywhere.}.
We also assume that the only external force is gravity, $\mathbf{F}=- g \mathbf{e}_3$ ($\mathbf{e}_j$ denoting the standard basis vectors).
In contrast to the rest of this paper, $x$ and $y$ are the horizontal variables and $z$ the vertical. It is convenient to introduce the notation
$\mathbf{x}'=(x,y)$, so that $\mathbf{x}=(\mathbf{x}', z)$ and the  fluid domain is given by
\[
\Omega= \{(\mathbf{x}', z) \in \mathbb{R}^2 \times \mathbb{R} : -d<z<\eta(\mathbf{x}')\}.
\]
We will also use the notation $\mathbf{u}=(u, v, w)$, but whenever it is convenient we will refer to components of vectors by their numbers.
By rotation invariance, one can always assume that the waves are travelling in the $x$-direction as in Section \ref{sec:walter}, but it is sometimes convenient to allow for a general direction.
We therefore assume the more general form $\eta=\eta(\mathbf{x}'+\mathbf{c}t)$, $\mathbf{u}=\mathbf{u}(\mathbf{x}'+\mathbf{c}t, z)$, $P=P(\mathbf{x}'+\mathbf{c}t, z)$, where we've also changed a sign for later convenience. The Euler equations for a travelling wave then take the form
\begin{alignat*}{2}
(\mathbf{u}+(\mathbf{c}, 0)) \cdot \nabla \mathbf{u}&=-\nabla P -g \mathbf{e}_3,\\
\nabla \cdot \mathbf{u}&=0.
\end{alignat*}
For notational simplicity
we will replace the relative velocity $\mathbf{u}+(\mathbf{c}, 0)$ by $\mathbf{u}$, which  eliminates the wave velocity $\mathbf{c}$ from the problem. Written out in full, the problem we want to solve is
\begin{subequations}
\label{eqn:erik:ww problem}
\begin{alignat}{2}
\mathbf{u} \cdot \nabla \mathbf{u}&=-\nabla P -g \mathbf{e}_3 && \quad \text{in $\Omega$,} \label{eqn:erik:Euler} \\
\nabla \cdot \mathbf{u}&=0 && \quad \text{in $\Omega$,} \\
\mathbf{u}\cdot \mathbf{n}&=0 && \quad \text{on $\partial \Omega$,} \\
P&=-\sigma  \nabla \cdot \Big(\tfrac{\nabla \eta}{\sqrt{1+|\nabla \eta|^2}}\Big) && \quad \text{on $S$,}  \label{eqn:erik:dynamic BC}
\end{alignat}
\end{subequations}
where the last condition is the dynamic boundary condition, which in the case of surface tension is given by the Laplace-Young equation saying that the jump in pressure from the air region into the water region is proportional to the mean curvature; here $\sigma\ge 0$ is  the (constant) coefficient of surface tension and we have set $P_{atm}=0$ without loss of generality. The vorticity vector will be denoted as
\[
\boldsymbol{\omega}\coloneqq \nabla \times \mathbf{u}.
\]

\subsubsection{Symmetries} Throughout the section, different symmetries will be important. Let $R_1$ and $R_2$ be the reflections of the first and second horizontal coordinates, respectively. The water wave problem has the symmetries
\[
T_j\colon (\eta, \mathbf{u}, P)\mapsto (\eta(R_j), R_j\mathbf{u}(R_j), P(R_j)), \quad j=1,2,
\]
and
\[
S\colon (\eta, \mathbf{u}, P)\mapsto (\eta, -\mathbf{u}, P).
\]
We can take advantage of these by considering solutions with the symmetries
\[
ST_1(\eta, \mathbf{u}, P)=(\eta, \mathbf{u}, P)=T_2(\eta, \mathbf{u}, P).
\]
This means in particular that $\eta$ and $u$ are even in $x$ and $y$,  $v$ is odd in both $x$ and $y$ and $w$ is odd in $x$ and even in $y$. We will simply call such solutions {\em symmetric}.
As we will see, in certain cases we only have the combined symmetry $ST_1T_2$, meaning that $w$ is odd in $\mathbf{x}'$, while all other variables are even. We call such solutions {\em weakly symmetric}.

\subsection{Irrotational theory}

\begin{figure}
  \centering
  \includegraphics[width=0.3\linewidth]{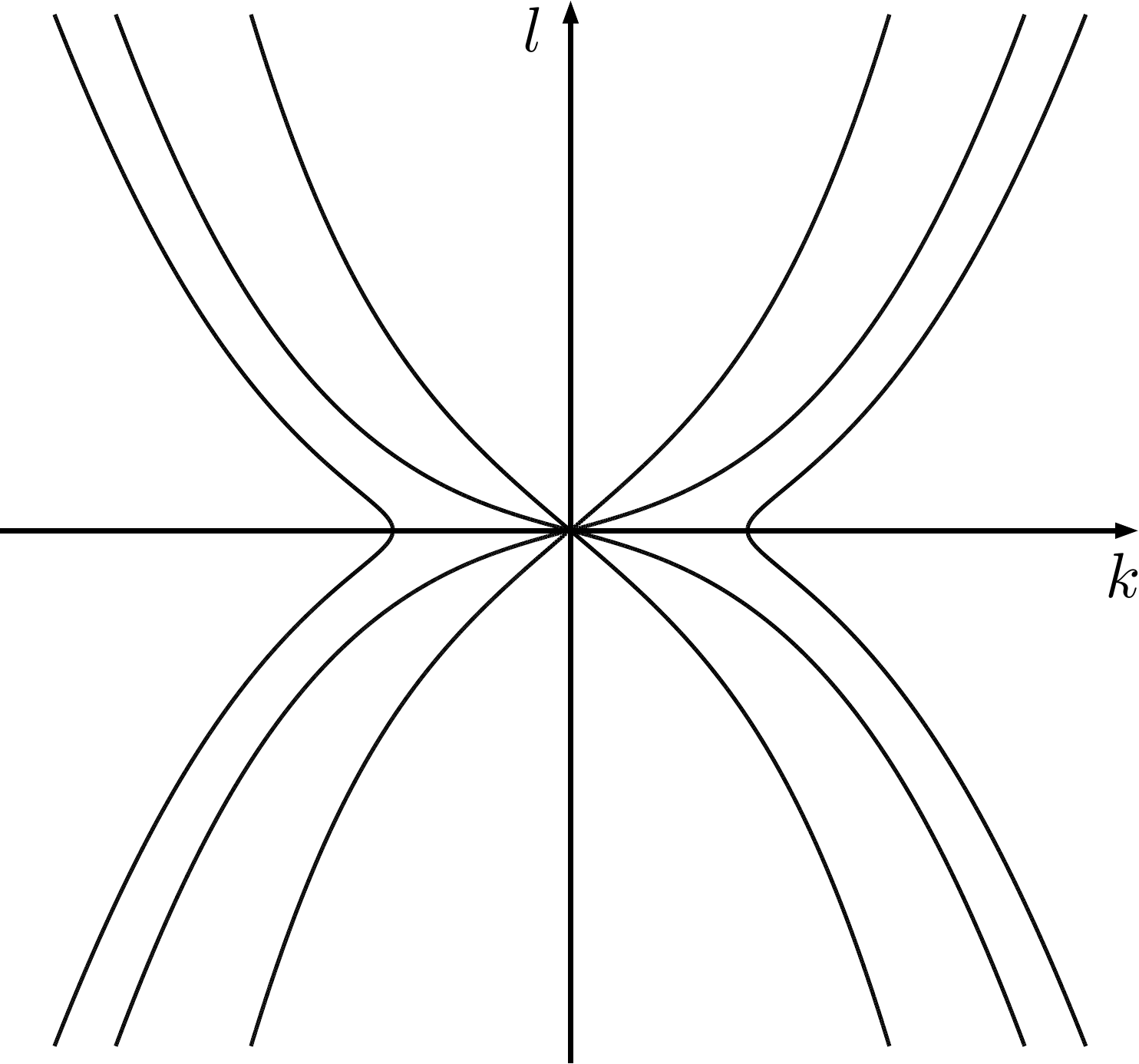} \hspace{5mm}
  \includegraphics[width=0.3\linewidth]{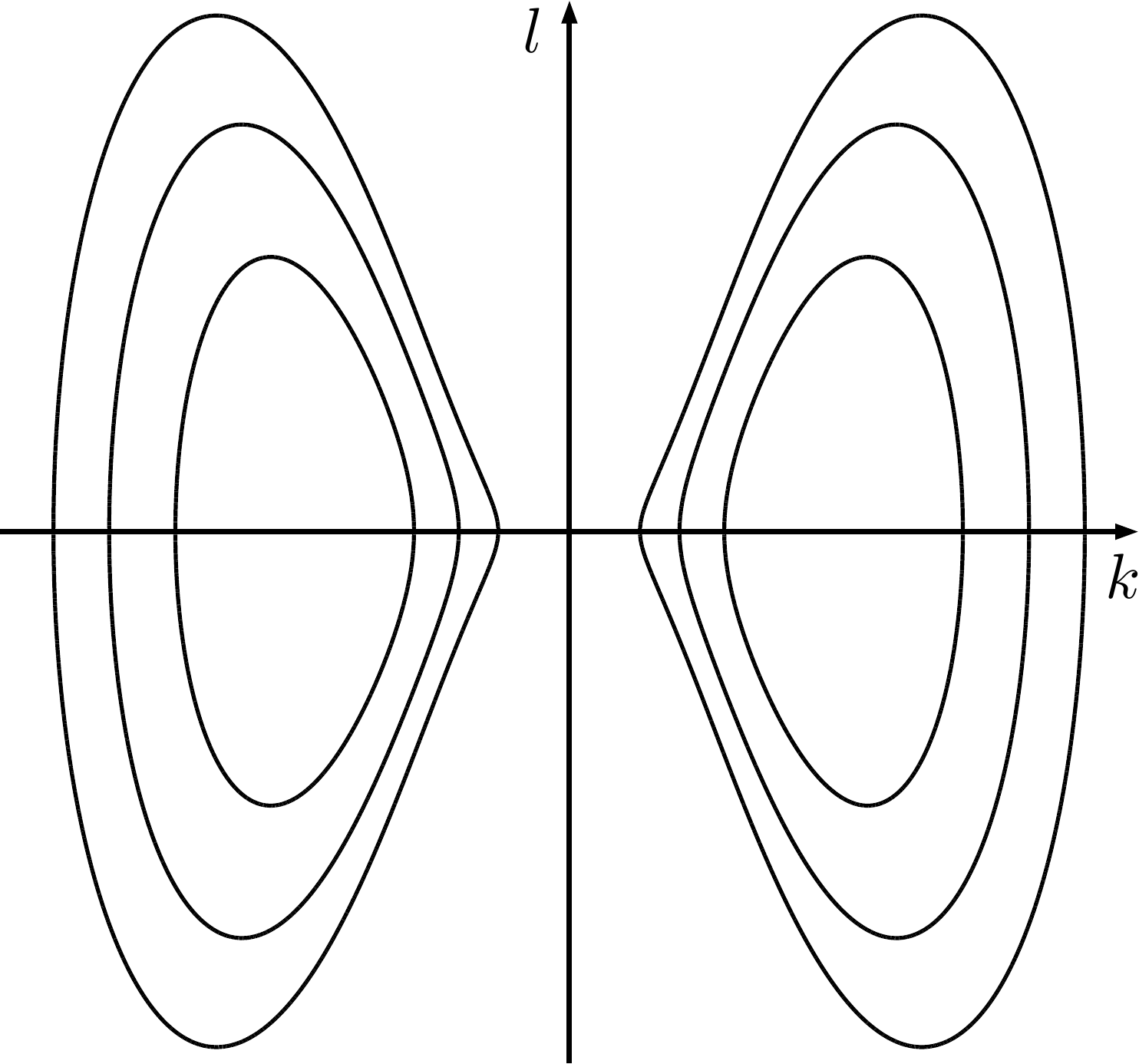}\hspace{5mm}
   \includegraphics[width=0.3\linewidth]{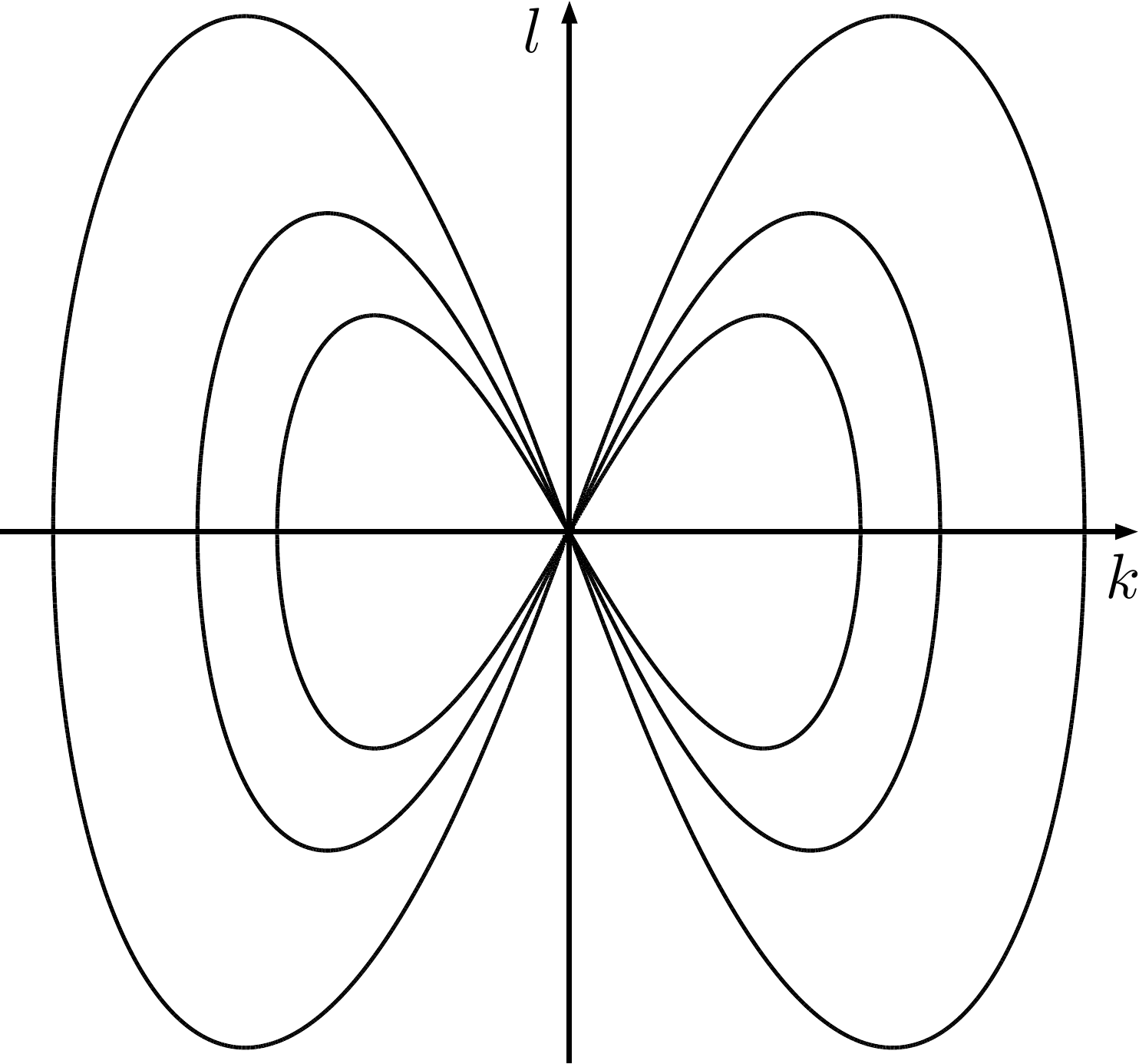}
  \caption{Illustration of the dispersion relation \eqref{eqn:erik:disprel}. The contours are the solutions to the dispersion relation in the $(k,l)$-plane for fixed phase speed $c_p$. From left to right: $\beta=0$, $0<\beta<1/3$ and $\beta>1/3$.}
  \label{fig:erik:disprel}
\end{figure}

In the irrotational case, $\boldsymbol{\omega}=0$, the equations can be simplified by the introduction of a velocity potential, just as in the two-dimensional setting. It is convenient to write $\mathbf{u}=\nabla \phi+(\mathbf{c},0)$, so that $\phi$ is the potential of the `absolute' rather than `relative' velocity field. In terms of $\phi$, the water wave problem reduces to the elliptic free boundary problem
\begin{subequations}
\label{eqn:erik:vp}
\begin{alignat}{2}
\Delta \phi &=0  &&\quad \text{in $\Omega$}, 	\label{eqn:erik:vp Laplace}\\
\phi_z&=0 &&\quad \text{on $B$},  \label{eqn:erik:vp kinematic bottom}\\
-\mathbf{c}\cdot \nabla \eta+\partial_\mathbf{n} \phi&=0  &&\quad \text{on $S$},  \\
(\mathbf{c},0)\cdot \nabla \phi+\frac12 |\nabla \phi|^2 + g \eta -\sigma  \nabla \cdot \left(\frac{\nabla \eta}{\sqrt{1+|\nabla \eta|^2}}\right) &=0 &&\quad \text{on $S$},
\end{alignat}
\end{subequations}
where $\mathbf{n} =(-\nabla \eta, 1)$ is a non-unit exterior normal vector at the surface. The problem can be further simplified by introducing the trace $\Phi=\phi|_{z=\eta}$ of the velocity potential at the free surface as a new variable. We then get the steady version of the Zakharov-Craig-Sulem formulation
\begin{equation}
\label{eqn:erik:ZCS}
\begin{aligned}
\mathbf{c}\cdot \nabla \Phi+\frac12 |\nabla \Phi|^2 -\frac{(\nabla \eta \cdot \nabla \Phi+G(\eta)\Phi)^2}{2(1+|\nabla \eta|^2)}+g\eta-\sigma \nabla \cdot \left(\frac{\nabla \eta}{\sqrt{1+|\nabla \eta|^2}}\right)&=0,\\
-\mathbf{c}\cdot \nabla \eta+G(\eta)\Phi&=0,
\end{aligned}
\end{equation}
where the Dirichlet-Neumann operator is defined by solving \eqref{eqn:erik:vp Laplace}, \eqref{eqn:erik:vp kinematic bottom} with the Dirichlet condition $\phi|_{S}=\Phi$ and setting $G(\eta)\Phi=\partial_{\mathbf{n}} \phi|_{S}$. Linearising and using the relation $G(0)e^{i\mathbf{k}\cdot \mathbf{x}'}=|\mathbf{k}| \tanh(|\mathbf{k}|d)e^{i\mathbf{k}\cdot \mathbf{x}'}$, we obtain the dispersion relation\footnote{There is a trivial element of the kernel with $(\eta,\Phi)=(0,1)$ corresponding to $\mathbf{k}=0$, which is due to the fact that $\Phi$ is only determined up to a constant. This is eliminated below in the periodic case by requiring that $\Phi$ has zero average and in the solitary case by working in homogeneous Sobolev spaces.},
\[
D_0(\mathbf{c}, \mathbf{k})=g + \sigma |\mathbf{k}|^2-\frac{(\mathbf{c} \cdot \mathbf{k})^2}{|\mathbf{k}|}  \coth (|\mathbf{k}|d)  = 0, \quad \mathbf{k}\ne 0.
\]
In the special case when $\mathbf{c}=c \mathbf{e}_1$ (so that the waves travel in the $x$-direction), we can solve for $c$ and obtain the expression
\begin{equation}
\label{eqn:erik:disprel}
c_p^2(\mathbf{k})=\frac{(g+\sigma|\mathbf{k}|^2) |\mathbf{k}| \tanh(|\mathbf{k}|d)}{k^2}, \quad \mathbf{k}=(k,l),
\end{equation}
for the linear phase speed, and when $l=0$ we recover the two-dimensional phase speed
\[
c_p(k)=\sqrt{\frac{(g+\sigma k^2)\tanh(k d)}{k}}.
\]
The dispersion relation is illustrated in Figure \ref{fig:erik:disprel}, by showing contour plots of $c_p(\mathbf{k})$ for different values of the nondimensional parameter $\beta\coloneqq \sigma/(gd^2)$.

\subsubsection{Doubly periodic waves}
\label{erik:sec:doubly periodic}
\begin{figure}
  \centering
  \includegraphics[width=0.4\linewidth]{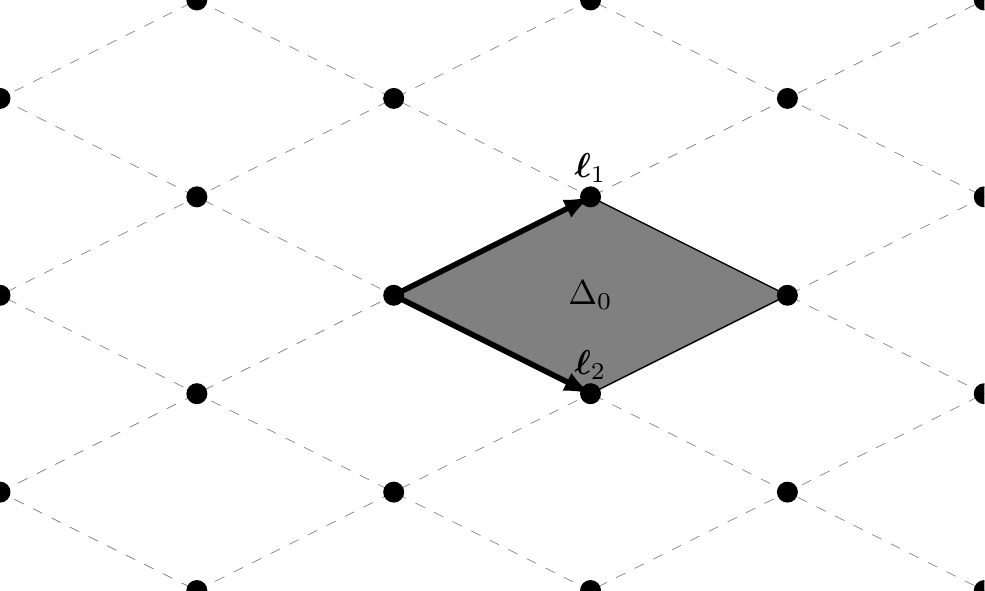} \hspace{1cm}
  \includegraphics[width=0.4\linewidth]{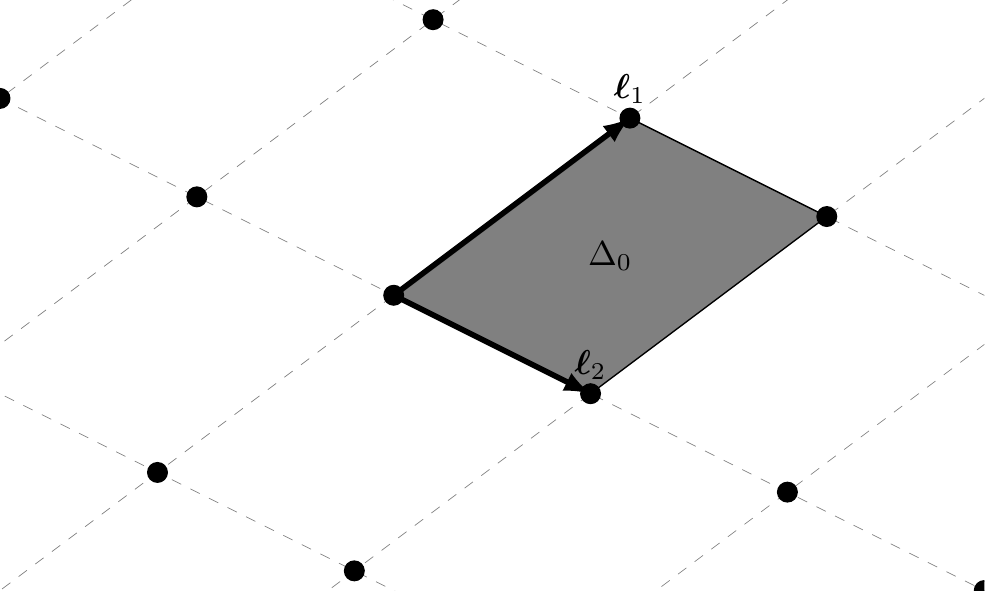}
  \caption{A symmetric lattice (left) and an asymmetric lattice (right). The shaded region is the fundamental domain $\Delta_0$.}
  \label{fig:erik:lattice}
\end{figure}

We now specialise to waves which are periodic in two linearly independent directions $\boldsymbol{\ell}_1, \boldsymbol{\ell}_2\in \mathbb{R}^2$, and hence with respect to the whole lattice
$\Lambda\coloneqq \{m_1\boldsymbol{\ell}_1+m_2 \boldsymbol{\ell}_2 \colon m_1, m_2 \in \mathbb{Z}\}\subset \mathbb{R}^2$ (Figure \ref{fig:erik:lattice}).
Then $\eta$ can be expanded as a Fourier series $\eta=\sum_{\mathbf{k}\in \Lambda'} \hat \eta_{\mathbf{k}} e^{i \mathbf{k} \cdot \mathbf{x}'}$ over the dual lattice
$\Lambda'= \{n_1 \mathbf{k}_1+n_2 \mathbf{k}_2 \colon n_1, n_2 \in \mathbb{Z}\}$ where $\mathbf{k}_i \cdot \boldsymbol{\ell}_j=2\pi \delta_{ij}$, $i, j=1,2$ (with similar expansions of the other variables). Note that the dimension of the kernel of the linearised problem for a given velocity vector $\mathbf{c}$ is given by the number of solutions to the equation $D_0(\mathbf{c}, \mathbf{k})=0$ in $\Lambda'\setminus\{0\}$. In particular, in the gravity-capillary case, we note that the number of solutions is finite since $D_0 \sim \sigma |\mathbf{k}|^2$ as $|\mathbf{k}|\to \infty$ (see also Figure \ref{fig:erik:disprel} (middle) and (right)).

While formal expansions of 2D periodic waves go back to Stokes in the mid 19th century, similar expansions of doubly periodic waves (often called `short crested waves' in the applied literature) were carried out first in the 1950's  by Fuchs \cite{Fuchs52} and Sretenskii \cite{Sretenskii53}. See Dias \& Kharif \cite[Section 6]{DiasKharif99} for more on formal expansions and numerical computations in the 20th century.

\begin{figure}
  \centering
  \includegraphics[width=0.4\linewidth]{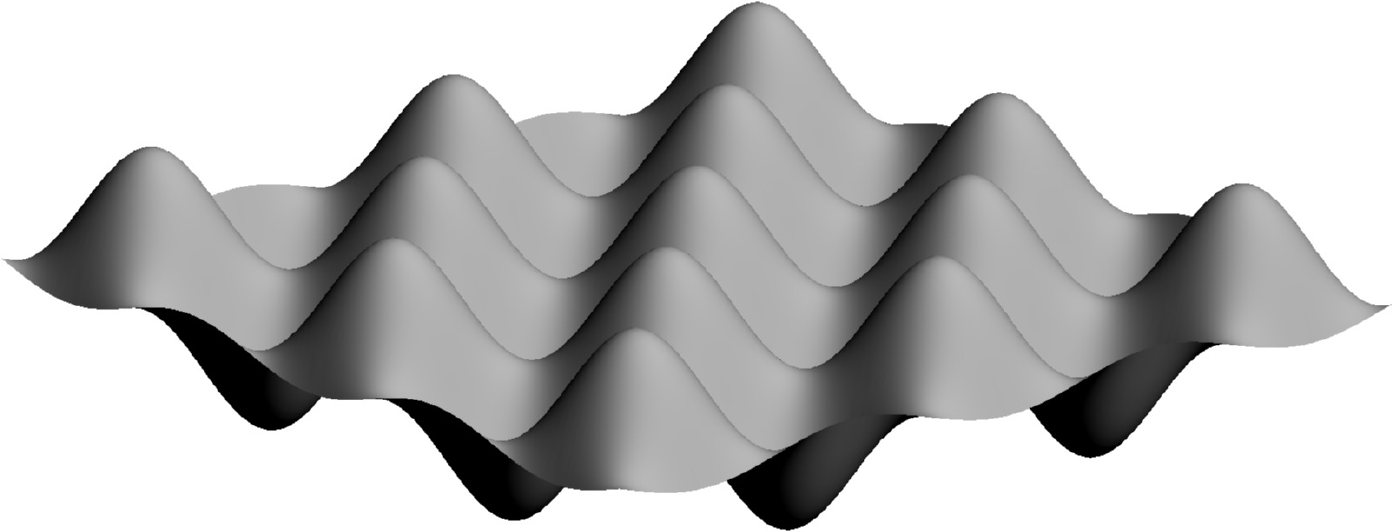} \hspace{1cm}
  \includegraphics[width=0.4\linewidth]{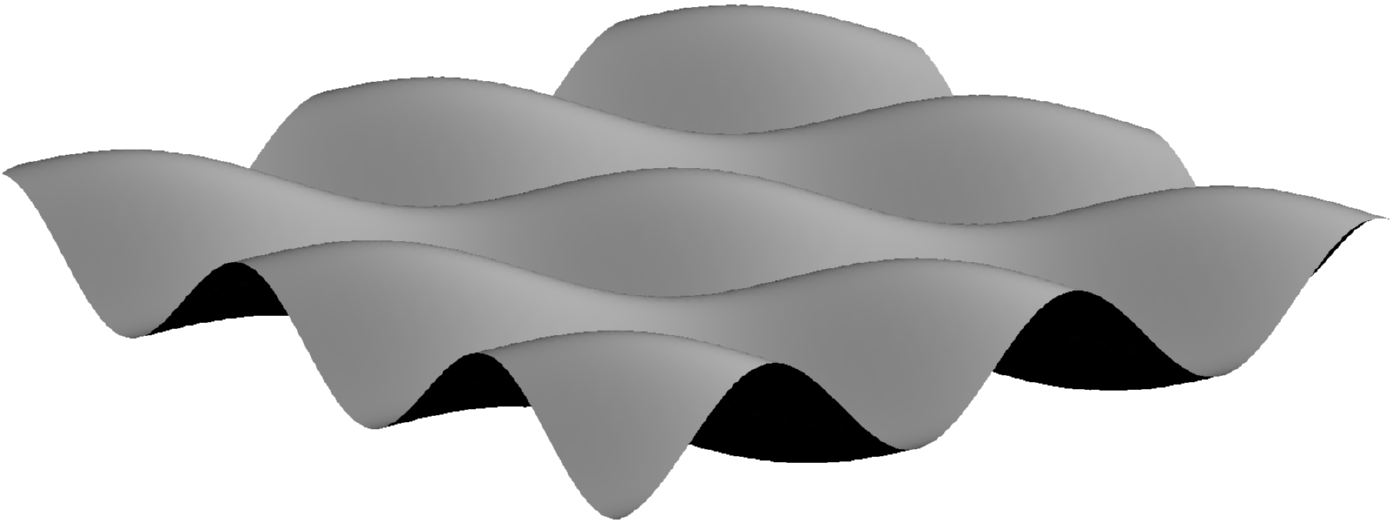}
  \caption{A small-amplitude symmetric doubly periodic wave (left) and an asymmetric doubly periodic wave (right).}
  \label{fig:erik:dp}
\end{figure}

The first rigorous construction of doubly periodic gravity-capillary waves is due to Reeder \& Shinbrot \cite{ReederShinbrot81}, who considered periodic lattices for which the fundamental domain is a `symmetric diamond' (see Figure \ref{fig:erik:lattice}). This means that $|\boldsymbol{\ell}_1|=|\boldsymbol{\ell}_2|$, and similarly for $\mathbf{k}_1$ and $\mathbf{k}_2$. After a rotation, we can assume that $\mathbf{k}_1=(\kappa_1, \kappa_2)$ and $\mathbf{k}_2=(\kappa_1, -\kappa_2)$ and take $\mathbf{c}=c\mathbf{e}_1$. Choosing the parameters and lattice appropriately, one can guarantee that $\pm \mathbf{k}_1$ and  $\pm \mathbf{k}_2$ are the only solutions to the dispersion relation. Hence the dimension of the kernel of the linearised problems is four, which is reduced to one by imposing the symmetries $ST_1$ and $T_2$, meaning in particular that $\eta$ is even in $x$ and $y$. This now allows a standard use of the Crandall-Rabinowitz local bifurcation theorem, although the proof in Reeder \& Shinbrot \cite{ReederShinbrot81} is formulated in different terms.
\begin{theorem}[\cite{ReederShinbrot81}]
\label{thm:erik:Reeder-Shinbrot}
Let $\lambda=\frac{1}{|\mathbf{k}_1|}\sqrt{\frac{g}{\sigma}}$ and $\mathbf{k}_1=|\mathbf{k}_1|(\cos\theta \sin \theta)$.
Assume that $(\lambda, \theta) \not \in M_d$, where $M_d\subset (0,\infty)\times (0, \pi/2) $ is a certain nowhere dense union of a countable collection of monotone curves. Then there exists a family of $\Lambda$-periodic and symmetric, smooth solutions $(\eta, \Phi, c)(\varepsilon)$ bifurcating at $(0, 0, c^\star)$, where $c^\star=c_p(\mathbf{k}_1)=c_p(\mathbf{k}_2)$.
The solutions are to leading order of the form
\[
\eta=\varepsilon \cos(\kappa_1 x)\cos(\kappa_2 y)+O(\varepsilon^2).
\]
\end{theorem}
The solutions are illustrated in Figure \ref{fig:erik:dp} (left). The set $M_d$ consists of parameters for which there are more solutions to the dispersion relation within the lattice $\Lambda'$ than $\pm \mathbf{k}_1$ and  $\pm \mathbf{k}_2$ at $c=c^\star$. For a more detailed description of this set we refer to \cite{ReederShinbrot81}. See also Sun \cite{Sun93} for a similar result  with pressure forcing.

Waves with an arbitrary fundamental domain were considered by Craig \& Nicholls \cite{CraigNicholls00} using a variational approach.
From the Hamiltonian formulation of the time-dependent water wave problem, it follows that travelling waves are critical points of the energy functional
\begin{align*}
\mathcal E(\eta, \Phi) &= \int_{\Delta_0} \left(\tfrac12 \Phi G(\eta)\Phi +\tfrac{1}2 g \eta^2 +\sigma (\sqrt{1+|\nabla \eta|^2}-1) \right) \,{\mathrm{d}\kern0.2pt} \mathbf{x}'\\
&= \int_{\Omega_0} \tfrac12  |\nabla \phi|^2\,  {\mathrm{d}\kern0.2pt} \mathbf{x} +\int_{\Delta_0} \left(\tfrac{1}2 g\eta^2 + \sigma (\sqrt{1+|\nabla \eta|^2}-1) \right) \, {\mathrm{d}\kern0.2pt} \mathbf{x}'
\end{align*}
subject to the constraints
\[
\mathcal I_j(\eta, \Phi)\coloneqq \int_{\Delta_0} \partial_j \eta \,\Phi \, {\mathrm{d}\kern0.2pt} \mathbf{x}'=-\int_{\Omega_0} \partial_j \phi\,  {\mathrm{d}\kern0.2pt} \mathbf{x}=2\mu_j,
\]
with the velocity vector $\mathbf{c}=(c_1, c_2)$ acting as a Lagrange multiplier. Here $\Delta_0$ is the fundamental domain of the lattice (see Figure \ref{fig:erik:lattice}) and $\Omega_0=\Omega \cap (\Delta_0 \times \mathbb{R})$ the corresponding part of the fluid domain.
This variational problem has a number of challenges. As always when dealing with water waves, one has to handle the unknown domain. This can be done either by using the formulation in terms of $\Phi$ rather than $\phi$, or by a transforming the fluid domain to a fixed domain. Secondly, the problem is quasilinear rather than semilinear, in the sense that the superquadratic part of the energy requires at least as much regularity as the quadratic part.  In fact, the quadratic part of the energy only controls the $H^1$ norm of $\eta$ (and a homogeneous $H^{1/2}$ norm of $\Phi$), whereas the Dirichlet-Neumann operator $G(\eta)$ requires $\eta$ to be at least in $W^{1,\infty}$. As a consequence, many variational methods fail to apply since they are developed for semilinear problems.
Instead, Craig \& Nicholls use a reduction method to first obtain a locally equivalent variational problem, which is more amenable to analysis. More specifically, they use a variational Lyapunov-Schmidt reduction near a velocity $\mathbf{c}^\star$ for which the kernel is nontrivial in order to reduce to a finite-dimensional problem. They work in the function space $X=H_{\text{per}}^{s+2}\times H_{\text{per},0}^{s+1}$, where  `per' stands for $\Lambda$-periodicity and $0$  for zero average, and $s>1$.
The idea is to split the Euler-Lagrange equation
\[
\mathcal J_{\mathbf{c}}'(u)\coloneqq \mathcal E'(\eta, \Phi)-c_1 \mathcal I_1'(\eta, \Phi)-c_2\mathcal I_2'(\eta, \Phi)=0,
\]
where we temporarily let $u=(\eta, \Phi)$, into the components
\begin{align*}
\Pi\mathcal J_{\mathbf{c}}'(v+w)&=0,\\
(I-\Pi)\mathcal J_{\mathbf{c}}'(v+w)&=0,
\end{align*}
in which $\Pi$ denotes $L^2$ orthogonal projection onto the kernel of the linearised operator $\mathcal J_{\mathbf{c}^\star}''(0) \colon X \to Y$, $Y=H^s_\text{per}\times H^s_{\text{per},0}$,  and $u=v+w$ with $v$ in the kernel $X_1$ and $w$ in the orthogonal complement $X_2$. Observe that $(I-\Pi)$ projects $Y$ onto the range $Y_1$ of $\mathcal J_{\mathbf{c}^\star}''(0)$. The second equation can be solved for $w=w(v, \mathbf{c})$ when $\mathbf{c}$ is close to $\mathbf{c}^\star$ and $v$ is small. Substituting this in the first equation, we obtain the reduced, finite-dimensional problem
\[
\Pi \mathcal J_{\mathbf{c}}'(v+w(v, \mathbf{c}))=0,
\]
which is equivalent to saying that $v$ is a critical point of the functional
\[
\tilde{\mathcal J}_{\mathbf{c}}(v)= \mathcal J_{\mathbf{c}}(v+w(v, \mathbf{c}))
\]
since
\[
\langle \tilde{\mathcal J}_\mathbf{c}'(v), \delta v\rangle= \langle \mathcal J'_{\mathbf{c}}(v+w(v, \mathbf{c})), \delta v\rangle
+\langle \mathcal J'_{\mathbf{c}}(v+w(v, \mathbf{c})), \mathrm{d}_v w(v, \mathbf{c}) \delta v \rangle
\]
and the second term vanishes since $\mathcal J'_{\mathbf{c}}(v+w(v, \mathbf{c})) \in X_1$ by choice of $w$ and since $\mathrm{d}_v w(v, \mathbf{c})\delta v \in X_2$.
Craig \& Nicholls show that it is in fact possible to choose $\mathbf{c}=\mathbf{c}(v)$ in such a way that $v$ is a critical point of the functional $\tilde {\mathcal E}(v)\coloneqq \mathcal E(v+w(v, \mathbf{c}(v)))$ subject to the constraints $\tilde {\mathcal I}_j(v)\coloneqq \mathcal I_j(v+w(v, \mathbf{c}(v)))=2\mu_j$, $j=1,2$. Furthermore, the original functionals are invariant under a `torus action' $u\mapsto u(\cdot+\boldsymbol{\xi})$, $\boldsymbol{\xi}\in \mathbb{R}^2 / \Lambda$, and the reduction can be done in such a way that the reduced functionals are also invariant. In the case that the kernel is four-dimensional, we obtain a two-dimensional constraint manifold $S(\boldsymbol{\mu})$ on which the reduced energy functional is constant due to the torus action (under a certain nondegeneracy condition). Thus each point on the constraint manifold is critical and gives rise to a solution to the original gravity-capillary problem.
\begin{theorem}[\cite{CraigNicholls00}]
\label{thm:erik:Craig-Nicholls}
Assume that the kernel is four-dimensional with $\eta$-component generated by the functions
$\{\cos(\mathbf{k}_1 \cdot \mathbf{x}'), \sin(\mathbf{k}_1 \cdot \mathbf{x}'), \cos(\mathbf{k}_2 \cdot \mathbf{x}'), \sin(\mathbf{k}_2 \cdot \mathbf{x}')\}$.
Then for each sufficiently small $\boldsymbol{\mu}$ which is not collinear with either $\mathbf{k}_1$ or $\mathbf{k}_2$, the constraint manifold $S(\boldsymbol{\mu})$ corresponds to a nontrivial $\Lambda$-periodic solution of  \eqref{eqn:erik:ZCS} and its translates.
\end{theorem}
By varying $\boldsymbol{\mu}$ we thus obtain a two-dimensional family of geometrically distinct solutions which to leading order are linear combinations of the trigonometric functions in the theorem; see Figure \ref{fig:erik:dp} (right) for an illustration.
Using index theoretical-methods, they similarly obtain at a point where the kernel is $2N$ dimensional $N-1$ geometrically distinct solutions for each $\boldsymbol{\mu}$ (under some nondegeneracy conditions).

Existence results for doubly periodic waves can also be obtained using  spatial-dynamics methods. Symmetric diamond waves were constructed by Groves \& Mielke \cite{GrovesMielke01}, while Groves \& Haragus \cite{GrovesHaragus03} considered arbitrary fundamental domains. Nilsson \cite{Nilsson19b} considered the same question for internal waves and corrected some mistakes in \cite{GrovesHaragus03}.

The case of pure gravity waves ($\sigma=0$) is much more challenging due to the appearance of small divisors.  Indeed, taking $\mathbf{c}=c\mathbf{e}_1$, this can be seen by noting that
\[
D_0(c\mathbf{e}_1, \mathbf{k})=
g-c^2 \frac{k^2}{|\mathbf{k}|} \coth(|\mathbf{k}|d) \underset{|\mathbf{k}| \to \infty}{\approx} g -c^2\frac{k^2}{|\mathbf{k}|},
\]
where $\mathbf{k}=(k,l)$.
In 2D (i.e.~for $l=0$), this behaves like $-c^2|k|$ for large $k$, but in 3D it can be can be made $O(|\mathbf{k}|^{-1/2})$ by choosing $\mathbf{k}\in \Lambda'$ appropriately (see e.g.~\cite[Theorem 2.7]{CraigNicholls00}). This can also be seen graphically by noting that the different contours in Figure \ref{fig:erik:disprel} (left) will come  close to the lattice points infinitely many times as $\mathbf{k}\to  \infty$. These small divisors signify that  bifurcation methods based on the standard implicit function theorem will fail, and must be replaced by an approach based on Nash-Moser theory.
This monumental task was first carried out by Iooss \& Plotnikov in the symmetric diamond case in \cite{IoossPlotnikov09} and then for a  general fundamental domain in \cite{IoossPlotnikov11} (both for infinite depth). Both works are based on the Zakharov-Craig-Sulem formulation and the main task is to show that the linearised operator around a {\em nontrivial} but  small  state $(\eta, \Phi)\in C^\infty$ can be inverted with a loss of derivatives.
Note that this operator is a pseudodifferential operator with variable coefficients. In both works, this  is accomplished by the {\em method of descent}, whereby the operator is transformed to one with constant coefficients plus a sufficiently smooth remainder, and by showing that the the transformed operator can be inverted with an acceptable loss of derivatives for a  large set of parameter values. The additional challenge in the asymmetric case is mainly that the change of variables is implicit and has to be incorporated into the Nash-Moser iteration.
In order to describe the main result, let $\mathbf{k}_j=|\mathbf{k}_j|(\cos \theta_j, \sin \theta_j)$, $j=1,2$, and $\tau_j=\tan \theta_j$.

\begin{theorem}[\cite{IoossPlotnikov11, IoossPlotnikov09}]
\label{thm:erik:Iooss-Plotnikov}
Choose integers $m\ge 34$ and $N\ge 2$ (even) and a real number $0<\delta<1$. Then there is a full measure subset $T\subset (\mathbb{R}^+)^2$ such that for any $\boldsymbol{\tau} \in T$ there exists a subset $E(\boldsymbol{\tau})$ of the quadrant $\{(\varepsilon_1^2, \varepsilon_2^2)\subset (\mathbb{R}^+)^2\}$ for which zero is a Lebesgue point and such that for any $\delta<\varepsilon_1/\varepsilon_2<1/\delta$ and $\boldsymbol{\varepsilon}\in E(\tau)$,  problem \eqref{eqn:erik:ZCS} with $\sigma=0$ and $d=\infty$ has a  weakly symmetric solution of the form $(\eta, \Phi, \mathbf{c})=(\eta_{2N}(\boldsymbol{\varepsilon}), \Phi_{2N}(\boldsymbol{\varepsilon}), \mathbf{c}(\boldsymbol{\varepsilon}))+|\boldsymbol{\varepsilon}|^N (\tilde \eta(\boldsymbol{\varepsilon}), \tilde \Phi(\boldsymbol{\varepsilon}), \tilde{\mathbf{c}}(\boldsymbol{\varepsilon}))$, where $(\eta_{2N}(\boldsymbol{\varepsilon}), \Phi_{2N}(\boldsymbol{\varepsilon}), \mathbf{c}(\boldsymbol{\varepsilon}))$ is a certain approximate solution with precision $|\boldsymbol{\varepsilon}|^{2N+1}$ and $(\tilde \eta(\boldsymbol{\varepsilon}), \tilde \Phi(\boldsymbol{\varepsilon}), \tilde{\mathbf{c}}(\boldsymbol{\varepsilon})) \in H_\text{per}^m\times H_{\text{per}, 0}^m\times \mathbb{R}^2$.
\end{theorem}

Here weakly symmetric means that $(\eta(-\mathbf{x}'), \Phi(-\mathbf{x}'))=(\eta(\mathbf{x}'), -\Phi(\mathbf{x}'))$ and one can check that this agrees with invariance under $ST_1T_2$. The symmetric diamond case is when $\tau_1=\tau_2$ and $\varepsilon_1=\varepsilon_2$, in which case the solution is (strongly) symmetric. The approximate solution exists whenever the dispersion relation only has the solutions $\pm \mathbf{k}_1$ and $\pm \mathbf{k}_2$ within the lattice $\Lambda'$ (which is true in the set $T$), and it has the form
\[
\eta_{2N}(\boldsymbol{\varepsilon})=\sum_{p+q=1}^{2N} \varepsilon_1^p \varepsilon_2^q \eta_{pq}, \quad
\Phi_{2N}(\boldsymbol{\varepsilon})=\sum_{p+q=1}^{2N} \varepsilon_1^p \varepsilon_2^q \Phi_{pq}, \quad
\mathbf{c}_{2N}(\boldsymbol{\varepsilon})=  (c_0,0) +\sum_{p+q=1}^{N} \varepsilon_1^{2p} \varepsilon_2^{2q} \mathbf{c}_{pq},
\]
with $c_0^2=g/(|\mathbf{k}_1|\cos^2(\theta_1))=g/(|\mathbf{k}_2|\cos^2(\theta_2))$ and $\eta_{10}=\cos(\mathbf{k}_1\cdot \mathbf{x}')$, $\eta_{01}=\cos(\mathbf{k}_2\cdot \mathbf{x}')$.

Many important questions remain about doubly periodic waves, perhaps the most important one being large-amplitude theory with or without surface tension. Formal expansions and numerical studies shed some light on what can be expected, although
computations in 3D are in general very expensive \cite{AkersReeger17, NichollsReitich06}. In particular, we note that Akers \& Reeger \cite{AkersReeger17} computed overhanging gravity-capillary waves (albeit with small negative gravity), while Craig \& Nicholls  \cite{CraigNicholls02} (among others) computed steep gravity waves.

\subsection{Fully localised 3D solitary waves}

We next turn our attention to {\em fully localised solitary waves}, that is, solutions which decay in all horizontal directions, $\eta(\mathbf{x}')\to 0$ as $|\mathbf{x}'|\to \infty$. In 2D, spatial dynamics has been a successful strategy for constructing small-amplitude solitary waves (see Section \ref{sec:miles}). However, this does not work as well in 3D since these methods typically only allow for one unbounded variable (see however \cite{HillLoydTurner21} for a spatial dynamics approach to localised radial patterns in the context of free surface ferrofluids). Instead, the constructions of fully localised solitary waves  so far rely on either variational methods or versions of the implicit function theorem.

\begin{figure}
  \centering
  \includegraphics[width=0.4\linewidth]{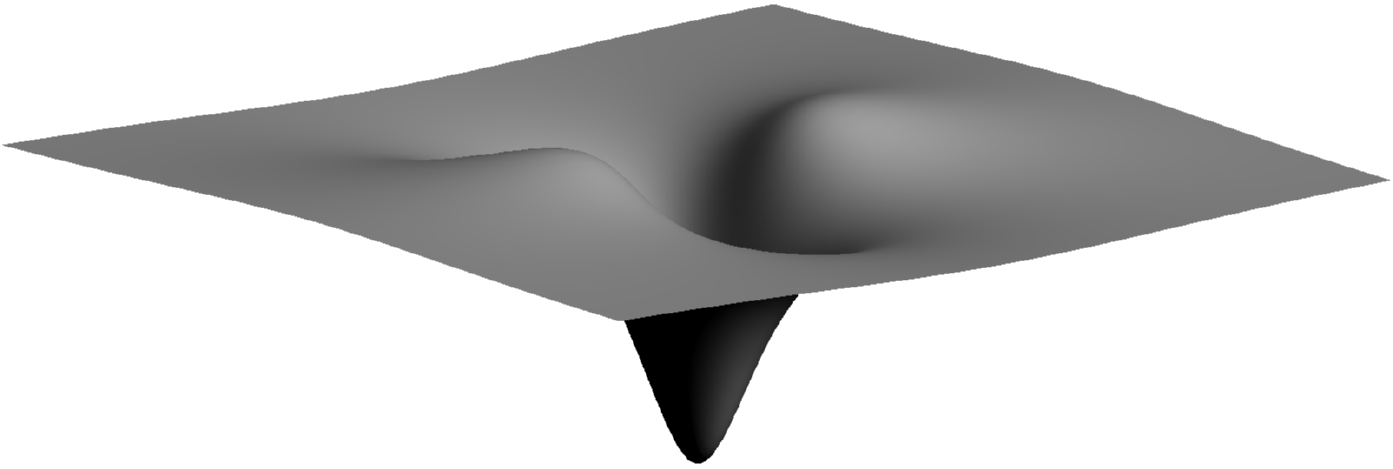}
  \caption{Fully localised solitary waves in the KP-I approximation.}
  \label{fig:erik:KPfloc}
\end{figure}

Before going into details we first look at some weakly nonlinear model equations. There are two main models. The first is the Kadomtsev-Petviashvili (KP) equation \cite{KadomtsevPetviashvili70}
\[
\left(-\zeta_T+\tfrac32 \zeta\zeta_X+\tfrac12\left(\tfrac13-\beta\right) \zeta_{XXX}\right)_X+\zeta_{YY}=0,
\]
which arises by a small-amplitude and anisotropic long-wave scaling
 \[
\eta=\varepsilon^2 d \zeta\left(\varepsilon d^{-1}(x+\sqrt{gd}\, t), \varepsilon^2d^{-1} y, \varepsilon^3 (g/d)^{1/2} t\right)+\mathcal{O}(\varepsilon^3)
\]
where $0<\varepsilon\ll 1$ (see e.g.~\cite{AblowitzSegur79, Lannes13}).
Fully localised solitary waves are know to exist for $\beta\coloneqq \sigma/(gd^2) >1/3$ (KP-I). In fact, there's an explicit family of solutions with algebraic decay, first found by Manakov et al.  \cite{ManakovZakharovBordagItsMatveev77}; see Figure \ref{fig:erik:KPfloc}.  Recently, the lumps have been shown to be nondegenerate and orbitally stable by Liu \& Wei \cite{LiuWei19}. However, uniqueness is currently unknown. Fully localised solitary waves to KP-I  have also been constructed as ground states for a variational problem by de Bouard \& Saut \cite{DeBouardSaut97} and Wang, Ablowitz \& Segur \cite{WangAblowitzSegur94}  and the set of ground states has been shown to be stable by de Bouard \& Saut \cite{DeBouardSaut96} and Liu \& Wang \cite{LiuWang97}; note that fully localised solitary waves $\zeta=\zeta(X+cT, Y)$ to the KP equation are formally critical points of the functional
\[
\int_{\mathbb{R}^2} \left(-\tfrac12 c\zeta^2+\tfrac14 \left(\beta-\tfrac13\right)\zeta_X^2+\tfrac12(\partial_X^{-1}\partial_Y \zeta)^2 +\tfrac14 \zeta^3\right)\, {\mathrm{d}\kern0.2pt}X\, {\mathrm{d}\kern0.2pt}Y.
\]
It is  unknown if the explicit lump solutions belong to the family of ground states.
For $\beta<1/3$ (KP-II) one does not expect to find solitary waves due to the fact that the quadratic part of the above functional is not positive definite. This was  shown rigorously in the case $c>0$ by de Bouard \& Saut \cite{DeBouardSaut97}.
Near $\beta=1/3$ the KP equation is replaced by a fifth order version \cite{Lannes13}.
For a nice review of the KP equation containing a much more detailed discussion, we refer to \cite{KleinSaut12}.

\begin{figure}
  \centering
  \includegraphics[width=\linewidth]{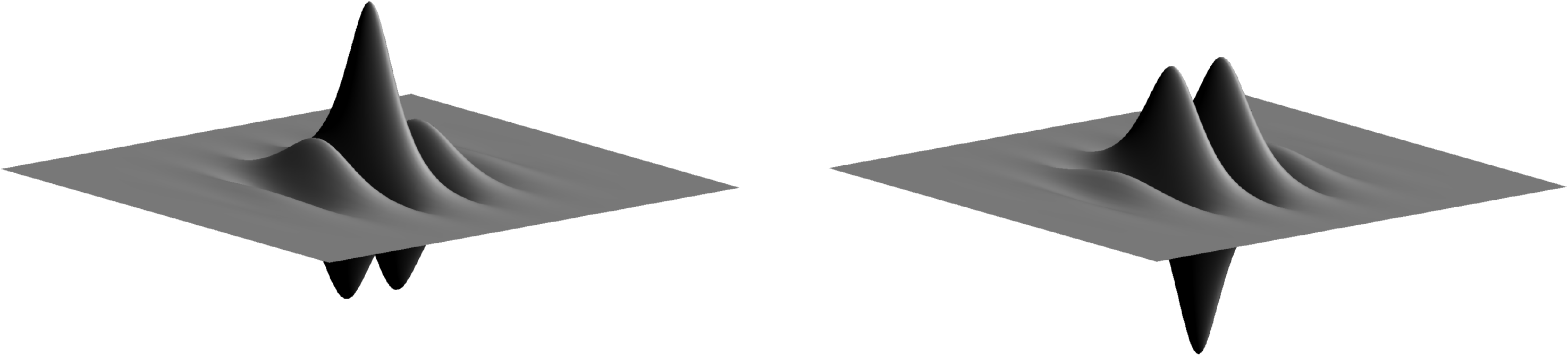}
  \caption{Fully localised solitary waves modelled by the Davey-Stewartson or nonlinear Schrödinger equations.}
  \label{fig:erik:DSfloc}
\end{figure}

The second main model equation is the Davey-Stewartson system, which has the form
\begin{align*}
i\zeta_T+ a_1\zeta _{XX}+a_2\zeta_{YY}&=- a_3|\zeta|^2\zeta-a_4\zeta \psi_X,\\
 a_5\psi_{XX}+\psi_{YY}&=(|\zeta|^2)_X,
\end{align*}
for certain real coefficients $a_1, \ldots, a_5$, and which arises from the modulation Ansatz
\begin{equation}
\label{eqn:erik:modulational}
\eta=\varepsilon d \zeta(\varepsilon d^{-1} (x+c_g t) , \varepsilon d^{-1} y, \varepsilon^2 (g/d)^{1/2} t)e^{ik_0 (x-c_p t)}+\text{c.c.}+\mathcal{O}(\varepsilon^2),
\end{equation}
where $c_p=c_p(k_0)$, $c_g=\omega'(k_0)$, is the corresponding group speed (with $\omega(k)=k c_p(k)$) and `c.c' stands for complex conjugate.
It was derived by Davey \& Stewartson \cite{DaveyStewartson74} and Benney \& Roskes \cite{BenneyRoskes69} in the absence of surface tension, which was later included by Djordjevic \& Redekopp \cite{DjordjevicRedekopp77};
see also \cite{AblowitzSegur79, Lannes13}. Note that $c_p\ne c_g$ in general, meaning that the envelope and carrier wave travel at different speeds and that the above formula doesn't represent a travelling wave.
The system can also be rewritten as a nonlocal version of the 2D cubic nonlinear Schrödinger equation,
\[
i\zeta_T+ a_1\zeta _{XX}+a_2\zeta_{YY}=- a_3|\zeta|^2\zeta-a_4\zeta E(|\zeta|^2),
\]
where
\[
\mathcal{F}(E(f))(\mathbf{k})=\frac{k^2}{a_5k^2+l^2}\hat f(\mathbf{k}).
\]
In the limit of infinite depth\footnote{This requires a different nondimensionalisation in \eqref{eqn:erik:modulational}; see Theorem \ref{thm:erik:floc infinite depth} below.}, the coefficient $a_4\to 0$ and one obtains the usual 2D cubic nonlinear Schrödinger equation.
Fully localised standing wave solutions $\zeta= e^{i\gamma T} \varphi(X,Y)$, $\psi=\psi(X,Y)$, to the Davey-Stewartson equation were constructed by variational methods  by Cipolatti \cite{Cipolatti92} in the case when all $a_j>0$ and $\gamma>0$, noting that $\varphi$ is a critical point of the functional
\begin{equation}
\label{eqn:erik:DS functional}
\int_{\mathbb{R}^2}\left(\tfrac{1}{2}a_1|\varphi_X|^2+\tfrac{1}{2}a_2|\varphi_Y|^2+\tfrac{1}{2} \gamma|\varphi|^2-\tfrac{1}4 a_3 |\varphi|^4-\tfrac{1}4 a_4 |\varphi|^2 E(|\varphi|^2)\right)\, {\mathrm{d}\kern0.2pt}X\, {\mathrm{d}\kern0.2pt}Y;
\end{equation}
see also \cite{PapanicolaouSulemSulemWang94} for an alternative existence proof.
These standing waves are however unstable due to finite-time blow-up of nearby solutions \cite{Cipolatti93, PapanicolaouSulemSulemWang94}.
The positivity conditions on the coefficients are in particular satisfied for $0<\beta<1/3$ at the $k_0$ which minimises the phase speed $c_p(k)$, in which case $c_p=c_g$ and we  expect to find fully localised solitary waves.
For the NLS equation, the existence of a positive ground state solution which is exponentially decaying and radial (after a rescaling of $X$ or $Y$) is classical; see e.g.~references in \cite{ChangGustafsonNakanishiTsai07}.
The ground state solution is known to be unique \cite{Kwong89} and suffers the same kind of instability as in the Davey-Stewartson case \cite{Weinstein83}. The fully localised waves in this approximation are illustrated in Figure~\ref{fig:erik:DSfloc}.

We now turn back to the irrotational water wave problem. First note that the non-existence of fully localised solitary waves for KP-II indicates that such solutions might not exist for the full water wave problem either. A first result in this direction was proved by Craig \cite{craig:nonexistence} who showed that
if $\sigma=0$, there are no fully localised solitary waves with $\eta$ of constant sign in 3D.
Recall from Section \ref{sec:miles} that in 2D, gravity solitary waves are positive {\em a priori}.
We also refer to \cite{chen2019existence, Wheeler18a} for similar results in infinite depth and with surface tension and/or localised vorticity. This still leaves the door open for sign-changing solutions.
We also note that, under certain conditions, fully localised solitary waves are ruled out by recent global well-posedness results for 3D water waves, since they imply that the solutions decay uniformly as $t\to \infty$; see \cite{Wang19} for 3D gravity waves on finite depth and  \cite{GermainMasmoudiShatah12, Wu11} for similar results on infinite depth. Note however that the conditions are quite stringent since they typically both require the waves to  decay sufficiently rapidly and to be sufficiently small.
In fact, there is a similar result for gravity-capillary waves on infinite depth in 3D \cite{DengIonescuPausaderPusateri17} even though fully localised waves exist in that case (see Theorem \ref{thm:erik:floc infinite depth} below). The results are consistent, however, since the $L^2$ norm of these waves does not vanish when the amplitude tends to $0$.

We now turn to existence results, taking $\mathbf{c}=c\mathbf{e}_1$ without loss of generality.
As in the periodic case, fully localised solitary waves are critical points of the energy
\begin{align*}
\mathcal E(\eta, \Phi)&= \int_{\mathbb{R}^2} \left(\tfrac12 \Phi G(\eta)\Phi +\tfrac{1}2 g \eta^2 +\sigma (\sqrt{1+|\nabla \eta|^2}-1) \right) \,{\mathrm{d}\kern0.2pt} \mathbf{x}'\\
& =      \int_{\Omega} \tfrac12  |\nabla \phi|^2\,  {\mathrm{d}\kern0.2pt} \mathbf{x} +\int_{\mathbb{R}^2} \left(\tfrac{1}2 g \eta^2 +\sigma (\sqrt{1+|\nabla \eta|^2}-1) \right) \, {\mathrm{d}\kern0.2pt} \mathbf{x}',
\end{align*}
subject to the constraint of fixed total  momentum in the $x$-direction,
\[
\mathcal I_1(\eta, \Phi)\coloneqq \int_{\mathbb{R}^2} \eta_x \Phi \, {\mathrm{d}\kern0.2pt} \mathbf{x}'=-\int_{\Omega} \phi_x\,  {\mathrm{d}\kern0.2pt} \mathbf{x}  =\text{const.}
\]
Equivalently, we can consider critical points of the augmented energy
\begin{align*}
 \mathcal E_c(\eta,\Phi)=\mathcal E(\eta, \Phi)-c\mathcal I_1(\eta, \Phi).
\end{align*}
In addition to the previously mentioned challenges with the variational formulation, there is now a problem with compactness due to the unbounded domain.

The first rigorous existence result is due to Groves \& Sun \cite{GrovesSun08} who proved the existence of fully localised solitary waves
with $\beta>1/3$ by considering the functional $\mathcal E_c$ in terms of $\eta$ and $\phi$ for $c$ close to the critical wave speed $\sqrt{gd}$.
The unknown fluid domain is mapped to the infinite `slab' $\mathbb{R}^2 \times (-1,0)$ by the transformation $z\mapsto (z-\eta)/(d+\eta)$. Note that the functional $\mathcal E_c$ is neither bounded from above nor from below, and therefore a natural idea is to try a mountain-pass approach. However, the quasilinear nature of the problem makes this challenging. To solve this issue, Groves \& Sun begin by first making a variational reduction of the problem in terms of the leading order part of the transformed velocity potential. This is similar to the approach by Craig \& Nicholls, with the difference that the reduced problem is still infinite-dimensional. However, it now has a semilinear structure, and is of mountain-pass type. Finally, the unbounded domain is dealt with by using the concentration-compactness principle. One technical aspect of the whole proof is that it relies on fractional Sobolev spaces based on $L^p$ with $p>2$.
\begin{theorem}[\cite{GrovesSun08}]
Fix $\beta>1/3$ and let $c_\varepsilon^2=(1-\varepsilon^2)gd$. For each $0<\varepsilon \ll 1$, there is a nontrivial symmetric fully localised solitary-wave solution of \eqref{eqn:erik:vp} with speed $c_\varepsilon$.
\end{theorem}

A few years later, Buffoni, Groves, Sun \& Wahl\'en \cite{BuffoniGrovesSunWahlen13} gave an alternative existence proof based on minimising the energy $\mathcal E(\eta, \Phi)$ subject to the constraint of fixed horizontal momentum $\mathcal I_1(\eta, \Phi)=2\mu$, with the wave speed $c$ appearing as a Lagrange multiplier. Since the energy and momentum are conserved quantities, this approach has the advantage of yielding stability of the set of minimisers, which is however conditional upon global well-posedness  (the perturbed solution remains close to the set of minimisers in the energy norm, {\em assuming} that it remains sufficiently bounded in a higher order Sobolev norm).
Another advantage compared to \cite{GrovesSun08} is that the proof uses $L^2$-based  Sobolev spaces.
Again, the quasilinear nature of the problem causes challenges, which can be overcome using a version of a penalisation approach due to Buffoni \cite{Buffoni04a}. Let us give a brief outline of the proof. We begin by fixing $0\ne \eta \in H^3(\mathbb{R}^2)$ and minimising the functional $\mathcal E(\eta, \cdot)$ over the set of $\Phi \in H_\star^{1/2}(\mathbb{R}^2)$ with $\mathcal I_1(\eta, \Phi)=2\mu$, where $H_\star^{1/2}(\mathbb{R}^2)$ is a certain homogeneous Sobolev space (see \cite{BuffoniGrovesSunWahlen13} for the precise definition). This problem has a unique global minimiser $\Phi_\eta=c G(\eta)^{-1}\eta_x$, where $c$ is an unknown Lagrange multiplier which depends on $\mu$ and $\eta$. We then minimise the functional
\begin{align*}
\mathcal J_\mu(\eta)\coloneqq \mathcal E(\eta, \Phi_\eta)=\mathcal K (\eta)+\frac{\mu^2}{\mathcal{L}(\eta)}
\end{align*}
over a punctured ball $B_R(0)\setminus \{0\}$ in $H^3(\mathbb{R}^2)$, where
\begin{align*}
\mathcal L(\eta)=\frac12 \int_{\mathbb{R}^2} \eta_x G(\eta)^{-1} \eta_x\, {\mathrm{d}\kern0.2pt} \mathbf{x}' \quad \text{and} \quad \mathcal K(\eta)=\int_{\mathbb{R}^2}\left\{\tfrac{1}2 g\eta^2+\sigma(\sqrt{1+|\nabla \eta|^2}-1) \right\}\,{\mathrm{d}\kern0.2pt}\mathbf{x}'.
\end{align*}
Since the functional $\mathcal J_\mu$ is not coercive on $H^3$, we add a penalisation term $\varrho(\|\eta\|_{H^3}^2)$ which vanishes for $\|\eta\|_{H^3}\le \tilde R$, with $\tilde R$ close to $R$, and blows up as $\|\eta\|_{H^3}\to R$. Ignoring problems of compactness, the penalised functional then automatically has a minimiser (note that $\mathcal J_\mu(\eta)\to \infty$ as $\|\eta\|_{H^3}\to 0$) and the trick is to show that the minimiser lies in the region   $\|\eta\|_{H^3}\le \tilde R$ unaffected by the penalisation. Using a test function inspired by the KP scaling, one can show that $\inf \mathcal J_{\mu, \varrho}(\eta)=O(\mu)$ and therefore the minimiser satisfies $\|\eta\|_{L^2}^2=O(\mu)$, which can be improved to $\|\eta\|_{H^3}^2=O(\mu)$ using elliptic regularity. This closes the argument, except for the problem with the lack of compactness due to the unbounded domain. Again, this is dealt with using concentration-compactness and we refer to  \cite{BuffoniGrovesSunWahlen13} for further details.
\begin{theorem}[\cite{BuffoniGrovesSunWahlen13}]
For any $R>0$ and $0<\mu\ll 1$, the value
\[
\inf\{ \mathcal H(\eta, \Phi)\colon \eta \in B_R(0)\subset H^3(\mathbb{R}^2),\ \Phi \in H_\star^{1/2}(\mathbb{R}^2),\ \mathcal I_1(\eta, \Phi)=2\mu\}
\]
is attained and any minimiser is a nontrivial fully localised solution to \eqref{eqn:erik:ZCS}.
Moreover, the set  of minimisers is conditionally stable.
\end{theorem}

In the case of small surface tension, $0<\beta<1/3$, we don't expect to see any small-amplitude solutions with the KP-scaling. Instead, we expect solutions resembling the modulational Ansatz \eqref{eqn:erik:modulational}. Note however that in this scaling the momentum is not small as the amplitude tends to $0$. In addition, evidence from the model equation \cite{Cipolatti93, PapanicolaouSulemSulemWang94} and from numerical simulations \cite{WangMilewski12} suggests that these solutions are unstable for small amplitudes. Therefore, we don't expect the approach used in \cite{BuffoniGrovesSunWahlen13} to work in this case. An existence result using a variational reduction approach was however proved by Buffoni, Groves \& Wahl\'en in \cite{BuffoniGrovesWahlen18}. The approach has  similarities to the one by Groves \& Sun \cite{GrovesSun08}, but contains some simplifications as well as some new ingredients due to the different scaling.
We fix $\beta \in (0, 1/3)$, so that the 2D linear phase speed $c_p(k)$ has a unique minimum $c_0=c_p(k_0)$ for some $k_0>0$ and set $c_\varepsilon^2=(1-\varepsilon^2)c_0^2$. Recall that at $k_0$, $c_g(k_0)=c_p(k_0)$.
Similar to  Buffoni, Groves, Sun \& Wahl\'en  \cite{BuffoniGrovesSunWahlen13}, we first eliminate $\Phi$ by noting that  $\mathcal E_{c_\varepsilon}(\eta, \Phi)$ has a unique critical point  $\Phi_\eta=c_\varepsilon G(\eta)^{-1}\eta_x$.
This time we end up with a functional of the form
\[
\mathcal J_\varepsilon (\eta)\coloneqq \mathcal E_{c_\varepsilon}(\eta, \Phi_\eta)=\mathcal K(\eta)-c_\varepsilon^2 \mathcal L(\eta),
\]
where $\mathcal K$ and $\mathcal L$ are as above. Again, this functional is quasilinear and we therefore perform a further reduction in order to obtain a semilinear problem. To motivate this reduction, note that we expect to find solitary waves resembling \eqref{eqn:erik:modulational} (with $c_g=c_p=c_0$), which in particular are concentrated in the Fourier plane around the points $\pm  (k_0,0)$. Therefore, we split
\[
		\eta=\eta_1+\eta_2,
\]
with $\eta_1=\chi(D)\eta$ and $\eta_2=(1-\chi(D))\eta$, where $\chi(\mathbf{k})$ is the  characteristic function of the set $|\mathbf{k}-(k_0,0)|\le \delta$ for a sufficiently small $\delta>0$. We similarly split the Euler-Lagrange equation into two parts,
\begin{align*}
\chi(D) \mathcal J_\varepsilon'(\eta_1+\eta_2)&=0,\\
(1-\chi(D)) \mathcal J_\varepsilon'(\eta_1+\eta_2)&=0.
\end{align*}
By the choice of $\chi$, it is not difficult to see that the second equation can be solved uniquely for $\eta_2$ when $\eta_1$ and $\eta_2$ are both assumed small. We substitute the solution $\eta_2(\eta_1)$ into the variational functional $\mathcal J_\varepsilon$ and obtain a reduced functional
\[
\tilde{\mathcal J}_\varepsilon(\eta_1)\coloneqq \mathcal J_\varepsilon(\eta_1+\eta_2(\eta_1)).
\]
Since $\hat \eta_1$ has support close to $\pm (k_0,0)$, it can be written uniquely as
\begin{equation}
\label{eqn:erik:eta_1}
\eta_1=\varepsilon \zeta(\varepsilon x, \varepsilon y) e^{ik_0x}+\varepsilon \overline{\zeta(\varepsilon x, \varepsilon y)}e^{-ik_0x}
\end{equation}
and  careful bookkeeping shows that
\begin{align*}
\mathcal T_\varepsilon(\zeta)&\coloneqq\varepsilon^{-2} \tilde{\mathcal J}_\varepsilon(\eta_1)\\
&=
\int_{\mathbb{R}^2}\left(\tfrac{1}{2}a_1|\varphi_X|^2+\tfrac{1}{2}a_2|\varphi_Y|^2+\tfrac{1}{2} \gamma|\varphi|^2-\tfrac{1}4 a_3 |\varphi|^4-\tfrac{1}4 a_4 |\varphi|^2 E(|\varphi|^2)\right)\, {\mathrm{d}\kern0.2pt}X\, {\mathrm{d}\kern0.2pt}Y\\
&\qquad+\mathcal{O}(\varepsilon^{1/2}\|\zeta\|_1^2),
\end{align*}
for suitable positive coefficients $a_1, \ldots, a_4$ and $\gamma$, where we note that $\mathcal T_0$ is the functional \eqref{eqn:erik:DS functional} considered by Cipolatti.
We study this functional in $B_R(0)\subset H_\varepsilon^1(\mathbb{R}^2)\coloneqq \chi(\varepsilon D)H^1(\mathbb{R}^2)$ and show using concentration-compactness methods that $\mathcal T_\varepsilon|_{N_\varepsilon}$, $0<\varepsilon \ll 1$, has a minimum, where
 \[
N_\varepsilon \coloneqq \{\zeta\in B_R(0)\setminus \{0\}\subset H_\varepsilon^1(\mathbb{R}^2) \colon \; \langle\mathcal{T}_\varepsilon'(\zeta), \zeta\rangle=0\}
\]
is the natural constraint set. A standard argument then shows that the Lagrange multiplier for the minimiser $\zeta$ of $\mathcal T_\varepsilon|_{N_\varepsilon}$ is zero, so that $\zeta$ is actually a critical point $\mathcal T_\varepsilon$ which gives rise to a fully localised solitary wave through formula \eqref{eqn:erik:eta_1} and the relation $\eta=\eta_1+\eta_2(\eta_1)$.
\begin{theorem}[\cite{BuffoniGrovesWahlen18}]
\label{thm:erik:BGW}
Fix $0<\beta<1/3$ and let $c_\varepsilon^2=(1-\varepsilon^2)c_0^2$. For each $0<\varepsilon\ll 1$, there exists a
nontrivial fully localised solitary-wave solution  of \eqref{eqn:erik:ZCS} with speed $c_\varepsilon$ corresponding to a
minimiser of $\mathcal T_\varepsilon$ on $N_\varepsilon$.
\end{theorem}

In fact, the same method can be used in the case $\beta>1/3$ to give an alternative existence proof to the ones in \cite{BuffoniGrovesSunWahlen13, GrovesSun08}. This has been carried out  by Ehrnstr\"om \& Groves \cite{EhrnstromGroves18} for the full-dispersion KP equation, which is a nonlocal model that can be seen as lying between the KP equation and the full water wave problem in terms of complexity.

Besides these variational approaches, it is also possible to prove existence using the implicit function theorem whenever the model equation has a nondegenerate fully localised solution. One such example is the gravity-capillary problem on infinite depth. The Davey-Stewartson equation then reduces to the 2D elliptic NLS equation, for which the unique ground state solution is known to be nondegenerate \cite{ChangGustafsonNakanishiTsai07, Kwong89, Weinstein85}.
\begin{theorem}[\cite{BuffoniGrovesWahlen21}]
\label{thm:erik:floc infinite depth}
Let $d=\infty$ and $c_\varepsilon^2=(1-\varepsilon^2)c_0^2$, where $c_0^2=2 \sqrt{g\sigma}$. For each sufficiently small value of $\varepsilon>0$ there exist two symmetric fully localised solutions of \eqref{eqn:erik:ZCS}  with speed $c_\varepsilon$ given by
\[
\eta(\mathbf{x}')=\pm \varepsilon  d_0 \zeta_0(\varepsilon d_0^{-1}  \mathbf{x}') \cos (d_0^{-1}x) +o(\varepsilon)
\]
uniformly over $\mathbf{x}'\in \mathbb{R}^2$, where $\zeta_0$ is the unique ground state of the the 2D nonlinear Schrödinger equation
\[
\zeta-\frac12 \zeta_{XX}-\zeta_{YY}-\frac{11}{16}|\zeta|^2\zeta=0
\]
and $d_0=\sqrt{\sigma/g}$.
\end{theorem}
The idea behind the proof is again to do a reduction to perturbation of the model equation. But instead of solving the reduced equation by variational methods, we solve it using a version of the implicit function theorem. An advantage compared to the variational existence proofs is that we get the existence of two different branches --- an elevation branch and a depression branch. This is predicted by numerical calculation for both finite and infinite depth \cite{KimAkylas05, ParauVandenBroeckCooker05a,WangMilewski12}.

There are many interesting open problems inspired in part by numerical simulations and results for model equations. Both for strong and weak surface tension (and infinite depth), there is so far no rigorous construction of large-amplitude solutions. In both these cases, there is also  evidence of other types of fully localised waves, corresponding to `excited states' of the model equations \cite{PelinovskyStepanyants93, Strauss77, WangMilewski12}, for which one could try to prove existence in the full water wave problem. It might also be possible to construct other small-amplitude solutions by `gluing' together distant copies of known fully localised solutions, as was  done in 2D \cite{BuffoniGroves99}. Finally, the stability of fully localised solitary waves is a very interesting topic. While both the model equations and numerical simulations suggest that all solutions with $0<\beta<1/3$ are unstable, there is evidence of a difference between the elevation and depression branches for larger amplitudes: the elevation waves remain unstable while the depression waves stabilise at some point \cite{WangMilewski12}. There is also numerical evidence that the instability might lead to the development of a breather solution, for which no rigorous existence theory is currently available.

\subsection{Waves with vorticity}

We now turn our attention to three-dimensional water waves over flows with non-zero vorticity. So far, there is very little theory for such waves, one of the main reasons being the lack of an analogue of the stream function formulation. The full steady Euler equations are of elliptic-hyperbolic type and most tools for proving existence of steady solutions to free-boundary problems (bifurcation theory, variational methods, spatial dynamics methods) are adapted to elliptic equations.
At a very basic level, if one linearises the water wave problem in the form \eqref{eqn:erik:ww problem} at a laminar flow ($\eta=0$, $\mathbf{u}=(U(z), V(z), 0)$), then the kernel of the linearised problem is infinite-dimensional due to the presence of many other nearby shear flows. This degeneracy is solved in 2D by prescribing the relation between the (scalar) vorticity and the stream function, but it is not immediately clear how to accomplish something similar in 3D.
Based on 2D theory (see Section \ref{sec:vera}) it might seem natural to assume constant vorticity, but it turns out that this essentially forces the solutions to be two-dimensional; see Wahl\'en \cite{Wahlen14} and an extension to non-steady waves by Martin \cite{Martin18}. We will  nevertheless present two different constructions of three-dimensional water waves with vorticity, both for doubly periodic solutions. The first is under the assumption that the (relative) velocity field is a Beltrami field, in which case the equations turn out to be elliptic and some methods  for irrotational waves can be adapted. The second is a construction of symmetric diamond waves under the assumption that the vorticity is small. It relies on an approach by Lortz for constructing magnetohydrostatic equilibria \cite{Lortz70}.
Before going into details, we make some general remarks about steady solutions to the Euler equations in three dimensions.

By the identity
\[
(\mathbf{u} \cdot \nabla)\mathbf{u}+ \mathbf{u} \times \boldsymbol{\omega} = (\nabla \cdot \mathbf{u})\mathbf{u} + \nabla\left(\frac{1}{2}|\mathbf{u}|^2\right)
\]
we have that \eqref{eqn:erik:Euler} is equivalent to
\begin{equation}
\label{eqn:erik:Bernoulli equation}
\mathbf{u} \times \boldsymbol{\omega}=\nabla H
\end{equation}
for divergence-free $\mathbf{u}$, where
\begin{equation}
\label{eqn:erik:Bernoulli function}
H=\tfrac{1}{2}|\mathbf{u}|^2 + P + gz
\end{equation}
is the {\em Bernoulli function}. Here \eqref{eqn:erik:Bernoulli equation} should be interpreted as saying that the left hand side is the gradient of some function. Given $H$, we can then define the pressure using \eqref{eqn:erik:Bernoulli function}. We can also use $H$ to rewrite the dynamic boundary condition \eqref{eqn:erik:dynamic BC} in the form
\[
\tfrac{1}{2}|\mathbf{u}|^2 -H + g\eta-\sigma \nabla \cdot \left(\frac{\nabla \eta}{\sqrt{1+|\nabla \eta|^2}}\right)=Q.
\]
Equation \eqref{eqn:erik:Bernoulli equation} shows that $H$ is constant along both the streamlines (the trajectories of $\mathbf{u}$) and the vortex lines (the trajectories of $\boldsymbol{\omega}$).
This in fact implies a kind of integrable structure under the assumption that $H$ is not constant, as observed by Arnold.
\begin{theorem}[\cite{Arnold65, Arnold66, ArnoldKhesin98}]
Let $\mathbf{u}$ be a solution to the steady Euler equations in a bounded domain $\Omega\subset \mathbb{R}^3$ with  $\mathbf{u}\cdot \mathbf{n}=0$ on $\partial \Omega$ and assume that both $\mathbf{u}$ and $\partial \Omega$ are analytic. If $\mathbf{u}$ and $\boldsymbol{\omega}$ are not everywhere collinear, then there exists an analytic set $C$ of codimension at least one, such that $\Omega \setminus C$ is a union of subdomains which each fall into one of two categories:
\begin{itemize}
\item the subdomain is fibered by invariant tori, on which the flow is conjugate to a linear flow (rational or irrational); or
\item the subdomain is fibered by invariant cylinders  with boundaries on $\partial \Omega$, on which the flow is periodic.
\end{itemize}
The invariant tori and cylinders are the level sets of $H$.
\end{theorem}
The analyticity assumptions can be relaxed with some modifications of the theorem (e.g.~assuming that $\mathbf{u}\times \boldsymbol{\omega}\ne 0$ throughout $\overline{\Omega}$).

Thus, we see that there is a dichotomy between flows with $\mathbf{u}$ and $\boldsymbol{\omega}$ collinear and not. In the latter case, the flow is very structured, whereas in the former case it could in principle be chaotic. Clearly, if $\mathbf{u}$ and $\boldsymbol{\omega}$ are everywhere collinear, then \eqref{eqn:erik:Bernoulli equation} is satisfied with $H$ constant. In particular, this is the case if the flow is irrotational. More generally, we can assume that $\boldsymbol{\omega}=f(\mathbf{x}) \mathbf{u}$ for some function $f(\mathbf{x})$. Such velocity fields are called {\em Beltrami fields} or {\em force-free fields}. Taking the divergence of both sides, and using $\nabla \cdot \mathbf{u}=0$, we see that $f$ has to satisfy $\mathbf{u}\cdot \nabla f=0$, that is, $f$ is constant along the streamlines. A special case is when $f(\mathbf{x})\equiv a$ is everywhere constant. Such fields are called {\em strong Beltrami fields} or {\em linear force-free fields}. Taking the curl of the vorticity and using that $\nabla \times (\nabla\times \mathbf{u}) =\nabla(\nabla \cdot \mathbf{u}) -\Delta \mathbf{u}$, we obtain that
\[
-\Delta \mathbf{u}=a^2 \mathbf{u},
\]
that is, strong Beltrami fields are eigenfunctions of the vector Laplacian. Hence, the equations are elliptic in this case, just like for irrotational flows.
 Beltrami fields with non-constant proportionality factor are much more difficult to construct than strong Beltrami fields (see e.g.~\cite{EncisoPeralta-Salas16}), due to the extra condition on $f$. The intuition that Beltrami fields can be chaotic has recently been verified rigorously; see e.g.~\cite{EncisoPeralta-Salas18, EncisoPeralta-SalasRomaniega20} and references therein. So far there are very few constructions of steady Euler flows with non-constant Bernoulli function without assuming some continuous symmetry which reduces the equations to a two-dimensional problem.

Besides ideal fluid flows, this problem has applications in nuclear fusion theory, where one tries to contain a plasma using external magnetic fields. Note that the equations governing the steady state of a perfectly conducting fluid are
\begin{align*}
\mathbf{J}\times \mathbf{B}&=\nabla p,\\
\nabla \cdot \mathbf{B}&=0,\\
\nabla \times \mathbf{B}&=\mu_0 \mathbf{J},
\end{align*}
where $\mathbf{B}$ is the magnetic field, $\mathbf{J}$ the current density, $p$ the pressure and $\mu_0$ the magnetic constant.
These are precisely the same equations as for steady ideal fluid flows if we identify $\mathbf{B}$ with $\mathbf{u}$, $\mu_0 \mathbf{J}$ with $\boldsymbol{\omega}$ and $\mu_0 p$ with $-H$.  It is beneficial if one can keep the plasma in (or at least close to) equilibrium. The simplest equilibria are axisymmetric, which essentially reduces the problem to 2D. The existence of asymmetric equilibria is a long-standing and notoriously difficult problem and some researchers have in fact argued that there may be no smooth solutions, at least except in very special situations \cite{Grad71}. One of the few rigorous existence result is by  Lortz \cite{Lortz70}, who constructed magnetohydrostatic equilibria in reflection symmetric toroidal domains. We will return to this below.

\subsubsection{Water waves on Beltrami flows}

In this section we summarise the recent paper \cite{LokharuSethWahlen20} on steady doubly-periodic water waves over Beltrami flows. The main assumption is that the velocity field is a (strong) Beltrami field, that is,  $\boldsymbol{\omega}=a \mathbf{u}$ for some real constant $a$. By the discussion above, we can then reformulate the water wave problem as
\begin{subequations}
	\label{eqn:erik:Beltrami problem}
	\begin{align}
	\nabla \times \mathbf{u} &=a \mathbf{u} && \text{in $\Omega$,} 	\\
	\nabla \cdot \mathbf{u}&=0 && \text{in $\Omega$,} \\
	\mathbf{u} \cdot \mathbf{n}&=0 && \text{on $\partial \Omega$,}\\
	\frac12 |\mathbf{u}|^2 + g \eta - \sigma \nabla \cdot \left(\frac{\nabla \eta}{\sqrt{1+|\nabla \eta|^2}}\right)&=Q && \text{on $S$,}
	\end{align}
where $a$ is a constant and $Q$ the Bernoulli constant.  In the case $a=0$, we recover the irrotational problem. We note here that for $a\ne 0$, problem \eqref{eqn:erik:Beltrami problem}  does {\em not} respect the symmetries $ST_1$ and $T_2$. There is however still the combined symmetry $ST_1T_2$, so we can look for weakly symmetric solutions (meaning in particular that $\eta(-\mathbf{x}')=\eta(\mathbf{x}')$, $\mathbf{u}(-\mathbf{x}',z)=(\mathbf{u}'(\mathbf{x}', z), -w(\mathbf{x}', z)$).

Assuming that $\eta\equiv 0$ and that $\mathbf{u}=\mathbf{U}(z)$ has no dependence on the horizontal coordinates, we find a two-parameter family of laminar flows
\[
\mathbf{U}_\mathbf{c}(z)=c_1 \mathbf{U}^{(1)}(z)+c_2 \mathbf{U}^{(2)}(z), \ \ \mathbf{c}=(c_1, c_2)\in \mathbb{R}^2, \text{ where}
\]
\[
\mathbf{U}^{(1)}(z)=(\cos(a z), -\sin(a z),0), \ \
\mathbf{U}^{(2)}(z)=(\sin(a z), \cos(a z),0),
\]
with $Q=Q(\mathbf{c}) \coloneqq \frac12(c_1^2 + c_2^2)$. Physically this means that the velocity field is constant at each vertical position but its direction varies with height (Figure \ref{fig:erik:Beltrami}). Although we are using the notation $\mathbf{c}$, this does not necessarily have to be the wave velocity. In the irrotational case, $\mathbf{U}_\mathbf{c}(z)\equiv (\mathbf{c}, 0)$ and in general $(\mathbf{c}, 0)$ is the relative velocity of the laminar flow at the surface.

\begin{figure}
  \centering
  \includegraphics[width=0.5\linewidth]{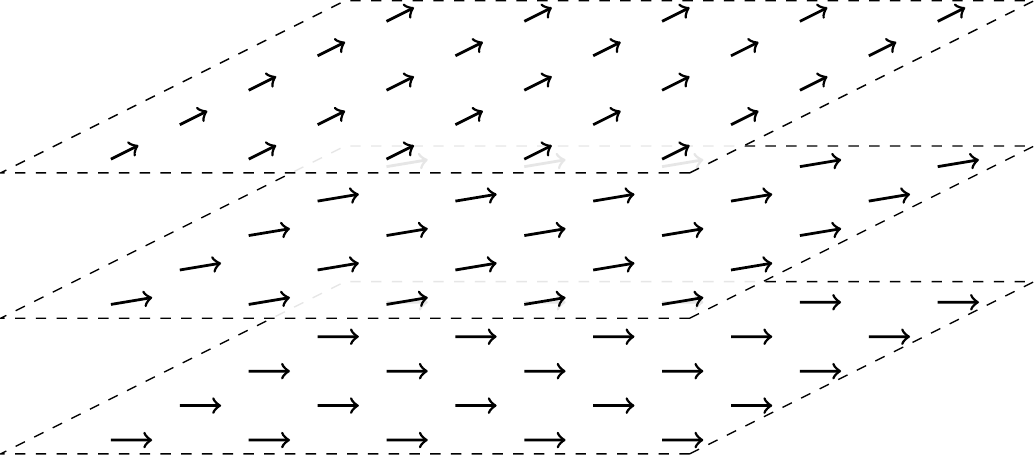}
  \caption{Laminar Beltrami flows from which small-amplitude doubly periodic waves bifurcate.}
  \label{fig:erik:Beltrami}
\end{figure}

The waves constructed in \cite{LokharuSethWahlen20} are  doubly periodic small-amplitude solutions of \eqref{eqn:erik:Beltrami problem}  bifurcating from a laminar flow $\mathbf{U}_{\mathbf{c}}$ for an appropriate value $\mathbf{c}=\mathbf{c}^\star$. The solutions are periodic with respect to a given lattice $\Lambda=\{m_1\boldsymbol{\ell}_1+m_2 \boldsymbol{\ell}_2 \colon m_1, m_2 \in \mathbb{Z}\}\subset \mathbb{R}^2$ and weakly symmetric. Note that the lattice is not required to be symmetric, that is, we allow $|\boldsymbol{\ell}_1|\ne |\boldsymbol{\ell}_2|$, as in Theorem \ref{thm:erik:Craig-Nicholls}.
To single out a two-parameter family of solutions, we require that $Q=Q(\mathbf{c})$ and impose the additional integral conditions
\begin{alignat}{2}
	\int_{\Omega_0}  u_j \,  {{\mathrm{d}\kern0.2pt}\kern0.2pt} \mathbf{x}&=\int_{\Omega_0} U_{\mathbf{c},j}\, {{\mathrm{d}\kern0.2pt}\kern0.2pt} \mathbf{x} &&\quad \text{for } j=1,2,
	\end{alignat}
\end{subequations}
where $\Omega_0$ is as in Section \ref{erik:sec:doubly periodic}. The appropriate value of $\mathbf{c}$ is found by solving the `dispersion equation'
\[
D_a(\mathbf{c}, \mathbf{k})\coloneqq g + \sigma |\mathbf{k}|^2-\frac{(\mathbf{c} \cdot \mathbf{k})^2}{|\mathbf{k}|^2} \kappa(|\mathbf{k}|) + a \frac{(\mathbf{c} \cdot \mathbf{k}) (\mathbf{c} \cdot \mathbf{k}_\perp) }{|\mathbf{k}|^2} = 0, \quad \mathbf{k}\ne 0,
\]
where
\begin{equation*}
\kappa(|\mathbf{k}|)=
\begin{cases}
\sqrt{a^2-|\mathbf{k}|^2} \cot(\sqrt{a^2-|\mathbf{k}|^2}d) & \text{if } |a|>|\mathbf{k}|,\\[0.5em]
\sqrt{|\mathbf{k}|^2-a^2}\coth(\sqrt{|\mathbf{k}|^2-a^2}d) & \text{if } |a|<|\mathbf{k}|,
\end{cases}
\end{equation*}
and $\mathbf{k}_\perp= (-l, k)$,
in which $\mathbf{k}=(k,l) \in \Lambda'$ denotes the wave vector. The bifurcation condition is that the dispersion equation has precisely four solutions in the dual lattice, and for simplicity we assume that these are given by $\pm \mathbf{k}_1$ and $\pm \mathbf{k}_2$ (note that $D_a(\mathbf{c}, -\mathbf{k})=D_a(\mathbf{c}, \mathbf{k})$). This means that the kernel of the linearised problem is two-dimensional when taking the symmetries into account. In contrast to the irrotational case, we cannot reduce to a one-dimensional kernel using the symmetries, even if we assume that the lattice is symmetric. Therefore, we  have to use two bifurcation parameters like in Theorem \ref{thm:erik:Iooss-Plotnikov}.

\begin{theorem}[\cite{LokharuSethWahlen20}]
\label{thm:erik:Beltrami}
Assume that
	\begin{itemize}
		\item[(i)] $\sqrt{a^2 - |\mathbf{k}|^2} \not \in \frac{\pi}{d}\mathbb Z_+$  for all  $\mathbf{k} \in \Lambda'$ such that $|\mathbf{k}|< |a|$;
		\item[(ii)] within the lattice $\Lambda'$, the dispersion relation with $\mathbf{c}=\mathbf{c}^\star$ has exactly four roots $\pm \mathbf{k}_1$ and $\pm \mathbf{k}_2$;
		\item[(iii)] the transversality condition $\nabla_{\mathbf{c}} D_a(\mathbf{c}^\star, \mathbf{k}_1) \nparallel \nabla_{\mathbf{c}} D_a(\mathbf{c}^\star, \mathbf{k}_2)$  holds.
	\end{itemize}
	Then there exists a two-parameter family of weakly symmetric doubly periodic solutions $(\eta, \mathbf{u}, \mathbf{c})(\boldsymbol{\varepsilon})\in  C^{2+\alpha}_\text{per}(\mathbb{R}^2)  \times C^{1+\alpha}_\text{per}(\overline{\Omega};\mathbb{R}^3)\times \mathbb{R}^2$, indexed by $\boldsymbol{\varepsilon}=(\varepsilon_1, \varepsilon_2)$, near $\mathbf{U}_{\mathbf{c}^\star}$ with
	\[
		\eta(\mathbf{x}')=\varepsilon_1 \cos(\mathbf{k}\cdot \mathbf{x}')+\varepsilon_2\cos(\mathbf{k}_2\cdot \mathbf{x}') +O(|\boldsymbol{\varepsilon}|^2) \text{ and } \mathbf{c}=\mathbf{c}^\star+O(|\boldsymbol{\varepsilon}|^2) \text{ as } \boldsymbol{\varepsilon}\to 0.
	\]
\end{theorem}

The proof involves reducing the problem to a single scalar equation for $\eta$ and using a multi-parameter local bifurcation argument based on Lyapunov-Schmidt reduction. A more detailed discussion of when the assumptions in the theorem are satisfied can be found in \cite{LokharuSethWahlen20}. Note in particular that it is allowed to take $a=0$, in which case we obtain an alternative existence proof of asymmetric irrotational solutions. In that case, assumptions (i) and (iii) are automatically satisfied. Since we can translate these solutions in $x$ and $y$, the result is consistent with Theorem \ref{thm:erik:Craig-Nicholls}.

Let us also mention that in addition to this existence result, there are two recent variational formulations for water waves over Beltrami flows. The first, by Lokharu \& Wahl\'en \cite{LokharuWahlen19} allows for overhanging waves, while the second, by Groves \& Horn \cite{GrovesHorn20}, includes a reduction to the surface which coincides with the steady Zakharov-Craig-Sulem formulation in the irrotational case. These  may be useful in future variational existence theories for doubly periodic waves and fully localised solitary waves. Open questions include the existence of large-amplitude doubly periodic waves and existence of any kind of fully localised solitary waves.

\subsubsection{Symmetric waves with small vorticity}

Finally, we report on a recent existence result for symmetric diamond waves with small vorticity, inspired by a construction of  magnetohydrostatic equilibria in reflection symmetric toroidal domains by Lortz \cite{Lortz70}. Lortz paper contains one of the few rigorous constructions of non-axisymmetric equilibria.  The main challenge in using his approach to construct water waves is how to incorporate the free boundary.

We will consider the same kind of lattice as in Reeder's and Shinbrot's work, that is, $\Lambda$ is generated by
$\boldsymbol{\ell}_1$ and $\boldsymbol{\ell}_2$, while $\Lambda'$ is generated by $\boldsymbol{k}_1$ and $\boldsymbol{k}_2$, where
\begin{align*}
\boldsymbol{k}_1&=(\kappa_1, \kappa_2), & \boldsymbol{k}_2&=(\kappa_1, -\kappa_2),\\
\boldsymbol{\ell}_1&=\pi(\kappa_1^{-1}, \kappa_2^{-1}), &  \boldsymbol{\ell}_2&=\pi(\kappa_1^{-1}, -\kappa_2^{-1}).
\end{align*}
For convenience, we set $\lambda_1=2\pi/\kappa_1$ and $\lambda_2=2\pi/\kappa_2$, so that $\boldsymbol{\ell}_1=(\lambda_1/2, \lambda_2/2)$. In addition, we will use the rectangular sublattice $\Gamma\coloneqq \lambda_1 \mathbb{Z}\times \lambda_2 \mathbb{Z}\subset \Lambda$ generated by $\boldsymbol{\ell}_1+ \boldsymbol{\ell}_2$ and $\boldsymbol{\ell}_1- \boldsymbol{\ell}_2$, and its dual lattice $\Gamma'=\kappa_1 \mathbb{Z}\times \kappa_2 \mathbb{Z}$ (note that $\Lambda' \subset \Gamma'$). In fact, we will first construct solutions which are periodic with respect to $\Gamma$ and then verify a posteriori that they are periodic with respect to $\Lambda$.
 We will look for solutions with the symmetries $ST_1$ and $T_2$ as in Theorem \ref{thm:erik:Reeder-Shinbrot}. In fact, this will be essential for the very formulation of the problem since it guarantees that the streamlines are periodic with respect to $\Gamma$.

 Recall that a  divergence-free vector field can be written locally as the cross product of two gradients \cite{Barbarosie11}, which in particular applies to the vorticity $\boldsymbol{\omega}$. The starting point for Lortz' method is the assumption that this can be done globally, with one of the functions being $H$. We thus write
\[
\boldsymbol{\omega}=\nabla H \times \nabla \tau.
\]
By standard vector calculus identities, we then get
\[
\mathbf{u}\times \boldsymbol{\omega}=\mathbf{u}\times (\nabla H \times \nabla \tau)=(\mathbf{u}\cdot \nabla \tau)\nabla H - (\mathbf{u}\cdot \nabla H)\nabla \tau.
\]
Since the Bernoulli function is supposed to be constant along the streamlines, we should therefore take $\mathbf{u}\cdot \nabla \tau=1$. If we set $\tau|_{x=0}=0$, $\tau(\mathbf{x})$ can be interpreted as the time it takes to go from the plane $\{x=0\}$ to the point $\mathbf{x}$ along a streamline (assuming that  $u\ne 0$  throughout $\overline{\Omega}$). It is easily seen that $\tau$ is $\lambda_2$-periodic in $y$.
Next, we introduce the function $q(\mathbf{x})=\tau(\mathbf{x}+\lambda_1 \mathbf{e}_1)-\tau(\mathbf{x})$ and note that $\mathbf{u}\cdot \nabla q=0$ by the periodicity of $\mathbf{u}$. Therefore, $q$ is constant on the streamlines. Forgetting about the definition of $H$ we can make the Ansatz $H=h(q)$ for some function $h$, in which case $H$ will automatically be constant along the streamlines and we obtain
\[
\mathbf{u}\times \boldsymbol{\omega}=\underbrace{(\mathbf{u}\cdot \nabla \tau)}_{=1}\nabla H - \underbrace{(\mathbf{u}\cdot \nabla H)}_{=0}\nabla \tau=\nabla H.
\]
The Euler equations are then satisfied if we define $P\coloneqq Q+h(q)-\frac12 |\mathbf{u}|^2-gz$ for some arbitrary constant $Q$.

Based on the above considerations, we look for solutions of the problem
\begin{subequations}
\label{eqn:erik:svw problem}
\begin{alignat}{2}
\label{eqn:erik:svw I}
\nabla \times \mathbf{u}  &= h'(q) \nabla q \times \nabla \tau  &  &\quad \text{in $\Omega$}, \\
\nabla \cdot \mathbf{u}   &= 0   &  &\quad \text{in $\Omega$}, \\
\mathbf{u} \cdot\mathbf{n} & = 0 &  &\quad \text{on $\partial \Omega$}, \\
\frac{1}{2}|\mathbf{u}|^2 -h(q)+ g\eta - \sigma \nabla \cdot \left(\frac{\nabla \eta}{\sqrt{1+|\nabla \eta|^2}}\right)&= Q &  &\quad \text{on $S$},
\end{alignat}
satisfying the integral condition
\begin{equation}
 \frac{1}{|\Omega_0|} \int_{\Omega_0} u \, {{\mathrm{d}\kern0.2pt}\kern0.2pt} \mathbf{x} = c.
\end{equation}
\end{subequations}
The number $c$ can be thought of as a wave speed and will be used as a bifurcation parameter; in particular, we consider it part of the solutions rather than a fixed constant.
The symmetry assumptions are crucial since they mean that both $q$ and $\nabla q\times \nabla \tau$ are periodic with respect to $\Gamma$ (and therefore also the right hand side of  \eqref{eqn:erik:svw I}). Indeed, they imply that the streamlines are periodic, so that  $\mathbf{x}$ and $\mathbf{x}+\lambda_1\mathbf{e}_1$ are on the same streamline and hence $q(\mathbf{x}+\lambda_1 \mathbf{e}_1)=q(\mathbf{x})$ (periodicity in the $y$-direction follows from the same property for $\tau$). From this, the periodicity of $\nabla q\times \nabla \tau$ follows in a straightforward way.
Observe that for any real number $c$, we obtain a trivial solution by setting
\[
\eta=0, \quad \mathbf{u}=c\mathbf{e}_1,\quad \tau = \frac{x}{c}, \quad q=\frac{\lambda_1}{c}
\]
if we also pick
\begin{equation*}
Q = Q(c) \coloneqq \frac{1}{2}c^2- h\left(\frac{\lambda_1}{c}\right).
\end{equation*}
Furthermore, it turns out that if we formally linearise the problem, we obtain the same linear dispersion relation as for the irrotational water wave problem, namely
\[
D_0(c\mathbf{e}_1, \mathbf{k})\coloneqq g+\sigma|\mathbf{k}|^2-\frac{c^2 k^2}{|\mathbf{k}|}\coth(|\mathbf{k}|d)=0,
\]
where $\mathbf{k}=(k, l)$.
In order to find doubly periodic waves, we pick a $c$ such that the dispersion relation is satisfied for the generators $\mathbf{k}_1=(\kappa_1, \kappa_2)$ and $\mathbf{k}_2=(\kappa_1, -\kappa_2)$ of the lattice $\Lambda'$.  In the theorem below, `per' refers to periodicity with respect to $\Lambda$ rather than $\Gamma$.

\begin{theorem}[\cite{SethVarholmWahlen21}]
Assume that within the lattice $\Gamma'$, the dispersion relation  with $c=c^\star$ has exactly the four roots $\pm \mathbf{k}_1$ and $\pm \mathbf{k}_2$. Then  for every $h\in C^{4+\alpha}(\mathbb{R})$  with $\| h \|_{4+\alpha}\ll 1$, there exists a family of symmetric   solutions $(\eta,\mathbf{u},c)(\varepsilon)\in C_\text{per}^{5+\alpha}(\mathbb{R}^2)\times C_\text{per}^{4+\alpha}(\overline{\Omega}; \mathbb{R}^3)\times \mathbb{R}$, $|\varepsilon|<\varepsilon_0$, to  problem  \eqref{eqn:erik:svw problem} bifurcating from the trivial solution $(0, c^\star\mathbf{e}_1, c^\star)$ at $\varepsilon=0$. Moreover,
\[
\eta=\varepsilon\cos(\kappa_1 x)\cos(\kappa_2 y)+o(\varepsilon) \text{ as } \varepsilon \to 0
\]
in $C_\text{per}^{4+\alpha}(\mathbb{R}^2)$.
If in addition $h'(\lambda_1/c^\star)\ne 0$ and $h'(\lambda_1/c^\star)\lambda_1+(c^\star)^3 \ne 0$, then the solutions have nonzero vorticity for $0<|\varepsilon| \ll 1$.
\end{theorem}

Note that the condition on the dispersion relation is essentially the same as in Theorem \ref{thm:erik:Reeder-Shinbrot}, except that we need to rule out solutions  in the larger lattice $\Gamma'$. Let us give some comments on the proof. We begin by looking for $\Gamma$-periodic, symmetric solutions. Just like in the proof of Theorem \ref{thm:erik:Beltrami}, the idea  is to first reduce to the boundary and then use a Crandall-Rabinowitz type argument. Since after restricting to symmetric solutions the kernel is one-dimensional, here one could actually use the standard Crandall-Rabinowitz theorem, were it not for the problem that the nonlinear operator is not Fr\'echet differentiable. This issue can be seen if we think of the solution of the problem
\[
\mathbf{u}\cdot \nabla \tau =1, \quad \tau|_{x_1}=0
\]
as a function $\tau(\mathbf{u})$ of $\mathbf{u}$ and try to differentiate with respect to $\mathbf{u}$. First note that if $\mathbf{u}$ belongs to $C^{m+\alpha}$, then so does $\tau(\mathbf{u})$. But unlike elliptic equations $\tau(\mathbf{u})$ does not gain any derivatives (except along the streamlines).
The formal derivative $\mathrm{d}_\mathbf{u} \tau(\mathbf{u})\mathbf{v}$ in the direction $\mathbf{v}$ is given by
\[
\mathbf{u} \cdot \nabla (\mathrm{d}_{\mathbf{u}} \tau(\mathbf{u})\mathbf{v}) =-\mathbf{v}\cdot \nabla \tau[\mathbf{u}], \quad  \mathrm{d}_\mathbf{u} \tau(\mathbf{u})\mathbf{v}|_{x_1}=0,
\]
and we note that the right hand side of the transport equation belongs to $C^{m-1+\alpha}$ if $\mathbf{u}, \mathbf{v}\in C^{m+\alpha}$. Hence, we don't expect $\tau(\mathbf{u})$ to be Fr\'echet differentiable as a map from $C^{m+\alpha}$ to $C^{m+\alpha}$, but only if we allow some loss of regularity when we take a derivative. By carefully keeping track of this loss of regularity in all the steps, we are able to modify the proof of the Crandall-Rabinowitz theorem so that it applies to our problem. Note that the loss of regularity is visible in the theorem, where the family $\eta(\varepsilon)$ lies in $C^{5+\alpha}$ but the asymptotic formula holds in $C^{4+\alpha}$. The  periodicity with respect to the lattice $\Lambda$ can be seen by using the fact that $(\eta(\varepsilon)(\boldsymbol{x}' +\boldsymbol{\ell}_j), \mathbf{u}(\varepsilon)(\boldsymbol{x}'+\boldsymbol{\ell}_j, z))$ is also a solution with the same projection onto the kernel. Finally, an expansion of  $\mathbf{u}(\varepsilon)$ in powers of $\varepsilon$ reveals that the vorticity is non-zero.

The existence of large-amplitude waves is of course an important open question here as well, as is the existence of pure gravity waves. But perhaps more fundamental is the possibility of finding solutions which bifurcate from other types of shear flows than uniform flows. It turns out that for strictly monotone shear flows, one is not free to choose the function $h$. It has to be of the form $h(q)=\frac{\lambda_1^2}{2q^2}$ and the kernel of the linearised problem at the shear flow is infinite-dimensional (due to the existence of an infinite-dimensional family of other nearby shear flows). Thus some new idea is needed. New ideas are also required to handle asymmetric solutions.

\section*{Acknowledgments}
The authors wish to thank the organizers of the Online Northeast PDE Seminar (J.~G\'omez-Serrano, B.~Pausader, F.~Pusateri, and I.~Tice) for initiating and running the lecture series that led to this article. We are especially grateful to B.~Pausader for his help in editing the manuscript.

The work of SVH was partially funded by the Austrian Science Fund (FWF), Grant Z 387-N.  The work of VMH was partially funded by the NSF through the award DMS-2009981. The work of SW was partially funded by the NSF through the award DMS-1812436.  The work EW was partially funded by the European Research Council (ERC) under the European Union's Horizon 2020 research and innovation programme (grant agreement no 678698) and the Swedish Research Council (grant nos.~621-2012-3753 and 2016-04999).

\bibliographystyle{amsplain}

\begin{thebibliography}{100}

\bibitem{AblowitzSegur79}
Mark~J. Ablowitz and Harvey Segur, \emph{On the evolution of packets of water
  waves}, J. Fluid Mech. \textbf{92} (1979), no.~4, 691--715.

\bibitem{AkersReeger17}
Benjamin~F. Akers and Jonah~A. Reeger, \emph{Three-dimensional overturned
  traveling water waves}, Wave Motion \textbf{68} (2017), 210 -- 217.

\bibitem{ambrose2016global}
David~M. Ambrose, Walter~A. Strauss, and J.~Douglas Wright, \emph{Global
  bifurcation theory for periodic traveling interfacial gravity-capillary
  waves}, Ann. Inst. H. Poincar\'{e} Anal. Non Lin\'{e}aire \textbf{33} (2016),
  no.~4, 1081--1101. \MR{3519533}

    \bibitem{amick1984semilinear}
C.~J. Amick, \emph{Semilinear elliptic eigenvalue problems on an infinite
  strip with an application to stratified fluids}, Ann. Scuola Norm. Sup. Pisa
  Cl. Sci. (4) \textbf{11} (1984), no.~3, 441--499. \MR{MR785621 (86i:35042)}

\bibitem{amick}
\bysame, Bounds for water waves.
\emph{Arch. Rational Mech. Anal.} \textbf{99}  (2) (1987),  91--114.

\bibitem{AFT1982}
C.~J. Amick, L.~E. Fraenkel, and J.~F. Toland, \emph{On the {S}tokes conjecture
  for the wave of extreme form}, Acta Math. \textbf{148} (1982), 193--214.
  \MR{666110}

 \bibitem{at:periodic}
C.~J. Amick and J.~F. Toland, \emph{On periodic water-waves and their
  convergence to solitary waves in the long-wave limit}, Philos. Trans. Roy.
  Soc. London Ser. A \textbf{303} (1981), no.~1481, 633--669. \MR{647410
  (83b:76009)}

\bibitem{at:finite}
\bysame, \emph{On solitary water-waves of finite amplitude}, Arch. Rational
  Mech. Anal. \textbf{76} (1981), no.~1, 9--95. \MR{629699 (83b:76017)}

\bibitem{amick1984limiting}
\bysame, \emph{The limiting form of internal waves}, Proc.
  Roy. Soc. London Ser. A \textbf{394} (1984), no.~1807, 329--344. \MR{763506}

\bibitem{at:homoclinic}
\bysame, \emph{Solitary waves with surface tension. {I}. {T}rajectories
  homoclinic to periodic orbits in four dimensions}, Arch. Rational Mech. Anal.
  \textbf{118} (1992), no.~1, 37--69. \MR{1151926}

\bibitem{amick1986global}
C.~J. Amick and R.~E.~L. Turner, \emph{A global theory of internal solitary
  waves in two-fluid systems}, Trans. Amer. Math. Soc. \textbf{298} (1986),
  no.~2, 431--484. \MR{MR860375 (87m:35210)}

\bibitem{amick1989small}
\bysame, \emph{Small internal waves in two-fluid systems}, Arch. Rational Mech.
  Anal. \textbf{108} (1989), no.~2, 111--139. \MR{1011554 (90h:76037)}




\bibitem{Arnold65}
Vladimir~I. Arnold, \emph{Sur la topologie des \'ecoulements stationnaires des
  fluides parfaits}, C. R. Acad. Sci. Paris \textbf{261} (1965), 17--20.
  \MR{0180949 (31 \#5179)}

  \bibitem{Arnold66}
\bysame, \emph{Sur la g\'eom\'etrie diff\'erentielle des groupes de {L}ie de
  dimension infinie et ses applications \`a l'hydrodynamique des fluides
  parfaits}, Ann. Inst. Fourier (Grenoble) \textbf{16} (1966), no.~fasc. 1,
  319--361. \MR{0202082 (34 \#1956)}

\bibitem{ArnoldKhesin98}
Vladimir~I. Arnold and Boris~A. Khesin, \emph{Topological methods in
  hydrodynamics}, Applied Mathematical Sciences, vol. 125, Springer-Verlag, New
  York, 1998. \MR{1612569 (99b:58002)}

   \bibitem{aabab}
A. I. Aptekareva \& N. G. Afendikovaa.
About the works of K. I. Babenko in the ield
of mechanics and applied mathematics
(on the 100th Anniversary of His Birth)
\emph{Mechanics of Solids}, \textbf{55} (7) (2020), 919–925.
Russian text: \emph{Prikladnaya Matematika i Mekhanika}, \textbf{84} (1) (2020),  3--12.

\bibitem{Babenko}
K.~I. Babenko, \emph{Some remarks on the theory of surface waves of finite
  amplitude}, Dokl. Akad. Nauk SSSR \textbf{294} (1987), no.~5, 1033--1037.
  \MR{898306}

\bibitem{BagriGroves14}
G.~S. Bagri and M.~D. Groves, \emph{A spatial dynamics theory for doubly periodic
  travelling gravity-capillary surface waves on water of infinite depth},
  Journal of Dynamics and Differential Equations (2014), 1--28 (English).

\bibitem{Balk}
 A. M. Balk. \emph{A Lagrangian for water waves}, Phys. Fluids \textbf{8(2)} (1996), 416--420.


\bibitem{Barbarosie11}
Cristian Barbarosie, \emph{Representation of divergence-free vector fields},
  Quart. Appl. Math. \textbf{69} (2011), no.~2, 309--316. \MR{2814529}

\bibitem{beale:existence}
J.~Thomas Beale, \emph{The existence of solitary water waves}, Comm. Pure Appl.
  Math. \textbf{30} (1977), no.~4, 373--389. \MR{0445136 (56 \#3480)}

\bibitem{benjamin:conjugate}
T.~Brooke Benjamin, \emph{A unified theory of conjugate flows}, Philos. Trans. Roy.
  Soc. London Ser. A \textbf{269} (1971), 587--643. \MR{446075}

\bibitem{benjamin1968gravity}
\bysame, \emph{Gravity currents and related phenomena}, J. Fluid
  Mech. \textbf{31} (1968), no.~2, 209--248.

\bibitem{benjamin:impulse}
\bysame, \emph{Impulse, flow force and variational principles}, IMA J. Appl.
  Math. \textbf{32} (1984), no.~1-3, 3--68. \MR{740456 (85h:70012a)}

\bibitem{BenneyRoskes69}
D.~J. Benney and G.~Roskes, \emph{Wave instabilities}, Stud. Appl. Math.
  \textbf{48} (1969), 377--385.

\bibitem{bona1983finite}
J.~L. Bona, D.~K. Bose, and R.~E.~L. Turner, \emph{Finite-amplitude steady
  waves in stratified fluids}, J. Math. Pures Appl. (9) \textbf{62} (1983),
  no.~4, 389--439. \MR{MR735931 (85e:76056)}

  \bibitem{bouss}
M. J. Boussinesq \emph{Essai sur la th\'eorie des eaux courantes}
M\'emoirs Acad. Sci. Inst. France (s\'erie 2)23 (1877), 1-68.

\bibitem{Bridges}
T. J. Bridges and A. Mielke, \emph{ A proof of the Benjamin-Feir instability}, Arch. Rational Mech. Anal. 133 (1995), 145--198.

\bibitem{Buffoni04a}
B.~Buffoni, \emph{Existence and conditional energetic stability of
  capillary-gravity solitary water waves by minimisation}, Arch. Ration. Mech.
  Anal. \textbf{173} (2004), no.~1, 25--68. \MR{2073504 (2005f:35008)}

 \bibitem{BuDaTo:98c}
B.~Buffoni, E.~N. Dancer, and J.~F. Toland, \emph{Sur les ondes de Stokes et une conjecture de Levi-Civita}, C. R. Acad. Sci. Paris S\'erie I \textbf{326} (1998), 1265--1268.

\bibitem{BDT2000b}
\bysame, \emph{The regularity and local
  bifurcation of steady periodic water waves}, Arch. Ration. Mech. Anal.
  \textbf{152} (2000), no.~3, 207--240. \MR{1764945}

\bibitem{BDT2000a}
\bysame, \emph{The sub-harmonic bifurcation of {S}tokes waves}, Arch. Ration.
  Mech. Anal. \textbf{152} (2000), no.~3, 241--271. \MR{1764946}

\bibitem{BuffoniGroves99}
B.~Buffoni and M.~D. Groves, \emph{A multiplicity result for solitary
  gravity-capillary waves in deep water via critical-point theory}, Arch.\
  Ration.\ Mech.\ Anal.\/ \textbf{146} (1999), 183--220.

\bibitem{BuffoniGrovesSunWahlen13}
B.~Buffoni, M.~D. Groves, S.~M. Sun, and E.~Wahl{\'e}n, \emph{Existence and
  conditional energetic stability of three-dimensional fully localised solitary
  gravity-capillary water waves}, J. Differential Equations \textbf{254}
  (2013), no.~3, 1006--1096. \MR{2997362}

\bibitem{BuffoniGrovesWahlen18}
Boris Buffoni, Mark~D. Groves, and Erik Wahl{\'e}n, \emph{A variational
  reduction and the existence of a fully localised solitary wave for the
  three-dimensional water-wave problem with weak surface tension}, Arch.
  Rational Mech. Anal. \textbf{228} (2018), no.~3, 773--820.

  \bibitem{BuffoniGrovesWahlen21}
\bysame, \emph{Fully localised
  three-dimensional gravity-capillary solitary waves on water of
  infinite-depth}, arXiv preprint
  arXiv:2108.05973 (2021)

   \bibitem{BSTa}
B. Buffoni, \'E.  S\'er\'e \&  J. F. Toland. \emph{Surface water
waves as saddle points of the energy}. Calc. Var. \textbf{17}
(2003), 199-220.

\bibitem{BSTb}
\bysame, \emph{Minimization
methods for quasi-linear problems with an application to periodic
water waves}. SIAM J. Math. Anal. {\bf 36} (4) (2005), 1080-1094.

\bibitem{BuTo:00}
B. Buffoni \& J. F. Toland. \emph{Dual free boundaries for
Stokes waves.} Comptes Rendus Acad. Sci. Paris, S\'erie I \textbf{332} (2000), 73-78.

\bibitem{bt:analytic}
\bysame, \emph{Analytic theory of global bifurcation},
  Princeton Series in Applied Mathematics, Princeton University Press,
  Princeton, NJ, 2003, An introduction. \MR{1956130 (2004b:47117)}

\bibitem{buhler2016wind}
Oliver B{\"u}hler, Jalal Shatah, Samuel Walsh, and Chongchun Zeng, \emph{On the
  wind generation of water waves}, Arch. Ration. Mech. Anal. \textbf{222}
  (2016), no.~2, 827--878. \MR{3544318}

\bibitem{burckel:book}
Robert~B. Burckel, \emph{An introduction to classical complex analysis. {V}ol.
  1}, Pure and Applied Mathematics, vol.~82, Academic Press, Inc. [Harcourt
  Brace Jovanovich, Publishers], New York-London, 1979. \MR{555733}

\bibitem{caffarelli2002monotonicity}
Luis~A. Caffarelli, David Jerison, and Carlos~E. Kenig, \emph{Some new
  monotonicity theorems with applications to free boundary problems}, Ann. of
  Math. (2) \textbf{155} (2002), no.~2, 369--404. \MR{1906591}

\bibitem{cao1996existence}
Daomin Cao, Norman~E. Dancer, Ezzat~S. Noussair, and Shunsen Yan, \emph{On the
  existence and profile of multi-peaked solutions to singularly perturbed
  semilinear {D}irichlet problems}, Discrete Contin. Dynam. Systems \textbf{2}
  (1996), no.~2, 221--236. \MR{1382508}

\bibitem{champneys:persistence}
A.~R. Champneys, \emph{Codimension-one persistence beyond all orders of
  homoclinic orbits to singular saddle centres in reversible systems},
  Nonlinearity \textbf{14} (2001), no.~1, 87--112. \MR{1808625}

\bibitem{CG1993}
G.~A. Chandler and I.~G. Graham, \emph{The computation of water waves modelled
  by {N}ekrasov's equation}, SIAM J. Numer. Anal. \textbf{30} (1993), no.~4,
  1041--1065. \MR{1231326}

\bibitem{ChangGustafsonNakanishiTsai07}
Shu-Ming Chang, Stephen Gustafson, Kenji Nakanishi, and Tai-Peng Tsai,
  \emph{Spectra of linearized operators for {NLS} solitary waves}, SIAM J.
  Math. Anal. \textbf{39} (2007/08), no.~4, 1070--1111. \MR{2368894
  (2008k:35400)}

 \bibitem{chensaff}
B. Chen \&  P. G. Saffman.
\emph{Numerical evidence for the existence of new types of gravity waves of permanent form on deep water.}
Stud. Appl. Math. \textbf{62} (1) (1980), 1--21.

\bibitem{chen2016continuous}
Robin~Ming Chen and Samuel Walsh, \emph{Continuous dependence on the density
  for stratified steady water waves}, Arch. Rational Mech. Anal. \textbf{219}
  (2016), no.~2, 741--792. \MR{3437862}

\bibitem{chen2021orbital}
\bysame, \emph{Orbital stability of internal waves}, arXiv preprint
  arXiv:2102.13590 (2021).

\bibitem{cww:strat}
Robin~Ming Chen, Samuel Walsh, and Miles~H. Wheeler, \emph{Existence and
  qualitative theory for stratified solitary water waves}, Ann. Inst. H.
  Poincar\'e Anal. Non Lin\'eaire \textbf{35} (2018), no.~2, 517--576.
  \MR{3765551}

\bibitem{chen2019center}
\bysame, \emph{Center manifolds without a phase space for quasilinear problems
  in elasticity, biology, and hydrodynamics}, arXiv preprint arXiv:1907.04370
  (2019).

\bibitem{chen2019existence}
\bysame, \emph{Existence, nonexistence, and asymptotics of deep water solitary
  waves with localized vorticity}, Arch. Ration. Mech. Anal. \textbf{234}
  (2019), no.~2, 595--633. \MR{3995048}

\bibitem{chen2020global}
\bysame, \emph{Global bifurcation for monotone fronts of elliptic equations},
  arXiv preprint arXiv:2005.00651 (2020).

\bibitem{chen2020large}
\bysame, \emph{Large-amplitude internal fronts in two-fluid systems}, Comptes
  Rendus. Math\'ematique \textbf{358} (2020), no.~9-10, 1073--1083 (en).

\bibitem{Cipolatti92}
Rolci Cipolatti, \emph{On the existence of standing waves for a
  {D}avey-{S}tewartson system}, Comm. Partial Differential Equations
  \textbf{17} (1992), no.~5-6, 967--988. \MR{1177301 (93h:35188)}

\bibitem{Cipolatti93}
\bysame, \emph{On the instability of ground states for a {D}avey-{S}tewartson
  system}, Ann. Inst. H. Poincar\'{e} Phys. Th\'{e}or. \textbf{58} (1993),
  no.~1, 85--104. \MR{1208793}

\bibitem{const-escher}
Adrian Constantin and Joachim Escher, \emph{Analyticity of periodic traveling free surface water waves with vorticity}, Ann. of Math. (2) 173 (2011), no. 1, 559--568.

\bibitem{CS2004}
Adrian Constantin and Walter Strauss, \emph{Exact steady periodic water
  waves with vorticity}, Comm. Pure Appl. Math. \textbf{57} (2004), no.~4,
  481--527. \MR{2027299 (2004i:76018)}

\bibitem{CS2010}
\bysame, \emph{Pressure beneath a Stokes wave},  Comm. Pure Appl. Math. 63 (2010), 533--557.

\bibitem{CS2011}
\bysame, \emph{Periodic traveling gravity water waves with discontinuous
  vorticity}, Arch. Ration. Mech. Anal. \textbf{202} (2011), no.~1, 133--175.
  \MR{2835865}

\bibitem{csv:critical}
Adrian Constantin, Walter Strauss and Eugen V\u{a}rv\u{a}ruc\u{a},
  \emph{Global bifurcation of steady gravity water waves with critical layers},
  Acta Math. \textbf{217} (2016), no.~2, 195--262. \MR{3689941}


\bibitem{cv:constant}
Adrian Constantin and Eugen V\u{a}rv\u{a}ruc\u{a}, \emph{Steady periodic water waves with
  constant vorticity: regularity and local bifurcation}, Arch. Ration. Mech.
  Anal. \textbf{199} (2011), no.~1, 33--67. \MR{2754336 (2012e:76017)}

\bibitem{cordoba2021existence}
Diego C{\'o}rdoba and Elena {Di Iorio}, \emph{Existence of gravity-capillary
  crapper waves with concentrated vorticity}, arXiv preprint arXiv:2106.15923
  (2021).

\bibitem{craig:nonexistence}
Walter Craig, \emph{Non-existence of solitary water waves in three dimensions},
  R. Soc. Lond. Philos. Trans. Ser. A Math. Phys. Eng. Sci. \textbf{360}
  (2002), no.~1799, 2127--2135, Recent developments in the mathematical theory
  of water waves (Oberwolfach, 2001). \MR{1949966 (2003m:76011)}

\bibitem{CGS}
W. Craig, P. Guyenne and C. Sulem, \emph{Water waves over a random bottom},  J. Fluid Mech. 640 (2009), 79--107.

\bibitem{CraigNicholls00}
Walter Craig and David~P. Nicholls, \emph{Travelling two and three dimensional
  capillary gravity water waves}, SIAM J. Math. Anal. \textbf{32} (2000),
  no.~2, 323--359. \MR{1781220 (2002c:76016)}

\bibitem{CraigNicholls02}
\bysame, \emph{Traveling gravity water waves in two and three dimensions}, Eur.
  J. Mech. B Fluids \textbf{21} (2002), no.~6, 615--641. \MR{1947187}

\bibitem{craigsulem}
W. Craig \& C. Sulem.
\emph{Numerical simulation of gravity waves.}
J. Comput. Phys. \textbf{108} (1) (1993), 73--83.

\bibitem{cs:sym}
Walter Craig and Peter Sternberg, \emph{Symmetry of solitary waves}, Comm.
  Partial Differential Equations \textbf{13} (1988), no.~5, 603--633.
  \MR{919444 (88m:35132)}

\bibitem{craik}
Alex D.~D. Craik, \emph{The origins of water wave theory}, Annual review of
  fluid mechanics. {V}ol. 36, Annu. Rev. Fluid Mech., vol.~36, Annual Reviews,
  Palo Alto, CA, 2004, pp.~1--28. \MR{2062306 (2005a:01012)}

  \bibitem{craik2}
A. Craik.
George Gabriel Stokes on water wave theory. \emph{Annual Review of Fluid Mechanics.} \textbf{37} (2005), 23-42.

\bibitem{rabinowitz:simple}
Michael~G. Crandall and Paul~H. Rabinowitz, \emph{Bifurcation from simple
  eigenvalues}, J. Functional Analysis \textbf{8} (1971), 321--340. \MR{0288640
  (44 \#5836)}

\bibitem{Crapper}
G.~D. Crapper, \emph{An exact solution for progressive capillary waves of
  arbitrary amplitude}, J. Fluid Mech. \textbf{2} (1957), 532--540. \MR{91075}

\bibitem{dancer}
E.~N. Dancer, \emph{Bifurcation theory for analytic operators}, Proc. London
  Math. Soc. (3) \textbf{26} (1973), 359--384. \MR{322615}

\bibitem{dancer:global}
\bysame, \emph{Global structure of the solutions of non-linear real analytic
  eigenvalue problems}, Proc. London Math. Soc. (3) \textbf{27} (1973),
  747--765. \MR{0375019}

\bibitem{darrigol:horse}
Olivier Darrigol, \emph{The spirited horse, the engineer, and the
  mathematician: water waves in nineteenth-century hydrodynamics}, Arch. Hist.
  Exact Sci. \textbf{58} (2003), no.~1, 21--95. \MR{2020055 (2004k:76002)}

\bibitem{DaveyStewartson74}
A.~Davey and K.~Stewartson, \emph{On three-dimensional packets of surface
  waves}, Proc. Roy. Soc. London Ser. A \textbf{338} (1974), 101--110.
  \MR{MR0349126 (50 \#1620)}

\bibitem{DeBouardSaut96}
Anne de~Bouard and Jean-Claude Saut, \emph{Remarks on the stability of generalized {KP}
  solitary waves}, Mathematical problems in the theory of water waves
  ({L}uminy, 1995), Contemp. Math., vol. 200, Amer. Math. Soc., Providence, RI,
  1996, pp.~75--84. \MR{1410501}

\bibitem{DeBouardSaut97}
\bysame, \emph{Solitary waves of generalized
  {K}adomtsev-{P}etviashvili equations}, Ann. Inst. H. Poincar\'{e} Anal. Non
  Lin\'{e}aire \textbf{14} (1997), no.~2, 211--236. \MR{1441393}

\bibitem{DengSun09}
Shengfu Deng and Shu-Ming Sun, \emph{Three-dimensional gravity-capillary waves
  on water---small surface tension case}, Phys. D \textbf{238} (2009), no.~17,
  1735--1751. \MR{2574421}

\bibitem{DengSun10}
\bysame, \emph{Exact theory of three-dimensional water waves at the critical
  speed}, SIAM J. Math. Anal. \textbf{42} (2010), no.~6, 2721--2761.
  \MR{2745790}

\bibitem{DengIonescuPausaderPusateri17}
Yu~Deng, Alexandru~D. Ionescu, Beno\^{\i}t Pausader, and Fabio Pusateri,
  \emph{Global solutions of the gravity-capillary water-wave system in three
  dimensions}, Acta Math. \textbf{219} (2017), no.~2, 213--402. \MR{3784694}

\bibitem{dias2003internal}
F.~Dias and J.-M. Vanden-Broeck, \emph{On internal fronts}, J. Fluid Mech.
  \textbf{479} (2003), 145--154. \MR{2011822}

\bibitem{di:handbook}
Fr\'{e}d\'{e}ric Dias and G\'{e}rard Iooss, \emph{Water-waves as a spatial
  dynamical system}, Handbook of mathematical fluid dynamics, {V}ol. {II},
  North-Holland, Amsterdam, 2003, pp.~443--499. \MR{1984157}

\bibitem{DiasKharif99}
Fr{\'e}d{\'e}ric Dias and Christian Kharif, \emph{Nonlinear gravity and
  capillary-gravity waves}, Annual review of fluid mechanics, {V}ol.\ 31, Annu.
  Rev. Fluid Mech., vol.~31, Annual Reviews, Palo Alto, CA, 1999, pp.~301--346.
  \MR{1670945 (99k:76024)}

\bibitem{DjordjevicRedekopp77}
V.~D. Djordjevi{\'c} and L.~G. Redekopp, \emph{On two-dimensional packets of
  capillary-gravity waves}, J. Fluid Mech. \textbf{79} (1977), no.~4, 703--714.
  \MR{0443555 (56 \#1924)}

\bibitem{dubreil}
M.~L. Dubreil-Jacotin, \emph{Sur la d\'etermination rigoureuse des ondes
  permanentes p\'eriodiques d'amplitude finie}, Journ. de Math. \textbf{13}
  (1934), 217--289.

\bibitem{dubreil1937theoremes}
\bysame, \emph{Sur les theoremes d'existence relatifs aux ondes
  permanentes periodiques a deux dimensions dans les liquides heterogenes}, J.
  Math. Pures Appl. \textbf{16} (1937), no.~9, 43--67.

\bibitem{DKSZ1996}
A.~I. Dyachenko, E.~A. Kuznetsov, M.~D. Spector, and V.~E. Zakharov,
  \emph{Analytical description of the free surface dynamics of an ideal fluid
  (canonical formalism and conformal mapping)}, Physics Letters A \textbf{221}
  (1996), no.~1, 73--79.

\bibitem{DLK2013}
S.~A. Dyachenko, P.~M. Lushnikov, and A.~O. Korotkevich, \emph{The complex
  singularity of a {S}tokes wave}, JETP Letters \textbf{98} (2013), no.~11,
  767--771.

\bibitem{DLK2016}
\bysame, \emph{Branch cuts of {S}tokes wave on deep water. {P}art {I}:
  {N}umerical solution and {P}ad\'e approximation}, Stud. Appl. Math.
  \textbf{137} (2016), no.~4, 419--472. \MR{3570991}

\bibitem{DH2019a}
Sergey~A. Dyachenko and Vera~Mikyoung Hur, \emph{Stokes waves with constant vorticity: {I}. {N}umerical
  computation}, Stud. Appl. Math. \textbf{142} (2019), no.~2, 162--189.
  \MR{3915685}

\bibitem{DH2019b}
\bysame, \emph{Stokes waves with constant vorticity: folds, gaps and fluid
  bubbles}, J. Fluid Mech. \textbf{878} (2019), 502--521. \MR{4010456}

\bibitem{DH2019c}
\bysame, \emph{Stokes waves in a constant
  vorticity flow}, Nonlinear Water Waves, Tutorials, Schools, and Workshops in
  the Mathematical Sciences, Birkh\"auser Basel, 2019, p.~18pp.

\bibitem{DHS2021}
Sergey~A. Dyachenko, Vera~Mikyoung Hur, and Denis~A. Silantyev, \emph{Almost
  extreme waves}, Preprint (2021).

\bibitem{EhrnstromGroves18}
Mats Ehrnstr\"{o}m and Mark~D. Groves, \emph{Small-amplitude fully localised
  solitary waves for the full-dispersion {K}adomtsev-{P}etviashvili equation},
  Nonlinearity \textbf{31} (2018), no.~12, 5351--5384. \MR{3871845}

\bibitem{ehrnstrom2019smooth}
Mats Ehrnstr{\"o}m, Samuel Walsh, and Chongchun Zeng, \emph{Smooth stationary
  water waves with exponentially localized vorticity}, to appear in J. Eur.
  Math. Soc. (2020).

\bibitem{EncisoPeralta-Salas16}
Alberto Enciso and Daniel Peralta-Salas, \emph{Beltrami fields with a
  nonconstant proportionality factor are rare}, Arch. Rational Mech. Anal.
  \textbf{220} (2016), no.~1, 243--260. \MR{3458163}

\bibitem{EncisoPeralta-Salas18}
\bysame, \emph{Existence of knotted vortex structures in stationary solutions
  of the {E}uler equations}, European {C}ongress of {M}athematics, Eur. Math.
  Soc., Z\"{u}rich, 2018, pp.~133--153. \MR{3887764}

\bibitem{EncisoPeralta-SalasRomaniega20}
Alberto Enciso, Daniel Peralta-Salas, and {\'A}lvaro Romaniega, \emph{Beltrami
  fields exhibit knots and chaos almost surely}, 2020.

\bibitem{escher2020stratified}
Joachim Escher, Patrik Knopf, Christina Lienstromberg, and Bogdan-Vasile
  Matioc, \emph{Stratified periodic water waves with singular density
  gradients}, Ann. Mat. Pura Appl. (4) \textbf{199} (2020), no.~5, 1923--1959.
  \MR{4142857}

\bibitem{escher2011stratified}
Joachim Escher, Anca-Voichita Matioc, and Bogdan-Vasile Matioc, \emph{On
  stratified steady periodic water waves with linear density distribution and
  stagnation points}, J. Differential Equations \textbf{251} (2011), no.~10,
  2932--2949. \MR{2831719 (2012h:35276)}

\bibitem{fs:nophase}
Gr\'{e}gory Faye and Arnd Scheel, \emph{Center manifolds without a phase
  space}, Trans. Amer. Math. Soc. \textbf{370} (2018), no.~8, 5843--5885.
  \MR{3803149}

\bibitem{filippov1960vortex}
I.~G. Filippov, \emph{Solution of the problem of the motion of a vortex under
  the surface of a fluid, for {F}roude numbers near unity}, J. Appl. Math.
  Mech. \textbf{24} (1960), 698--716. \MR{0128208}

\bibitem{filippov1961motion}
\bysame, \emph{On the motion of a vortex below the surface of a liquid}, J.
  Appl. Math. Mech. \textbf{25} (1961), 357--365. \MR{0151069}

\bibitem{fh:existence}
K.~O. Friedrichs and D.~H. Hyers, \emph{The existence of solitary waves}, Comm.
  Pure Appl. Math. \textbf{7} (1954), 517--550. \MR{0065317 (16,413f)}

\bibitem{Fuchs52}
R.~A. Fuchs, \emph{On the theory of short-crested oscillatory waves}, U.S.
  Natl. Bur. Stand. Circ. \textbf{521} (1952), 187--200.

\bibitem{gallay2011interaction}
Thierry Gallay, \emph{Interaction of vortices in weakly viscous planar flows},
  Arch. Ration. Mech. Anal. \textbf{200} (2011), no.~2, 445--490. \MR{2787587}

\bibitem{GermainMasmoudiShatah12}
P.~Germain, N.~Masmoudi, and J.~Shatah, \emph{Global solutions for the gravity
  water waves equation in dimension 3}, Ann. of Math. (2) \textbf{175} (2012),
  no.~2, 691--754. \MR{2993751}

\bibitem{glass2018point}
Olivier Glass, Alexandre Munnier, and Franck Sueur, \emph{Point vortex dynamics
  as zero-radius limit of the motion of a rigid body in an irrotational fluid},
  Inventiones mathematicae (2018).

\bibitem{Grad71}
Harold Grad, \emph{Mathematical problems arising in plasma physics}, Actes du
  {C}ongr\`es {I}nternational des {M}ath\'{e}maticiens ({N}ice, 1970), {T}ome
  3, 1971, pp.~105--113. \MR{0421280}

\bibitem{Groves07}
M.~D. Groves, \emph{Three-dimensional travelling gravity-capillary water
  waves}, GAMM-Mitt. \textbf{30} (2007), no.~1, 8--43. \MR{MR2324393
  (2008f:76024)}

\bibitem{GrovesHaragus03}
M.~D. Groves and M.~Haragus, \emph{A bifurcation theory for three-dimensional
  oblique travelling gravity-capillary water waves}, J. Nonlinear Sci.
  \textbf{13} (2003), no.~4, 397--447. \MR{MR1992382 (2004e:76020)}

\bibitem{GrovesHorn20}
M.~D. Groves and J.~Horn, \emph{A variational formulation for steady surface
  water waves on a beltrami flow}, Proceedings of the Royal Society A:
  Mathematical, Physical and Engineering Sciences \textbf{476} (2020),
  no.~2234, 20190495.

\bibitem{GrovesMielke01}
M.~D. Groves and A.~Mielke, \emph{A spatial dynamics approach to
  three-dimensional gravity-capillary steady water waves}, Proc. Roy. Soc.
  Edinburgh Sect. A \textbf{131} (2001), no.~1, 83--136. \MR{MR1820296
  (2002k:76016)}

\bibitem{GrovesSun08}
M.~D. Groves and S.-M. Sun, \emph{Fully localised solitary-wave solutions of
  the three-dimensional gravity-capillary water-wave problem}, Arch. Ration.
  Mech. Anal. \textbf{188} (2008), no.~1, 1--91. \MR{MR2379653 (2009c:76009)}

\bibitem{GrovesSunWahlen16a}
M.~D. Groves, S.~M. Sun, and E.~Wahl{\'e}n, \emph{A dimension-breaking
  phenomenon for water waves with weak surface tension}, Arch.\ Ration.\ Mech.\
  Anal.\/ \textbf{220} (2016), no.~2, 747--807.

\bibitem{groves:survey}
Mark~D. Groves, \emph{Steady water waves}, J. Nonlinear Math. Phys. \textbf{11}
  (2004), no.~4, 435--460. \MR{2097656 (2006a:76014)}

\bibitem{Haragus15}
Mariana Haragus, \emph{Transverse dynamics of two-dimensional gravity-capillary
  periodic water waves}, J. Dynam. Differential Equations \textbf{27} (2015),
  no.~3-4, 683--703. \MR{3435127}

\bibitem{hi:book}
Mariana Haragus and G\'{e}rard Iooss, \emph{Local bifurcations, center
  manifolds, and normal forms in infinite-dimensional dynamical systems},
  Universitext, Springer-Verlag London, Ltd., London; EDP Sciences, Les Ulis,
  2011. \MR{2759609}

\bibitem{haziot2021stratified}
Susanna~V. Haziot, \emph{Stratified large-amplitude steady periodic water waves
  with critical layers}, Comm. Math. Phys. \textbf{381} (2021), no.~2,
  765--797. \MR{4207457}

\bibitem{helmholtz1858integrale}
H.~Helmholtz, \emph{{\"U}ber {I}ntegrale der hydrodynamischen {G}leichungen,
  welche den {W}irbelbewegungen entsprechen}, J. Reine Angew. Math. \textbf{55}
  (1858), 25--55. \MR{1579057}

\bibitem{henry2014global}
David Henry and Anca-Vocihita Matioc, \emph{Global bifurcation of
  capillary--gravity-stratified water waves}, Proc. Roy. Soc. Edinburgh Sect. A
  \textbf{144} (2014), no.~4, 775--786. \MR{3233756}

\bibitem{HillLoydTurner21}
Dan~J. Hill, David J.~B. Lloyd, and Matthew~R. Turner, \emph{Localised radial
  patterns on the free surface of a ferrofluid}, 2021.

\bibitem{Hur2006}
Vera~Mikyoung Hur, \emph{Global bifurcation theory of deep-water waves with
  vorticity}, SIAM J. Math. Anal. \textbf{37} (2006), no.~5, 1482--1521.
  \MR{2215274}

\bibitem{Hur2011}
\bysame, \emph{Stokes waves with vorticity}, J. Anal. Math. \textbf{113}
  (2011), 331--386. \MR{2788362}

\bibitem{HVB2020}
Vera~Mikyoung Hur and Jean-Marc Vanden-Broeck, \emph{A new application of
  {C}rapper's exact solution to waves in constant vorticity flows}, Eur. J.
  Mech. B Fluids \textbf{83} (2020), 0--4. \MR{4102022}

\bibitem{HW2020}
Vera~Mikyoung Hur and Miles~H. Wheeler, \emph{Exact free surfaces in constant
  vorticity flows}, J. Fluid Mech. \textbf{896} (2020), R1, 10. \MR{4111654}

\bibitem{HW2021}
\bysame, \emph{Overhanging and touching waves in constant vorticity flows},
  Preprint (2021).


\bibitem{Ioossp1}
G. Iooss and P.~I. Plotnikov
\emph{Existence of multimodal standing gravity waves.}
J. Math. Fluid Mech. \textbf{7} (3) (2005), 349--364.

\bibitem{Ioossp2}
\bysame, \emph{Multimodal standing gravity waves: a completely resonant system}. J. Math. Fluid Mech. \textbf{7} (1) (2005), 110--126.

\bibitem{IoossPlotnikov11}
\bysame, \emph{Asymmetrical three-dimensional
  travelling gravity waves}, Arch. Ration. Mech. Anal. \textbf{200} (2011),
  no.~3, 789--880. \MR{2796133 (2012f:35430)}

\bibitem{IoossPlotnikov09}
G{\'e}rard Iooss and Pavel~I. Plotnikov, \emph{Small divisor problem in the
  theory of three-dimensional water gravity waves}, Mem. Amer. Math. Soc.
  \textbf{200} (2009), no.~940, viii+128. \MR{2529006 (2011e:76021)}

  \bibitem{ioosspt}
G. Iooss, P. I.  Plotnikov and J. F.  Toland.
\emph{Standing waves on an infinitely deep perfect fluid under gravity}.
Arch. Ration. Mech. Anal. \textbf{177} (3) (2005),  367--478.

\bibitem{james1997small}
Guillaume James, \emph{Small amplitude steady internal waves in stratified
  fluids}, Ann. Univ. Ferrara Sez. VII (N.S.) \textbf{43} (1997), 65--119
  (1998). \MR{1686749}

\bibitem{james2001internal}
\bysame, \emph{Internal travelling waves in the limit of a discontinuously
  stratified fluid}, Arch. Rational Mech. Anal. \textbf{160} (2001), no.~1,
  41--90. \MR{1864121 (2002h:76031)}

\bibitem{jtw:capgrav}
Mathew~A. Johnson, Tien Truong, and Miles~H. Wheeler, \emph{Solitary waves in a
  Whitham equation with small surface tension}, arXiv preprint arXiv:2103.02675
  (2021).

\bibitem{KadomtsevPetviashvili70}
B.B. Kadomtsev and V.I. Petviashvili, \emph{{On the stability of solitary waves
  in weakly dispersing media}}, Sov. Phys., Dokl. \textbf{15} (1970), 539--541
  (English. Russian original).

  \bibitem{keady}
  G. Keady and J. Norbury, \emph{{ On the existence theory for irrotational water waves}},
  Math. Proc. Cambridge Philos. Soc. 83 (1978), no. 1, 137-157.

\bibitem{KimAkylas05}
B. Kim and T.~R. Akylas, \emph{On gravity-capillary lumps}, J. Fluid Mech.
  \textbf{540} (2005), 337--351. \MR{MR2263130 (2007e:76027)}

\bibitem{kirchgassner1993structure}
K.~Kirchg{{\"a}}ssner and K.~Lankers, \emph{Structure of permanent waves in
  density-stratified media}, Meccanica \textbf{28} (1993), 269--276.

\bibitem{kirchgassner1982wavesolutions}
K. Kirchg{{\"a}}ssner, \emph{Wave-solutions of reversible systems and
  applications}, J. Differential Equations \textbf{45} (1982), no.~1, 113--127.
  \MR{662490 (83j:35063)}

\bibitem{kirchhoff1876vorlesungen}
G.~R. Kirchhoff, \emph{Vorlesungen {\"u}ber mathematische physik:
  mechanik}, vol.~1, Teubner, 1876.

\bibitem{kirrmann1991reduction}
P. Kirrmann, \emph{Reduktion nichtlinearer elliptischer systeme in
  zylindergebeiten unter verwendung von optimaler regularit{\"a}t in
  h{\"o}lder-r{\"a}umen}, Ph.D. thesis, Universit{\"a}t Stuttgart, 1991.

\bibitem{KleinSaut12}
C.~Klein and J.-C. Saut, \emph{Numerical study of blow up and stability of
  solutions of generalized {K}adomtsev-{P}etviashvili equations}, J. Nonlinear
  Sci. \textbf{22} (2012), no.~5, 763--811. \MR{2982052}

\bibitem{klw:nosubcrit}
V. Kozlov, E. Lokharu, and M.~H. Wheeler, \emph{Nonexistence of
  {S}ubcritical {S}olitary {W}aves}, Arch. Ration. Mech. Anal. \textbf{241}
  (2021), no.~1, 535--552. \MR{4271966}

   \bibitem{krasovskii}
Yu. P. Krasovskii.
 \emph{On the theory of steady waves of finite amplitude}. U.S.S.R. Comput. Math. Math. Phys. \textbf{1} (1961), 996--1018.

 \bibitem{kuznetsov}
N. G. Kuznetsov. A tale of two Nekrasov equations. arXiv:2009.01754 (2020).

\bibitem{Kwong89}
Man~Kam Kwong, \emph{Uniqueness of positive solutions of {$\Delta u-u+u^p=0$}
  in {${\bf R}^n$}}, Arch. Rational Mech. Anal. \textbf{105} (1989), no.~3,
  243--266. \MR{969899 (90d:35015)}

\bibitem{Laget1997interfacial}
O.~Laget and F.~Dias, \emph{Numerical computation of capillary-gravity
  interfacial solitary waves}, J. Fluid Mech. \textbf{349} (1997), 221--251.
  \MR{1480072}

\bibitem{lamb}
Horace Lamb, \emph{Hydrodynamics}, sixth ed., Cambridge Mathematical Library,
  Cambridge University Press, Cambridge, 1993, With a foreword by R. A.
  Caflisch [Russel E. Caflisch]. \MR{1317348 (96f:76001)}

\bibitem{lamb2000three}
Kevin~G. Lamb, \emph{Conjugate flows for a three-layer fluid}, Phys. Fluids
  \textbf{12} (2000), no.~9, 2169--2185. \MR{1780436 (2001d:76045)}

\bibitem{lankers1997fast}
K.~Lankers and G.~Friesecke, \emph{Fast, large-amplitude solitary waves in the
  2d euler equations for stratified fluids}, Nonlinear Anal. \textbf{29}
  (1997), no.~9, 1061--1078.

\bibitem{lannes}
D. Lannes. \emph{Well-posedness of the water-waves equations}. J. Amer. Math. Soc. \textbf{18} (3) (2005), 605--654.

\bibitem{lannes2013stability}
\bysame, \emph{A stability criterion for two-fluid interfaces and
  applications}, Arch. Ration. Mech. Anal. \textbf{208} (2013), no.~2,
  481--567. \MR{3035985}

\bibitem{Lannes13}
\bysame, \emph{The water waves problem. mathematical analysis and asymptotics},
  Mathematical Surveys and Monographs, vol. 188, American Mathematical Society,
  Providence, RI, 2013. \MR{3060183}

\bibitem{LannesFloating}
\bysame, \emph{On the dynamics of floating structures}, Ann. PDE 3 (2017), no. 1, Paper No. 11, 81 pp.

\bibitem{lavrentiev}
M.~A. Lavrentiev, \emph{I. {O}n the theory of long waves. {II}. {A}
  contribution to the theory of long waves}, Amer. Math. Soc. Translation
  \textbf{1954} (1954), no.~102, 53. \MR{0061952 (15,906a)}

\bibitem{le2018transmission}
Hung Le, \emph{Elliptic equations with transmission and {W}entzell boundary
  conditions and an application to steady water waves in the presence of wind},
  Discrete Contin. Dyn. Syst. \textbf{38} (2018), no.~7, 3357--3385.
  \MR{3809086}

\bibitem{le2019existence}
\bysame, \emph{On the existence and instability of solitary water waves with a
  finite dipole}, SIAM J. Math. Anal. \textbf{51} (2019), no.~5, 4074--4104.
  \MR{4019191}

\bibitem{lc}
T.~Levi-Civita, \emph{Determinazione rigorosa delle onde irrotazionali
  periodiche in acqua profonda}, Rend. Accad. Lincei \textbf{33} (1924),
  141--150.

\bibitem{li1998dirichlet}
Yanyan Li and Louis Nirenberg, \emph{The {D}irichlet problem for singularly
  perturbed elliptic equations}, Comm. Pure Appl. Math. \textbf{51} (1998),
  no.~11-12, 1445--1490. \MR{1639159}

\bibitem{lin1988large}
C.-S. Lin, W.-M. Ni, and I.~Takagi, \emph{Large amplitude stationary solutions
  to a chemotaxis system}, J. Differential Equations \textbf{72} (1988), no.~1,
  1--27. \MR{929196}

\bibitem{lin:instability}
Zhiwu Lin, \emph{On linear instability of 2{D} solitary water waves}, Int.
  Math. Res. Not. IMRN (2009), no.~7, 1247--1303. \MR{2495304 (2010f:35310)}

\bibitem{LiuWei19}
Yong Liu and Juncheng Wei, \emph{Nondegeneracy, {M}orse index and orbital
  stability of the {KP}-{I} lump solution}, Arch. Ration. Mech. Anal.
  \textbf{234} (2019), no.~3, 1335--1389. \MR{4011698}

\bibitem{LiuWang97}
Yue Liu and Xiao-Ping Wang, \emph{Nonlinear stability of solitary waves of a
  generalized {K}adomtsev-{P}etviashvili equation}, Comm. Math. Phys.
  \textbf{183} (1997), no.~2, 253--266. \MR{1461958}

\bibitem{LokharuSethWahlen20}
E.~Lokharu, D.~S.~Seth, and E.~Wahl{\'e}n, \emph{An existence theory for
  small-amplitude doubly periodic water waves with vorticity}, Arch. Rational
  Mech. Anal. \textbf{238} (2020), no.~2, 607--637.

\bibitem{LokharuWahlen19}
E.~Lokharu and E.~Wahl\'{e}n, \emph{A variational principle for
  three-dimensional water waves over {B}eltrami flows}, Nonlinear Anal.
  \textbf{184} (2019), 193--209. \MR{3921103}

\bibitem{lokharu:bounds}
Evgeniy Lokharu, \emph{On bounds for steady waves with negative vorticity}, J.
  Math. Fluid Mech. \textbf{23} (2021), no.~2, Paper No. 37, 11. \MR{4228659}

\bibitem{lombardi:book}
Eric Lombardi, \emph{Oscillatory integrals and phenomena beyond all algebraic
  orders}, Lecture Notes in Mathematics, vol. 1741, Springer-Verlag, Berlin,
  2000, With applications to homoclinic orbits in reversible systems.
  \MR{1770093}

\bibitem{LHF1977a}
M.~S. Longuet-Higgins and M.~J.~H. Fox, \emph{Theory of the almost-highest
  wave: the inner solution}, J. Fluid Mech. \textbf{80} (1977), no.~4,
  721--741. \MR{452143}

\bibitem{LHF1977b}
\bysame 
\emph{Theory of the almost-highest
  wave. II. Matching and analytic extension}, J. Fluid Mech. \textbf{85} (1978), no.~4,
  769--786. \MR{502858}

   \bibitem{LongH:78}
M.~S. Longuet-Higgins, \emph{Some new relations between
{S}tokes's coefficients
  in the theory of gravity waves}. J. Inst. Maths. Applics. \textbf{22} (1978),  261--273.

\bibitem{Lortz70}
Dietrich Lortz, \emph{\"{U}ber die {E}xistenz toroidaler magnetohydrostatischer
  {G}leichgewichte ohne {R}otationstransformation}, Z. Angew. Math. Phys.
  \textbf{21} (1970), 196--211. \MR{261852}

\bibitem{Lushnikov2016}
Pavel~M. Lushnikov, \emph{Structure and location of branch point singularities
  for {S}tokes waves on deep water}, J. Fluid Mech. \textbf{800} (2016),
  557--594. \MR{3521330}

\bibitem{LDS2017}
Pavel~M. Lushnikov, Sergey~A. Dyachenko, and Denis~A. Silantyev, \emph{New
  conformal mapping for adaptive resolving of the complex singularities of
  {S}tokes wave}, Proc. A. \textbf{473} (2017), no.~2202, 20170198, 19.
  \MR{3672665}

\bibitem{makarenko1992bore}
N.~I. Makarenko, \emph{Smooth bore in a two-layer fluid}, Free boundary
  problems in continuum mechanics ({N}ovosibirsk, 1991), Internat. Ser. Numer.
  Math., vol. 106, Birkh{\"a}user, Basel, 1992, pp.~195--204. \MR{1229538
  (94d:76097)}

\bibitem{makarenko1999conjugate}
\bysame, \emph{Conjugate flows and smooth bores in a weakly stratified fluid},
  Prikl. Mekh. Tekhn. Fiz. \textbf{40} (1999), no.~2, 69--78. \MR{1711960}

\bibitem{ManakovZakharovBordagItsMatveev77}
S.V. Manakov, V.E. Zakharov, L.A. Bordag, A.R. Its, and V.B. Matveev,
  \emph{Two-dimensional solitons of the kadomtsev-petviashvili equation and
  their interaction.}, Physics Letters A \textbf{63} (1977), no.~3, 205--206.

\bibitem{marchioro1993vortices}
Carlo Marchioro and Mario Pulvirenti, \emph{Vortices and localization in
  {E}uler flows}, Comm. Math. Phys. \textbf{154} (1993), no.~1, 49--61.
  \MR{1220946}

\bibitem{marchioro1994book}
\bysame, \emph{Mathematical theory of incompressible nonviscous fluids},
  Applied Mathematical Sciences, vol.~96, Springer-Verlag, New York, 1994.
  \MR{1245492 (94k:76001)}

\bibitem{Martin18}
Calin~Iulian Martin, \emph{Non-existence of time-dependent three-dimensional
  gravity water flows with constant non-zero vorticity}, Physics of Fluids
  \textbf{30} (2018), no.~10, 107102.

\bibitem{mcleod:froude}
J.~B. McLeod, \emph{The {F}roude number for solitary waves}, Proc. Roy. Soc.
  Edinburgh Sect. A \textbf{97} (1984), 193--197. \MR{751191 (85g:76011)}

\bibitem{McLeod}
\bysame, \emph{The {S}tokes and {K}rasovskii conjectures for the wave of
  greatest height}, Stud. Appl. Math. \textbf{98} (1997), no.~4, 311--333.
  \MR{1446239}

\bibitem{mho:saddlecenter}
A.~Mielke, P.~Holmes, and O.~O'Reilly, \emph{Cascades of homoclinic orbits to,
  and chaos near, a {H}amiltonian saddle-center}, J. Dynam. Differential
  Equations \textbf{4} (1992), no.~1, 95--126. \MR{1150399}

\bibitem{mielke:reduction}
Alexander Mielke, \emph{Reduction of quasilinear elliptic equations in
  cylindrical domains with applications}, Math. Methods Appl. Sci. \textbf{10}
  (1988), no.~1, 51--66. \MR{929221 (89d:35063)}

\bibitem{mielke1995homoclinic}
\bysame, \emph{Homoclinic and heteroclinic solutions in two-phase flow},
  Proceedings of the {IUTAM}/{ISIMM} {S}ymposium on {S}tructure and {D}ynamics
  of {N}onlinear {W}aves in {F}luids ({H}annover, 1994), Adv. Ser. Nonlinear
  Dynam., vol.~7, World Sci. Publ., River Edge, NJ, 1995, pp.~353--362.
  \MR{1685879}

\bibitem{miles1957windwaves1}
John~W. Miles, \emph{On the generation of surface waves by shear flows}, J.
  Fluid Mech. \textbf{3} (1957), 185--204. \MR{0091692 (19,1004e)}

\bibitem{miles1959windwaves2}
\bysame, \emph{On the generation of surfaces waves by shear flows. {II}}, J.
  Fluid Mech. \textbf{6} (1959), 568--582. \MR{0108951 (21 \#7663)}

\bibitem{miles:survey}
\bysame, \emph{Solitary waves}, Annual review of fluid mechanics, {V}ol. 12,
  Annual Reviews, Palo Alto, Calif., 1980, pp.~11--43. \MR{565388 (82e:76018)}

\bibitem{nakayama2020breathers}
K.~Nakayama and K.~G. Lamb, \emph{Breathers in a three-layer fluid}, J. Fluid
  Mech. \textbf{903} (2020), A40, 20. \MR{4157455}

\bibitem{Nekrasov}
A.~I. Nekrasov, \emph{On steady waves}, Izv. Ivanovo-Voznesensk. Politekhn.
  In-ta \textbf{3} (1921).

\bibitem{NgSt}
H. Q. Nguyen and W. A. Strauss, \emph{Proof of modulational instability of Stokes waves in deep water}, arXiv:2007.05018 (2020).

\bibitem{ni1991shape}
Wei-Ming Ni and Izumi Takagi, \emph{On the shape of least-energy solutions to a
  semilinear {N}eumann problem}, Comm. Pure Appl. Math. \textbf{44} (1991),
  no.~7, 819--851. \MR{1115095}



\bibitem{ni1998location}
Wei-Ming Ni, Izumi Takagi, and Juncheng Wei, \emph{On the location and profile
  of spike-layer solutions to a singularly perturbed semilinear {D}irichlet
  problem: intermediate solutions}, Duke Math. J. \textbf{94} (1998), no.~3,
  597--618. \MR{1639546}

\bibitem{ni1995location}
Wei-Ming Ni and Juncheng Wei, \emph{On the location and profile of spike-layer
  solutions to singularly perturbed semilinear {D}irichlet problems}, Comm.
  Pure Appl. Math. \textbf{48} (1995), no.~7, 731--768. \MR{1342381}

  \bibitem{nicholls}
David~P. Nicholls and Fernando Reitich, \emph{A new approach to analyticity of Dirichlet-Neumann operators}.
Proc. Roy. Soc. Edinburgh Sect. A \textbf{131} (6) (2001), 1411--1433.

\bibitem{NichollsReitich06}
\bysame, \emph{Stable, high-order computation of
  traveling water waves in three dimensions}, Eur. J. Mech. B Fluids
  \textbf{25} (2006), no.~4, 406--424. \MR{2241379}

\bibitem{nilsson2017internal}
Dag~Viktor Nilsson, \emph{Internal gravity-capillary solitary waves in finite
  depth}, Math. Methods Appl. Sci. \textbf{40} (2017), no.~4, 1053--1080.
  \MR{3610717}

\bibitem{Nilsson19b}
\bysame, \emph{Three-dimensional internal gravity-capillary waves in finite
  depth}, Math. Methods Appl. Sci. \textbf{42} (2019), no.~12, 4113--4145.
  \MR{3978643}

\bibitem{PapanicolaouSulemSulemWang94}
G.~C. Papanicolaou, C.~Sulem, P.-L. Sulem, and X.~P. Wang, \emph{The focusing
  singularity of the {D}avey-{S}tewartson equations for gravity-capillary
  surface waves}, Phys. D \textbf{72} (1994), no.~1-2, 61--86. \MR{1270855}

\bibitem{ps:stability}
Robert~L. Pego and Shu-Ming Sun, \emph{Asymptotic linear stability of solitary
  water waves}, Arch. Ration. Mech. Anal. \textbf{222} (2016), no.~3,
  1161--1216. \MR{3544325}

\bibitem{PelinovskyStepanyants93}
D.~Pelinovsky and Y.~Stepanyants, \emph{New multisoliton solutions of the
  {K}adomtsev-{P}etviashvili equation}, JETP Lett. \textbf{57} (1993), 24--28.

\bibitem{Plotnikov1982}
P.~I. Plotnikov, \emph{Proof of the {S}tokes conjecture in the theory of
  surface waves}, Dokl. Akad. Nauk SSSR \textbf{269} (1983), no.~1, 80--83.
  \MR{699357}

\bibitem{plotnikov:turning}
\bysame, \emph{Nonuniqueness of solutions of a problem on solitary waves, and
  bifurcations of critical points of smooth functionals}, Izv. Akad. Nauk SSSR
  Ser. Mat. \textbf{55} (1991), no.~2, 339--366. \MR{1133302 (93d:35132)}

\bibitem{PT2004}
P.~I. Plotnikov and J.~F. Toland, \emph{Convexity of {S}tokes waves of extreme
  form}, Arch. Ration. Mech. Anal. \textbf{171} (2004), no.~3, 349--416.
  \MR{2038344}

\bibitem{pommerenke:book}
Ch. Pommerenke, \emph{Boundary behaviour of conformal maps}, Grundlehren der
  Mathematischen Wissenschaften [Fundamental Principles of Mathematical
  Sciences], vol. 299, Springer-Verlag, Berlin, 1992. \MR{1217706}

\bibitem{ParauVandenBroeckCooker05a}
E.~I. P\u{a}r\u{a}u, J.-M. Vanden-Broeck, and M.~J. Cooker, \emph{Nonlinear
  three-dimensional gravity-capillary solitary waves}, J. Fluid Mech.
  \textbf{536} (2005), 99--105. \MR{MR2263898 (2007e:76032)}

\bibitem{rabinowitz:global}
Paul~H. Rabinowitz, \emph{Some global results for nonlinear eigenvalue
  problems}, J. Functional Analysis \textbf{7} (1971), 487--513. \MR{0301587
  (46 \#745)}

\bibitem{ReederShinbrot81}
John Reeder and Marvin Shinbrot, \emph{Three-dimensional, nonlinear wave
  interaction in water of constant depth}, Nonlinear Anal. \textbf{5} (1981),
  no.~3, 303--323. \MR{607813 (82i:35168)}

\bibitem{rouhi1993hamiltonian}
Ali Rouhi and Jon Wright, \emph{Hamiltonian formulation for the motion of
  vortices in the presence of a free surface for ideal flow}, Physical Review E
  \textbf{48} (1993), no.~3, 1850.

\bibitem{grue2002solitary}
Per-Olav Rus{\aa}s and John Grue, \emph{Solitary waves and conjugate flows in a
  three-layer fluid}, Eur. J. Mech. B Fluids \textbf{21} (2002), no.~2,
  185--206. \MR{1933937 (2003h:76026)}

\bibitem{russell}
J.~Scott Russell, \emph{Report on waves}, 14th meeting of the British
  Association for the Advancement of Science, vol. 311--390, 1844.

\bibitem{schn-wayne}
G. Schneider and C. E. Wayne, \emph{The long-wave limit for the water wave problem. I. The case of zero surface tension} Comm. Pure Appl. Math. 53 (2000), 1475--1535.

\bibitem{SethVarholmWahlen21}
Douglas~Svensson Seth, Kristoffer Varholm, and Erik Wahl{\'e}n, \emph{Symmetric
  doubly periodic gravity-capillary waves with small vorticity}, Preprint.

     \bibitem{sharg}
E. Shargorodsky.
\emph{An estimate for the Morse index of a Stokes wave}. Arch. Ration. Mech. Anal. \textbf{209} (1) (2013), 41--59.

\bibitem{ShargTol:01}
  E. Shargorodsky and J. F. Toland. \emph{A Riemann-Hilbert
    problem and the {B}ernoulli boundary condition in the
    variational theory of {S}tokes waves}. Annales Inst.
  H. Poincar\'e \textbf{20} (1) (2003), 37--52.

\bibitem{ShargTol:03}
  \bysame,  \emph{Riemann-Hilbert
    theory for problems with vanishing coefficients that arise
    in nonlinear hydrodynamics}. Jour. Funct. Anal.
  \textbf{197}  (2003), 283--300.

\bibitem{bernoulli}
\bysame, \emph{Bernoulli free-boundary problems}. Mem. Amer. Math. Soc. \textbf{196} No 914, (2008). ISBN: 978-0-8218-4189-1.

\bibitem{shatah2013travelling}
Jalal Shatah, Samuel Walsh, and Chongchun Zeng, \emph{Travelling water waves
  with compactly supported vorticity}, Nonlinearity \textbf{26} (2013), no.~6,
  1529--1564. \MR{3053431}

\bibitem{shatah2011interface}
Jalal Shatah and Chongchun Zeng, \emph{Local well-posedness for fluid interface
  problems}, Arch. Rational Mech. Anal. \textbf{199} (2011), no.~2, 653--705.
  \MR{2763036 (2012c:35336)}

\bibitem{SS1985}
J.~A. Simmen and P.~G. Saffman, \emph{Steady deep-water waves on a linear shear
  current}, Stud. Appl. Math. \textbf{73} (1985), no.~1, 35--57. \MR{797557}

\bibitem{sinambela2020large}
Daniel Sinambela, \emph{Large-amplitude solitary waves in two-layer density
  stratified water}, SIAM J. Math. Anal., \textbf{53} (2021), no.~4, 4812--4864.

\bibitem{Sretenskii53}
L.~Sretenskii, \emph{Spatial probelm of determination of steady waves of finite
  amplitude ({R}ussian)}, Dokl. Akad. Nauk SSSR (N.S.) \textbf{89} (1953),
  25--28.

\bibitem{starr}
Victor~P. Starr, \emph{Momentum and energy integrals for gravity waves of
  finite height}, J. Mar. Res. \textbf{6} (1947), 175--193.

\bibitem{Stokes1847}
George~Gabriel Stokes, \emph{On the theory of oscillatory waves}, Trans. Camb.
  Philos. Soc. \textbf{8} (1847), 441--455.

\bibitem{Stokes1880}
\bysame, \emph{Considerations relative to the greatest height of oscillatory
  irrotational waves which can be propagated without change of form},
  Mathematical and Physical Papers, vol.~1, Cambridge University Press, 1880,
  pp.~314--326.

\bibitem{Strauss77}
Walter~A. Strauss, \emph{Existence of solitary waves in higher dimensions},
  Comm. Math. Phys. \textbf{55} (1977), no.~2, 149--162. \MR{0454365 (56
  \#12616)}

\bibitem{struik}
D.~J. Struik, \emph{D\'etermination rigoureuse des ondes irrotationelles
  p\'eriodiques dans un canal \`a profondeur finie}, Math. Ann. \textbf{95}
  (1926), no.~1, 595--634. \MR{1512296}

\bibitem{su2020longtime}
Qingtang Su, \emph{Long time behavior of 2{D} water waves with point vortices},
  Comm. Math. Phys. \textbf{380} (2020), no.~3, 1173--1266. \MR{4179726}

\bibitem{sun2002solitary}
S.~M. Sun, \emph{Solitary internal waves in continuously stratified fluids of
  great depth}, Phys. D \textbf{166} (2002), no.~1-2, 76--103. \MR{1913747
  (2003d:76029)}

\bibitem{Sun93}
T.~Sun, \emph{Three-dimensional steady water waves generated by partially
  localized pressure disturbances}, SIAM J. Math. Anal. \textbf{24} (1993),
  no.~5, 1153--1178.

\bibitem{sp:steep}
A.~F. Teles~da Silva and D.~H. Peregrine, \emph{Steep, steady surface waves on
  water of finite depth with constant vorticity}, J. Fluid Mech. \textbf{195}
  (1988), 281--302. \MR{985439 (90a:76061)}

\bibitem{TdSP1988}
\bysame, \emph{Steep, steady surface waves on water of finite depth with
  constant vorticity}, J. Fluid Mech. \textbf{195} (1988), 281--302.
  \MR{985439}

\bibitem{terkrkorov1958vortex}
A.~M. Ter-Krikorov, \emph{Exact solution of the problem of the motion of a
  vortex under the surface of a liquid}, Izv. Akad. Nauk SSSR Ser. Mat.
  \textbf{22} (1958), 177--200. \MR{0108145}

\bibitem{ter1963theorie}
\bysame, \emph{Th{\'e}orie exacte des ondes longues stationnaires dans
  un liquide h{\'e}t{\'e}rog{\`e}ne}, J. M{\'e}canique \textbf{2} (1963),
  351--376. \MR{MR0160400 (28 \#3613)}

   \bibitem{stokes}
J. F.  Toland,
\emph{Stokes Waves}.Topological Methods in Nonlinear Analysis,
Journal of the Juliusz Schauder Center.
\textbf{7} (1) (1996), 1–48 \& \textbf{8} (2), 413--414.

\bibitem{Toland:00}
\bysame, \emph{Regularity of Stokes waves in Hardy
spaces and in
  spaces of distributions}. Jour. Math. Pure et Appl. \textbf{79} (9) (2000), 901-917.

\bibitem{To:01}
\bysame, \emph{On a pseudo-differential equation for Stokes
waves.} Arch. Rational mech. Anal. \textbf{162} (2002), 179-189.

\bibitem{fareast}
J. F. Toland and  G. Iooss.
\emph{Riemann-Hilbert and variational structure for standing waves}.
Far East J. Appl. Math. \textbf{15} (3) (2004), 459--488.

\bibitem{jww:whitham}
T. Truong, E. Wahl{\'e}n, and M.~H. Wheeler, \emph{Global bifurcation of
  solitary waves for the Whitham equation}, arXiv preprint arXiv:2009.04713
  (2020).

\bibitem{turner1981internal}
R.~E.~L. Turner, \emph{Internal waves in fluids with rapidly varying density},
  Ann. Scuola Norm. Sup. Pisa Cl. Sci. (4) \textbf{8} (1981), no.~4, 513--573.
  \MR{MR656000 (83j:76027)}

\bibitem{turner1984variational}
\bysame, \emph{A variational approach to surface solitary waves}, J.
  Differential Equations \textbf{55} (1984), no.~3, 401--438.

\bibitem{varholm2016solitary}
Kristoffer Varholm, \emph{Solitary gravity-capillary water waves with point
  vortices}, Discrete Contin. Dyn. Syst. \textbf{36} (2016), no.~7, 3927--3959.
  \MR{3485858}

\bibitem{varholm:global}
\bysame, \emph{Global bifurcation of waves with multiple critical layers}, SIAM
  J. Math. Anal. \textbf{52} (2020), no.~5, 5066--5089. \MR{4164492}

\bibitem{varholm2020stability}
Kristoffer Varholm, Erik Wahl\'{e}n, and Samuel Walsh, \emph{On the stability
  of solitary water waves with a point vortex}, Comm. Pure Appl. Math.
  \textbf{73} (2020), no.~12, 2634--2684. \MR{4164269}

\bibitem{vonkarman1940engineer}
Theodore von K\'{a}rm\'{a}n, \emph{The engineer grapples with non-linear
  problems}, Bull. Amer. Math. Soc. \textbf{46} (1940), 615--683. \MR{3131}

\bibitem{wahlen:crit}
Erik Wahl{\'e}n, \emph{Steady water waves with a critical layer}, J.
  Differential Equations \textbf{246} (2009), no.~6, 2468--2483. \MR{2498849
  (2010i:76026)}

\bibitem{Wahlen14}
\bysame, \emph{Non-existence of three-dimensional travelling water waves with
  constant non-zero vorticity}, Journal of Fluid Mechanics \textbf{746} (2014).

\bibitem{walsh:stratified}
Samuel Walsh, \emph{Stratified steady periodic water waves}, SIAM J. Math.
  Anal. \textbf{41} (2009), no.~3, 1054--1105. \MR{2529956 (2011a:35430)}

\bibitem{walsh2014local}
\bysame, \emph{Steady stratified periodic gravity waves with surface
  tension {I}: {L}ocal bifurcation}, Discrete Contin. Dyn. Syst. Ser. A
  \textbf{34} (2014), no.~8, 3287--3315.

\bibitem{walsh2014global}
\bysame, \emph{Steady stratified periodic gravity waves with surface tension
  {II}: {G}lobal bifurcation}, Discrete Contin. Dyn. Syst. Ser. A \textbf{34}
  (2014), no.~8, 3241--3285.

\bibitem{walsh2013wind}
Samuel Walsh, Oliver B{{\"u}}hler, and Jalal Shatah, \emph{Steady water waves
  in the presence of wind}, SIAM J. Math. Anal. \textbf{45} (2013), no.~4,
  2182--2227. \MR{3073648}

\bibitem{wang2017small}
Ling-Jun Wang, \emph{Small-amplitude solitary and generalized solitary
  traveling waves in a gravity two-layer fluid with vorticity}, Nonlinear Anal.
  \textbf{150} (2017), 159--193. \MR{3584939}

\bibitem{WangAblowitzSegur94}
X.~P. Wang, M.~J. Ablowitz, and H.~Segur, \emph{Wave collapse and instability
  of solitary waves of a generalized {K}adomtsev-{P}etviashvili equation},
  Phys. D \textbf{78} (1994), no.~3-4, 241--265. \MR{1302410}

\bibitem{Wang19}
Xuecheng Wang, \emph{Global solution for the 3{D} gravity water waves system
  above a flat bottom}, Adv. Math. \textbf{346} (2019), 805--886. \MR{3914181}

\bibitem{WangMilewski12}
Zhan Wang and Paul~A. Milewski, \emph{Dynamics of gravity-capillary solitary
  waves in deep water}, J. Fluid Mech. \textbf{708} (2012), 480--501.
  \MR{2975453}

\bibitem{Weinstein83}
Michael~I. Weinstein, \emph{Nonlinear {S}chr\"{o}dinger equations and sharp
  interpolation estimates}, Comm. Math. Phys. \textbf{87} (1982/83), no.~4,
  567--576. \MR{691044}

\bibitem{Weinstein85}
\bysame, \emph{Modulational stability of ground states of nonlinear
  {S}chr\"{o}dinger equations}, SIAM J. Math. Anal. \textbf{16} (1985), no.~3,
  472--491. \MR{783974}

\bibitem{wheeler:solitary}
Miles~H. Wheeler, \emph{Large-amplitude solitary water waves with vorticity},
  SIAM J. Math. Anal. \textbf{45} (2013), no.~5, 2937--2994. \MR{3106477}

\bibitem{wheeler:pressure}
\bysame, \emph{Solitary water waves of large amplitude generated by surface
  pressure}, Arch. Ration. Mech. Anal. \textbf{218} (2015), no.~2, 1131--1187.
  \MR{3375547}

\bibitem{Wheeler18a}
\bysame, \emph{Integral and asymptotic properties of solitary waves in deep
  water}, Comm. Pure Appl. Math. \textbf{71} (2018), 1941--1956.

\bibitem{wheeler2019stratified}
\bysame, \emph{On stratified water waves with critical layers and {C}oriolis
  forces}, Discrete Contin. Dyn. Syst. \textbf{39} (2019), no.~8, 4747--4770.
  \MR{3986308}

\bibitem{Wu97}
Sijue Wu, \emph{Well-posedness in Sobolev spaces of the full water wave problem in 2-D}, Invent. Math. 130 (1997), no. 1, 39--72.

\bibitem{Wu11}
\bysame, \emph{Global wellposedness of the 3-{D} full water wave problem},
  Invent. Math. \textbf{184} (2011), no.~1, 125--220. \MR{2782254}

\bibitem{yanowitch1962gravity}
M.~Yanowitch, \emph{Gravity waves in a heterogeneous incompressible fluid},
  Comm. Pure Appl. Math. \textbf{15} (1962), 45--61. \MR{MR0165792 (29 \#3072)}

  \bibitem{zak}
 V. E. Zakharov. Stability of periodic waves of finite amplitude on the surface of deep fluid. Jour. Appl. Mech. Tech. Phys. 2 (1968), 190-194.

\end{thebibliography}

\end{document}